\documentclass[12pt]{article}
\usepackage{amssymb,amsmath}
\usepackage[dvips]{graphicx}
\usepackage{psfrag}

\numberwithin{equation}{section}

\newcommand{\tw}{\rm t w}

\newcommand{\bfone}{{\bf 1}}

\newcommand{\colim}{{\rm c o l i m}}

\newcommand{\Spec}{\rm S p e c}

\newcommand{\otimesO}{\otimes_{{\cal O}_X}}

\newcommand{\Coder}{\rm C o d e r}

\newcommand{\Vb}{V^{\bullet}}

\newcommand{\Cbu}{C^{\bullet}}
\newcommand{\Cbd}{C_{\bullet}}

\newcommand{\nCbu}{C^{\bullet}_{{\rm norm}}}
\newcommand{\nCbd}{C_{\bullet}^{{\rm norm}}}

\newcommand{\ccCb}{\check{{\mathcal C}}^{\bullet}}
\newcommand{\ccC}{\check{{\mathcal C}}}

\newcommand{\cpa}{\check{\partial}}

\newcommand{\Cyl}{{\rm C y l}\,}
\newcommand{\tCyl}{\widetilde{\rm C y l}\,}
\newcommand{\Conf}{{\rm Conf}\,}
\newcommand{\Der}{{\rm Der}\,}
\newcommand{\Hom}{{\rm Hom}\,}

\newcommand{\cHom}{\mathcal Hom\,}

\newcommand{\Alg}{{\rm Alg}\,}
\newcommand{\Coalg}{{\rm Coalg}\,}
\newcommand{\restr}{{\rm restr}\,}

\newcommand{\sgn}{{\rm s g n}}

\newcommand{\lan}{\langle}
\newcommand{\ran}{\rangle}

\newcommand{\Lie}{{\bf Lie}}
\newcommand{\comm}{{\bf comm}}
\newcommand{\ass}{{\bf assoc}}
\newcommand{\coLie}{{\bf coLie}}
\newcommand{\cocomm}{{\bf cocomm}}

\newcommand{\Ger}{{\bf e_2}}
\newcommand{\calc}{{\bf calc}}
\newcommand{\pcalc}{{\bf pcalc}}

\newcommand{\KS}{{\bf K S}}
\newcommand{\bB}{{\bf B}}

\newcommand{\Ho}{{\rm Ho}\,}

\newcommand{\Hoger}{\Ho({\bf e_2})}
\newcommand{\Hocalc}{\Ho({\bf calc})}

\newcommand{\bul}{{\bullet}}
\newcommand{\Omb}{{\Omega^{-\bullet}}}
\newcommand{\OmbX}{{\Omega_X^{-\bullet}}}

\newcommand{\al}{{\alpha}}
\newcommand{\la}{{\lambda}}
\newcommand{\io}{{\iota}}

\newcommand{\Om}{{\Omega}}
\newcommand{\Si}{{\Sigma}}
\newcommand{\si}{{\sigma}}
\newcommand{\ga}{{\gamma}}
\newcommand{\vf}{{\varphi}}
\newcommand{\ve}{{\varepsilon}}
\newcommand{\ka}{{\kappa}}
\newcommand{\vr}{{\varrho}}
\newcommand{\G}{{\Gamma}}

\newcommand{\ml}{{\mathfrak{l}}}

\newcommand{\ma}{{\mathfrak{a}}}

\newcommand{\mc}{{\mathfrak{c}}}

\newcommand{\md}{{\mathfrak{d}}}

\newcommand{\tm}{{\tilde{m}}}

\newcommand{\bS}{{\bf S}}

\newcommand{\bs}{{\bf s}\,}
\newcommand{\bx}{{\bf x}}
\newcommand{\bq}{{\bf q}}

\newcommand{\pa}{{\partial}}

\newcommand{\cY}{{\cal Y}}
\newcommand{\cF}{{\cal F}}
\newcommand{\cD}{{\cal D}}

\newcommand{\cG}{{\cal G}}
\newcommand{\cH}{{\cal H}}
\newcommand{\cP}{{\cal P}}
\newcommand{\cQ}{{\cal Q}}

\newcommand{\cB}{{\cal B}}
\newcommand{\cM}{{\cal M}}
\newcommand{\cL}{{\cal L}}
\newcommand{\cI}{{\cal I}}
\newcommand{\cC}{{\cal C}}

\newcommand{\cR}{{\cal R}}
\newcommand{\cV}{{\cal V}}
\newcommand{\cW}{{\cal W}}
\newcommand{\cO}{{\cal O}}
\newcommand{\ocC}{\overline{\cal C}}
\newcommand{\ocO}{\overline{\cal O}}

\newcommand{\bbC}{{\mathbb C}}
\newcommand{\bbR}{{\mathbb R}}
\newcommand{\bbZ}{{\mathbb Z}}
\newcommand{\bbK}{{\mathbb K}}
\newcommand{\bbH}{{\mathbb H}}

\newcommand{\bbL}{{\mathbb L}}
\newcommand{\bbG}{{\mathbb G}}

\newcommand{\bbF}{{\mathbb F}}
\newcommand{\bbT}{{\mathbb T}}

\newcommand{\La}{{\Lambda}}

\newcommand{\te}{\theta}

\newcommand{\de}{{\delta}}
\newcommand{\hde}{\widehat{\delta}}
\newcommand{\D}{{\Delta}}
\newcommand{\Ups}{{\Upsilon}}

\newcommand{\tcV}{\widetilde{{\cal V}}}
\newcommand{\tcW}{\widetilde{{\cal W}}}
\newcommand{\tnu}{{\widetilde{\nu}}}
\newcommand{\tDe}{{\widetilde{\Delta}}}

\newcommand{\tW}{{\widetilde{W}}}

\newcommand{\tf}{{\widetilde{f}}}

\newcommand{\tga}{{\widetilde{\gamma}}}

\date{}
\newtheorem{defi}{Definition}
\newtheorem{pred}{Proposition}

\newtheorem{teo}{Theorem}
\newtheorem{cor}{Corollary}

\title{Formality of the homotopy calculus algebra
of Hochschild (co)chains}

\author{
Vasiliy Dolgushev,
Dmitry Tamarkin, and
Boris Tsygan}

\begin{document}

\maketitle

\begin{center}
{\it To Mikhail Olshanetsky on the occasion of his 70th birthday.}
\end{center}

\begin{abstract}
The Kontsevich-Soibelman
solution of the cyclic version of
Deligne's conjecture and the formality of
the operad of little discs on a cylinder
provide us with a natural homotopy
calculus structure on the pair
$(\Cbu(A), \Cbd(A))$
``Hochschild cochains $+$ Hochschild chains''
of an associative algebra $A$.
We show that for an arbitrary smooth algebraic variety
$X$ over a field $\bbK$ of characteristic zero
the sheaf $(\Cbu(\cO_X), \Cbd(\cO_X))$
of homotopy calculi is formal.
This result was announced in paper
\cite{TT1} by the second and the third author.
\end{abstract}

\tableofcontents

\section{Introduction}
The standard Cartan calculus on polyvector fields and exterior forms
can be naturally extended to the Hochschild cohomology $HH^{\bul}(A,A)$ and the
Hochschild homology $HH_{\bul}(A,A)$ of an arbitrary  associative algebra $A$
\cite{DGT}, \cite{Rinehart}. This calculus is induced by simple operations
on Hochschild (co)chains,  and the identities of this
algebraic structure hold for these operations up to homotopy.

The Kontsevich-Soibelman
proof of the cyclic version of
Deligne's conjecture \cite{K-Soi1} and the formality of
the operad of little discs on a cylinder\footnote{See
Proposition 11.3.3 on page 50 in \cite{K-Soi1}.}
imply that this nice collection of the operations on
the pair $(\nCbu(A), \nCbd(A))$
``(normalized) Hochschild cochains $+$
(normalized) Hochschild chains''
can be extended to
an $\infty$- or homotopy calculus structure.

This homotopy calculus structure on the pair
$(\nCbu(A), \nCbd(A))$ is a natural generalization
of the homotopy Gerstenhaber algebra structure on
the cochains $\nCbu(A)$\,.
In paper \cite{BLT} we proved the formality of this
homotopy Gerstenhaber algebra on $\nCbu(A)$ for an
arbitrary regular commutative algebra $A$ over a field
$\bbK$ of characteristic zero.
In this paper we extend this result to
the homotopy calculus algebra on the pair
$(\nCbu(A), \nCbd(A))$\,.

As well as in \cite{BLT}
we also consider the situation when the algebra $A$
is replaced by the structure sheaf $\cO_X$ of a smooth algebraic
variety $X$ over the field $\bbK$\,. More precisely, we consider the
homotopy calculus algebra
on the pair $(\nCbu(\cO_X), \nCbd(\cO_X))$
where $\nCbu(\cO_X)$ and $\nCbd(\cO_X)$ is, respectively, the sheaf of
(normalized) Hochschild cochains and the sheaf of (normalized) Hochschild
chains of $\cO_X$\,. In this paper we show that the sheaf of homotopy calculi
$(\nCbu(\cO_X), \nCbd(\cO_X))$ is formal.

If $A$ is an associative
algebra (with unit),
the pair $(\nCbu(A), \nCbd(A))$
is also equipped with an algebraic structure
defined by a degree $-1$ Lie bracket
on $\nCbu(A)$\,, a degree $-1$ Lie module structure on
$\nCbd(A)$ over $\nCbu(A)$\,, and Connes' operator
on $\nCbd(A)$ which is compatible with the Lie module
structure. In the paper we refer to such algebra
structures as $\La\Lie^+_{\de}$-algebra.
(See Definition \ref{lie+de}.)

In paper \cite{Tsygan} the third author
conjectured that if $A$ is a regular
commutative algebra then this $\La\Lie^+_{\de}$-algebra
structure on $(\nCbu(A), \nCbd(A))$ is formal.
This conjecture was proved in \cite{W} (at least in the
case $\bbR \subset \bbK$) by Willwacher who used
the constructions of B. Shoikhet \cite{Sh} and the
first author \cite{FTHC}.

In general $Ho(\La\Lie^+_{\de})$-part of the homotopy calculus structure
on $(\nCbu(A), \nCbd(A))$ derived from \cite{K-Soi1} may not coincide
with the $\La\Lie^+_{\de}$-algebra on the pair $(\nCbu(A), \nCbd(A))$\,.
However, we show that this homotopy calculus algebra
on $(\nCbu(A), \nCbd(A))$ is quasi-isomorphic to another
homotopy calculus algebra on $(\nCbu(A), \nCbd(A))$ whose
$Ho(\La\Lie^+_{\de})$-part is the ordinary $\La\Lie^+_{\de}$-algebra
given by the above Lie bracket on $\nCbu(A)$, the Lie algebra
module on $\nCbd(A)$ over $\nCbu(A)$ and Connes' operator on
$\nCbd(A)$\,. In this sense, the formality of the homotopy calculus
algebra on $(\nCbu(A), \nCbd(A))$ is a generalization of
Willwacher's cyclic formality theorem \cite{W}.

The organization of the paper is as follows.
In Section 2 we fix the notation and recall required results
about (co)operads and (co)algebras. Section 3 is devoted to
$\infty$- or homotopy versions for the algebras over
the operads $\calc$\,, $\Ger$\,, and $\Lie^+_{\de}$\,.
In Section 4 we recall the Kontsevich-Soibelman
operad and the operad $\Cyl$ of little discs on a cylinder.
We show that the homology operad $H_{-\bul}(\Cyl, \bbK)$ of $\Cyl$
with the reversed grading is the operad of calculi.
Finally we recall required results from \cite{K-Soi1} and prove
a useful property of the Kontsevich-Soibelman operad.
Section 5 is devoted to properties of the homotopy calculus
algebra on the pair $(\nCbu(A), \nCbd(A))$\,.
In Section 6 we formulate and prove the main result of this
paper. (See Theorem \ref{main} on page \pageref{main}.)
In the concluding section we discussion applications and
generalizations of Theorem \ref{main}. We also discuss recent
articles related to our main result.

~\\
{\bf Acknowledgment.}
A bigger part of this work was done when V.D.
was a Boas Assistant Professor of Mathematics Department
at Northwestern University. During these two years V.D. benefited
from working at Northwestern so much that he
feels as if he finished one more graduate school.
V.D. cordially thanks Mathematics
Department at Northwestern University for this time.
The results of this work were presented at the conference
Poisson 2008 in Lausanne. We would like to thank
the participants of this conference for questions and
useful comments.
V.D. would like to thank Pavel Snopok for showing him
a very convenient drawing program ``Inkscape''.
D.T. and B.T. are supported by
NSF grants. The work of V.D. is partially supported by the Grant
for Support of Scientific
Schools NSh-8065.2006.2.

\section{Preliminaries}
\subsection{(Co)operads and (co)algebras}
Most of the notation and conventions
for (co)operads and their (co)algebras
are borrowed from \cite{BLT}.

Depending on a context our underlying
symmetric monoidal category is either
the category of graded vector spaces, or
the category of chain complexes, or the category
of compactly generated topological spaces,
or the category of sets.
By suspension $\bs \cV$ of a graded vector space
(or a chain complex) $\cV$
we mean $\ve \otimes \cV$, where $\ve$ is a
one-dimensional vector space placed in degree
$+1$\,. For a vector $v\in \cV$ we denote by $|v|$
its degree. The symmetric group of
permutations of $n$ letters is denoted by $S_n$\,.
The underlying field $\bbK$
has characteristic zero.

For an operad $\cO$
we denote by $\Alg_{\cO}$ the category of algebras
over the operad $\cO$. Dually, for a cooperad
$\cC$ we denote by $\Coalg_{\cC}$ the category
of nilpotent\footnote{For the definition of
{\it nilpotent} coalgebra see section $2.4.1$ in
\cite{Hinich}.}
coalgebras over the cooperad $\cC$\,.
By {\it corestriction} we mean the canonical
map
\begin{equation}
\label{corest}
\rho_{\cV} : \bbF_{\cC}(\cV) \to \cV
\end{equation}
from the free coalgebra
$\bbF_{\cC}(\cV)$ to the vector space of its
cogenerators $\cV$\,. We often omit the subscript
in the notation $\rho_{\cV}$ for the corestriction.

For a polynomial functor $\cP$ we denote by
$\bbT(\cP)$ (resp. $\bbT^*(\cP)$) the free operad
(resp. the free cooperad) (co)generated by $\cP$\,.
The notation $\bullet$ is reserved for the
monoidal product of the polynomial functors.
Thus, if $\cP$ and $\cQ$ are polynomial
functors then
\begin{equation}
\label{bullet}
\cP \bullet \cQ(n) =
\bigoplus_{k_1+ \dots + k_m = n}
\cP(m) \otimes_{S_m} ( \cQ(k_1) \otimes
\dots \otimes \cQ(k_m) )\,.
\end{equation}
This formula can be easily generalized to
the colored polynomial functors.

By ``suspension'' of a (co)operad $\cO$ of graded vector spaces
(or chain complexes)
we mean the (co)operad $\La(\cO)$
whose $m$-th vector space is
\begin{equation}
\label{susp-op}
\La(\cO)(m) = \bs^{1-m} \cO(m) \otimes \sgn_{m}\,,
\end{equation}
where $\sgn_{m}$ is the sign representation of
the symmetric group $S_m$\,.

For a commutative algebra $\cB$ and a $\cB$-module
$\cV$ we denote by $S_{\cB}(\cV)$ the symmetric algebra
of $\cV$ over $\cB$\,.
$S_{\cB}^m(\cV)$ stands for the $m$-th component
of this algebra. If $\cB=\bbK$ then $\cB$ is omitted
from the notation. The abbreviation ``DGLA'' stands for
differential graded Lie algebra.

We denote by $\ast$ the polynomial functor
\begin{equation}
\label{ast}
\ast(n) =
\begin{cases}
\bbK\,, \qquad {\rm if} ~~ n = 1 \,, \\
0 \,, \qquad {\rm otherwise}\,.
\end{cases}
\end{equation}
This functor carries the unique structure of
the operad (resp. the cooperad) such that
$\ast$ is the initial (resp. the terminal) object
in the category of operads (resp. cooperads) of
graded vector spaces or chain
complexes. There is an obvious generalization
of $\ast$ (\ref{ast}) to the category of sets and to the
category of topological spaces. However, we will need $\ast$
only for linear (co)operads, i.e. the (co)operads in
the category of graded vector spaces or the
category of chain complexes.

All the linear operads (resp. linear cooperads), we consider, are equipped with
an augmentation (resp. coaugmentation). In other words,
for every operad $\cO$ we will have a chosen morphism
of operads:
\begin{equation}
\label{aug}
\tau\, :\, \cO \to \ast\,.
\end{equation}
Dually for every cooperad $\cC$ we will have a chosen
morphism of cooperads
\begin{equation}
\label{coaug}
\ka\, :\, \ast \to \cC\,.
\end{equation}

We are going to deal with $2$-colored (co)operads.
Throughout the paper we label the two
colors of all $2$-colored (co)operads by $\mc$ and $\ma$\,.
For example,
the notation $\Lie^+$ is reserved for the
$2$-colored operad which governs the pairs ``Lie algebra
$\cV$ and a Lie algebra module $\cW$ over $\cV$\,.''
Vectors of the Lie algebra $\cV$ are colored by $\mc$
and vectors of the module $\cW$ are colored by $\ma$\,.

For a linear $2$-colored operad $\cO$ we will denote by
$\cO^{\mc}(n,k)$ (resp. $\cO^{\ma}(n,k)$) the vector space of
operations producing a vector with the color $\mc$ (resp. $\ma$)
from $n$ vectors with the color $\mc$ and
$k$ vectors with the color $\ma$\,. We use the same notation
for the linear $2$-colored cooperads and
for topological $2$-colored operads.

The polynomial functor $\ast$ (\ref{ast}) has the obvious generalization
to the category of linear $2$-colored (co)operads:
\begin{equation}
\label{ast-color}
\begin{array}{c}
\ast^{\mc}(n,k) =
\begin{cases}
\bbK\,, \qquad {\rm if} ~~ (n,k)= (1,0) \,, \\
0 \,, \qquad {\rm otherwise}\,.
\end{cases} \\[0.5cm]
\ast^{\ma}(n,k) =
\begin{cases}
\bbK\,, \qquad {\rm if} ~~ (n,k)= (0,1)\,, \\
0 \,, \qquad {\rm otherwise}\,.
\end{cases}
\end{array}
\end{equation}

For a linear operad $\cO$ we denote by
$Bar(\cO)$ its bar construction. Dually, for a linear
cooperad $\cC$ we denote by
$Cobar(\cC)$ its cobar construction.

We recall that,
as a cooperad
of graded vector spaces, $Bar(\cO)$ is freely
generated by the polynomial functor $\bs^{-1}\ocO$\,,
where $\ocO$ is the kernel of the augmentation
(\ref{aug}). Dually, as an operad
of graded vector spaces, $Cobar(\cC)$ is freely
generated by the polynomial functor $\bs\ocC$\,,
where $\ocC$ is cokernel of the coaugmentation
(\ref{coaug}). The differential $\pa^{Bar}$ on the
operad $Bar(\cO)$ is defined using the multiplication
of the operad $\cO$ and the differential $\pa^{Cobar}$
on the cooperad $Cobar(\cC)$ is defined using the
comultiplication of the cooperad $\cC$\,. See Chapter 3
in \cite{Fresse} or Section 2 in \cite{GJ} for details.

For a quadratic operad $\cO$ there is a natural
sub-cooperad $\cO^{\vee}$ of $Bar(\cO)$ which satisfies the
property:
$$
\pa^{Bar} \Big|_{\cO^{\vee}} = 0\,.
$$
The details of the construction of $\cO^{\vee}$ can
be found in Section 5.2 in \cite{Fresse}.
Following \cite{GK} we call $\cO^{\vee}$ the Koszul dual
cooperad of $\cO$\,.

For a linear operad $\cO$ (resp. linear cooperad $\cC$) and
a vector space $\cV$ we denote by
$\bbF_{\cO}(\cV)$ (resp. by $\bbF_{\cC}(\cV)$) the
free algebra (resp. free coalgebra) over the
operad $\cO$ (resp. cooperad $\cC$).
For a linear $2$-colored (co)operad $\cO$ the
functor\footnote{$\bbF_{\cO}$ is called the Schur functor.}
$\bbF_{\cO}$
splits according to the colors $\mc$ and $\ma$ as
$$
\bbF_{\cO}(\cV, \cW) =
\bbF_{\cO}(\cV, \cW)_{\mc} \oplus
\bbF_{\cO}(\cV, \cW)_{\ma}\,,
$$
where
$$
\bbF_{\cO}(\cV, \cW)_{\mc} =
\bigoplus_{n,k}\cO^{\mc}(n,k)
\otimes_{S_n \times S_k} \cV^{\otimes \, n} \otimes \cW^{\otimes \, k}\,,
$$
and
$$
\bbF_{\cO}(\cV, \cW)_{\ma} =
\bigoplus_{n,k}\cO^{\ma}(n,k)
\otimes_{S_n \times S_k} \cV^{\otimes \, n} \otimes \cW^{\otimes \, k}\,.
$$

We need to recall some facts about
algebras over the operad $Cobar(\cC)$
for a coaugmented cooperad $\cC$\,.

Since $Cobar(\cC)$ is freely
generated by the suspension $\bs\ocC$ of the
cokernel $\ocC$
of the coaugmentation (\ref{coaug})
a $Cobar(\cC)$-algebra
structure on a chain complex $\cV$ is uniquely
determined by the restriction of the
multiplication map
$$
\mu: \bbF_{Cobar(\cC)}(\cV) \to \cV
$$
to the subspace
$$
\bbF_{\bs\ocC}(\cV) \subset \bbF_{Cobar(\cC)}(\cV)\,.
$$
In other words, a $Cobar(\cC)$-algebra
structure on $\cV$ is uniquely determined by
a degree $1$ map from $\bbF(\ocC)(\cV)$
to $\cV$\,.

It turns out that the maps from $\bbF(\ocC)(\cV)$
to $\cV$ have a elegant description in terms
of coderivations of the free coalgebra
$\bbF_{\cC}(\cV)$\,. To recall this description
we introduce the Lie subalgebra
\begin{equation}
\label{Coder-prime}
\Coder'(\bbF_{\cC}(\cV)) =
\{Q \in \Coder(\bbF_{\cC}(\cV)) ~|~
Q\Big|_{\cV}  = 0 \}\,,
\end{equation}
where $\cV$ is considered as a subspace
of $\cC(1)\otimes \cV$ via the coaugmentation
(\ref{coaug})\,. In other words, the elements
of $\Coder'(\bbF_{\cC}(\cV))$ are coderivations
of the free coalgebra
$\bbF_{\cC}(\cV)$ which can be factored through
the projection
$$
\bbF_{\cC}(\cV) \to \bbF_{\ocC}(\cV)\,.
$$
It is not hard to see that the subspace
(\ref{Coder-prime}) is closed under the commutator
and the differentials coming from $\cC$ and $\cV$\,.
Thus the graded vector space $\Coder'(\bbF_{\cC}(\cV))$ is
in fact a DGLA.

Let us recall from \cite{GJ} the
following proposition
\begin{pred}[Proposition 2.14 \cite{GJ}]
\label{coder-cofree}
For a coaugmented cooperad $\cC$
the composition with the corestriction
(\ref{corest}) $\rho_{\cV}: \bbF_{\cC}(\cV) \to \cV$
induces an isomorphism of graded vector spaces
\begin{equation}
\label{der-free}
Coder'(\bbF_{\cC}(\cV)) \cong
\Hom(\bbF_{\ocC}(\cV), \cV)\,,
\end{equation}
where, as above, $\ocC$ is the cokernel of the
coaugmentation (\ref{coaug}) of $\cC$\,.
\end{pred}
Due to this proposition a $Cobar(\cC)$-algebra
structure on a chain complex $\cV$ is uniquely
determined by a degree $1$ coderivation
\begin{equation}
\label{deriv-Q}
Q \in Coder'(\bbF_{\cC}(\cV))\,.
\end{equation}

According to Proposition 2.15 from \cite{GJ}
the compatibility of the $Cobar(\cC)$-algebra
structure on $\cV$
with the total differential on $Cobar(\cC)$
and the differential on $\cV$ is equivalent
to the Maurer-Cartan equation for the corresponding
derivation (\ref{deriv-Q}):
\begin{equation}
\label{MC-deriv-Q}
[\pa^{\cC} + \pa^{\cV},Q] + \frac{1}{2}[Q, Q] =0\,,
\end{equation}
where $\pa^{\cC}$ is the differential
on $\bbF_{\cC}(\cV)$ induced by the one on
the cooperad $\cC$ and $\pa^{\cV}$ comes from
that on $\cV$\,.

In other words,
\begin{pred}[Proposition 2.15, \cite{GJ}]
\label{GJ}
There is a natural
bijection between the Maurer-Cartan elements of the
DGLA $Coder'(\bbF_{\cC}(\cV))$ and
and
the $Cobar(\cC)$-algebra structures
on $\cV$.
\end{pred}

If we have a map
\begin{equation}
\label{mu}
\mu : \cC_1 \to \cC_2
\end{equation}
of coaugmented cooperads then
the corresponding map between
the cobar constructions
$$
Cobar(\mu) : Cobar(\cC_1) \to Cobar(\cC_2)
$$
allows us to pull
$Cobar(\cC_2)$-algebra structure on $\cV$
to a $Cobar(\cC_1)$-algebra on $\cV$\,.

We claim that
\begin{pred}
\label{how-to-pull}
If $Q_1$ is a Maurer-Cartan element of the DGLA
$Coder'(\bbF_{\cC_1}(\cV))$ corresponding to
a $Cobar(\cC_1)$-algebra structure on $\cV$
and $Q_2$ is a Maurer-Cartan element of the DGLA
$Coder'(\bbF_{\cC_2}(\cV))$ corresponding to
a $Cobar(\cC_2)$-algebra structure on $\cV$
then
\begin{equation}
\label{pull-MC}
\rho_{\cV} \circ Q_1 = \rho_{\cV} \circ Q_2 \circ \bbF(\mu)\,,
\end{equation}
where the map
$$
\bbF(\mu) : \bbF_{\ocC_1}(\cV) \to  \bbF_{\ocC_2}(\cV)
$$
is induced by (\ref{mu})\,.
\end{pred}
{\bf Proof.} Let
$$
\nu_2 :
\bbF_{Cobar(\cC_2)}(\cV) \to \cV
$$
be the $Cobar(\cC_2)$-algebra structure
on $\cV$\,. Then the $Cobar(\cC_1)$-algebra structure
on $\cV$
$$
\nu_1 :
\bbF_{Cobar(\cC_1)}(\cV) \to \cV
$$
is obtained by composing the map $\nu_2$
with the map
$$
\bbF(Cobar(\mu)) :
\bbF_{Cobar(\cC_1)}(\cV) \to
\bbF_{Cobar(\cC_2)}(\cV)\,.
$$
It is not hard to see that
the restriction of
$\nu_1$ to the subspace
$\bbF_{\bs\ocC_1}(\cV)$ coincides with the
composition of the maps
$$
\bbF_{\bs \ocC_1}(\cV) \stackrel{\bbF(\mu)}{\to}
\bbF_{\bs\ocC_2}(\cV)
$$
and
$$
\nu_2 \Big|_{\bbF_{\bs\ocC_2}(\cV)} :
\bbF_{\bs\ocC_2}(\cV) \to \cV\,.
$$
Thus the proposition follows from the
equation
$$
\nu_i \Big|_{\bbF_{\bs\ocC_i}(\cV)}
= \rho_{\cV}\circ Q_i \circ \si\,,
$$
where $\rho_{\cV}$ is the corestriction
(\ref{corest}) and $\si$ is the suspension
isomorphism
$$
\si: \bbF_{\bs\ocC_i}(\cV) \to \bbF_{\ocC_i}(\cV)\,,
$$
and $i=1,2$\,. $\Box$

We will freely use Propositions
\ref{coder-cofree}, \ref{GJ} and
\ref{how-to-pull} for
colored cooperads.

We remark that
all $2$-colored (co)operads, we consider, satisfy the following
property: {\it an argument with the color $\ma$ can enter
an operation at most once.
If an argument with this
color enters an operation then
the resulting color is also $\ma$\,.
Otherwise the resulting color is $\mc$\,. }
In other words, for every $n$
\begin{equation}
\label{P}
\begin{array}{cc}
\cO^{\mc}(n,k) = \cO^{\ma}(n,k) = \{ \bf 0\}\,, &
\forall ~~ k > 1\,,\\[0.3cm]
\cO^{\ma}(n,0) = \{ \bf 0\}\,, & \cO^{\mc}(n,1) = \{ \bf 0\}
\end{array}
\end{equation}
for the (co)operads of graded vector spaces or chain complexes and
\begin{equation}
\label{P-top}
\begin{array}{cc}
\cO^{\mc}(n,k) = \cO^{\ma}(n,k) = \emptyset\,, &
\forall ~~ k > 1\,,\\[0.3cm]
\cO^{\ma}(n,0) = \emptyset\,, & \cO^{\mc}(n,1) = \emptyset
\end{array}
\end{equation}
for the (co)operads of topological spaces or sets.

It is not hard to see that bar and cobar constructions
the (co)operads of graded vector spaces or chain complexes
preserve property (\ref{P}).

Let us recall that
\begin{defi}[M. Gerstenhaber, \cite{G}]
\label{G-alg}
A graded vector space $\cV$ is a Gerstenhaber algebra
if it is equipped with a commutative and associative
product $\wedge$ of degree $0$
and a Lie bracket $[\,,\,]$ of degree $-1$\,.
These operations have to be compatible in the sense of the
following Leibniz rule
\begin{equation}
\label{cup-Lie}
[a, b \wedge  c] =
[a, b] \wedge c +
(-1)^{(|a|+1)|b|} b \wedge [a, c]\,,
\end{equation}
where $a,b,c$ are homogeneous vectors of $\cV$\,.
\end{defi}

\begin{defi}
\label{precalc}
A precalculus is a pair of a Gerstenhaber algebra
$(\cV, \wedge, [,])$ and a graded vector
space $\cW$ together with
\begin{itemize}
\item a module structure
$i_{\bul} \,:\, \cV \otimes \cW
\mapsto \cW $ of the graded commutative algebra
$\cV$ on $\cW$\,,

\item an action $l_{\bul} \,:\, \bs^{-1} \cV\otimes \cW
\mapsto \cW$ of the Lie algebra
$\bs^{-1} \cV$ on
$\cW$ which are compatible in
the sense of the following equations
\begin{equation}
\label{l-i}
i_a l_b - (-1)^{|a|(|b|+1)}
l_b i_a = i_{[a,b]}\,,
\end{equation}
and
\begin{equation}
\label{l-cup}
l_{a\wedge b} = l_a i_{b} + (-1)^{|a|}i_a l_b \,.
\end{equation}
\end{itemize}
\end{defi}
Furthermore,
\begin{defi}
\label{calc}
A calculus is a precalculus
$(\cV,\cW, [,], \wedge, i_{\bul}, l_{\bul})$
with a degree $-1$ unary operation $\de$ on $\cW$
such that
\begin{equation}
\label{l-i-delta}
\delta \, i_{a} - (-1)^{|a|} i_{a} \, \delta =
l_a\,,
\end{equation}
and\footnote{Although $\de^2=0$\,, the operation $\de$ is
never considered as a part of the differential on $\cW$\,.}
\begin{equation}
\label{de-2}
\de^2 = 0\,.
\end{equation}
\end{defi}
We call $l$ and $i$ the Lie derivative and
the contraction, respectively.

We will use the following list of
(co)operads:
\begin{itemize}

\item $\Lie$ (resp. $\coLie$) is
the operad of Lie algebras (resp. the cooperad
of Lie coalgebras),

\item $\comm$ (resp. $\cocomm$)
is the operad of commutative (associative)
algebras (resp. the
operad of cocommutative coassociative
coalgebras),

\item $\Ger$ denotes the operad of Gerstenhaber
algebras, (see Definition \ref{G-alg}),

\item $\KS$ denotes the operad of M. Kontsevich
and Y. Soibelman. This
operad\footnote{In \cite{K-Soi1} this operad
is denoted by $P$.} is described in
sections 11.1, 11.2 and 11.3 of \cite{K-Soi1},

\item $\Lie^+$ (resp. $\coLie^+$) denotes the
$2$-colored operad of pairs ``Lie algebra $+$ its
module'' (resp. the $2$-colored cooperad of pairs
``Lie coalgebra $+$ its comodule''),

\item $\comm^+$ (resp. $\cocomm^+$) denotes the
$2$-colored operad of pairs ``commutative algebra $+$ its
module'' (resp. the $2$-colored cooperad of pairs
``cocommutative coalgebra $+$ its comodule''),

\item $\pcalc$ denotes the $2$-colored
operad of precalculi, (see Definition
\ref{precalc}),

\item $\calc$ denotes the $2$-colored
operad of calculi, (see Definition \ref{calc}),

\item $\ass$ is the
non-symmetric
operad of sets controlling unital monoids;
each set $\ass(n)$, $n\geq 0$, is a point.

\end{itemize}

It is not hard to show that for the vector space of the
free calculus algebra generated by the pair $(\cV, \cW)$ we have
\begin{equation}
\label{free-calc}
\bbF_{\calc}(\cV, \cW) \cong \bbF_{\comm^+}(\bbF_{\La\Lie^+}(\cV,
\cW\oplus \bs^{-1}\, \cW))\,.
\end{equation}
In other words, for the color components
we have the isomorphisms of graded vector spaces:
\begin{equation}
\label{free-calc-mc}
\bbF_{\calc}(\cV, \cW)_{\mc} \cong \bbF_{\comm}( \bbF_{\La\Lie}(\cV) )\,,
\end{equation}
and
\begin{equation}
\label{free-calc-ma}
\bbF_{\calc}(\cV, \cW)_{\ma} \cong
\bbF_{\comm^+}(\bbF_{\La\Lie}(\cV),   \bbF_{\La\Lie^+}(\cV,
\cW \oplus \bs^{-1}\, \cW)_{\ma} )_{\ma} \,.
\end{equation}

\subsection{Hochschild (co)chain complexes}
For an associative algebra $A$
$$
\Cbu(A) = \Hom (A^{\otimes \bul}, A)
$$
denotes the Hochschild cochain complex
and
$$
\Cbd(A)= A\otimes A^{\otimes (-\bul)}
$$
stands for the Hochschild chain
complex of $A$ with the reversed grading.

For the normalized versions of these
complexes we reserve the notation:
$$
\nCbu(A) = \{ P\in
Hom (A^{\otimes \bul}, A)\,\, |\,\,
P(\dots ,1, \dots) = 0 \}
$$
and
$$
\nCbd(A)= A\otimes (A/\bbK\, 1)^{\otimes (-\bul)}\,.
$$

\begin{itemize}

\item The notation $\pa^{Hoch}$ is reserved both for the Hochschild
coboundary operator on $\nCbu(A)$
and Hochschild boundary operator on $\nCbd(A)$
$$
(\pa^{Hoch} P)(a_0, a_1, \dots, a_k) = a_0 P(a_1, \dots, a_k) -
P(a_0 a_1, \dots, a_k) + P(a_0, a_1 a_2, a_3, \dots , a_k) -
\dots
$$
$$
+ (-1)^{k} P(a_0, \dots , a_{k-2}, a_{k-1} a_k) +
(-1)^{k+1} P(a_0, \dots , a_{k-2}, a_{k-1}) a_k
$$

$$
\pa^{Hoch}(a_0, a_1, \dots, a_m) =
(a_0 a_1, a_2, \dots, a_m) - (a_0, a_1 a_2, a_3, \dots,  a_m) +
\dots +
$$
$$
(-1)^{m-1}(a_0, \dots, a_{m-2}, a_{m-1} a_m) +
(-1)^m (a_m a_0, a_1, a_2, \dots, a_{m-1})\,,
$$
$$
a_i\in A\,, \qquad P \in C^k_{{\rm norm}}(A)\,.
$$

\item The notation $\cup$ is reserved for the cup-product
on $\nCbu(A)$
\begin{equation}
\label{cup}
P_1\cup P_2 (a_1, a_2, \dots, a_{k_1 + k_2}) =
P_1(a_1, \dots, a_{k_1}) P_2(a_{k_1+1}, \dots, a_{k_1+k_2})\,,
\end{equation}
$$
P_i \in C^{k_i}_{{\rm norm}}(A)\,.
$$

\item $[\, , \, ]_G$ stands for the Gerstenhaber bracket
on $\nCbu(A)$
$$
[Q_1, Q_2]_{G} =
$$
\begin{equation}
\label{Gerst}
\sum_{i=0}^{k_1}(-1)^{i k_2}
Q_1(a_0,\,\dots , Q_2 (a_i,\,\dots,a_{i+k_2}),\, \dots,
a_{k_1+k_2}) -
(-1)^{k_1 k_2} (1 \leftrightarrow 2)\,,
\end{equation}
$$
Q_i \in  C^{k_i+1}_{{\rm norm}}(A)\,.
$$

\item $I_P(c)$ is the contraction
of a Hochschild cochain $P \in C^k_{{\rm norm}}(A)$ with a
Hochschild chain $c=(a_0, a_1, \dots, a_m)$
\begin{equation}
\label{I-P}
I_P (a_0, a_1, \dots, a_m) =
\begin{cases}
(a_0 P(a_1, \dots, a_k), a_{k+1}, \dots , a_m)\,, \qquad {\rm if} ~~ m\ge k \,, \\
0 \,, \qquad {\rm otherwise}\,.
\end{cases}
\end{equation}

\item $L_Q(c)$ denotes the Lie derivative
of a Hochschild chain $c = (a_0, a_1, \dots, a_m)$ along a
Hochschild cochain $Q\in C^{k+1}_{{\rm norm}}(A)$
\begin{equation}
\label{L-Q}
L_{Q}(a_0, a_1, \dots, a_m)=
\sum_{i=0}^{m-k}(-1)^{ki} (a_0 , \dots,
Q(a_i, \dots, a_{i+k}), \dots, a_m) +
\end{equation}
$$
\sum_{j=m-k}^{m-1}(-1)^{m(j+1)} ( Q(a_{j+1}, \dots, a_m, a_0, \dots, a_{k+j-m}),
a_{k+j+1-m}, \dots , a_j )\,.
$$

\item $B: \nCbd(A) \to C^{{\rm norm}}_{\bul-1}(A)$
denotes Connes' operator
\begin{equation}
\label{B}
B(a_0, a_1, \dots, a_m)=
\sum_{i=0}^{m}(-1)^{m i}
(1, a_i , \dots, a_m, a_0, a_1, \dots, a_{i-1} )\,.
\end{equation}

\end{itemize}

The notation $HH^{\bul}(A)$ (resp. $HH_{\bul}(A)$) is used
for the Hochschild cohomology (resp. homology groups) of $A$
with coefficients in $A$
$$
HH^{\bul}(A) = H^{\bul}(\nCbu(A), \pa^{Hoch})\,,
$$
$$
HH_{\bul}(A) = H^{\bul}(\nCbd(A), \pa^{Hoch})\,.
$$

To describe algebraic structures on
pairs $(\nCbu(A), \nCbd(A))$ and
$(HH^{\bul}(A), HH_{\bul}(A))$ we use the
the language of operads.
Thus,
the Gerstenhaber bracket $[\,,\,]_G$ equips
the cochain complex $\nCbu(A)$ with an algebra
structure over the operad $\La\Lie$ and
the Lie derivative (\ref{L-Q}) equips the
pair $(\nCbu(A), \nCbd(A))$ with the algebra
structure over the operad $\La\Lie^+$\,.

In order to add Connes' operator (\ref{B}) into
this operadic picture we give one more
definition
\begin{defi}
\label{lie+de}
We say that the pair of graded vector
spaces $(\cV,\cW)$ is an algebra over the operad
$\Lie^+_{\de}$ if $\cV$ is a Lie algebra,
$\cW$ is a module over $\cV$ and $\cW$ is
equipped with a degree $-1$ unary operation
$\de$ satisfying the equations
\begin{equation}
\label{Lie-de-2}
\de^2 = 0\,,
\end{equation}
and
\begin{equation}
\label{l-delta}
[\delta, l_a] = 0\,, \qquad \forall ~~ a\in \cV\,,
\end{equation}
where $l$ is the action of $\cV$ on $\cW$\,.
\end{defi}
Adding Connes' operator $B$ into the picture
we may say that the pair $(\nCbu(A), \nCbd(A))$
is a $\La\Lie^+_{\de}$-algebra.

The operations (\ref{cup}), (\ref{Gerst}), (\ref{I-P}),
(\ref{L-Q}) and (\ref{B}) are closed with respect to the
(co)boundary operator $\pa^{Hoch}$\,.

According to \cite{G} the operations
$\cup$ (\ref{cup}) and $[\,,\,]_G$ (\ref{Gerst})
induce on $HH^{\bul}(A)$ the structure of a
Gerstenhaber algebra.
Furthermore, it is known \cite{DGT} that the operations
 (\ref{cup}), (\ref{Gerst}), (\ref{I-P}),
(\ref{L-Q}) and (\ref{B}) induce on the pair
$(HH^{\bul}(A), HH_{\bul}(A))$ the structure of
the calculus algebra.

\section{The operads $\Hocalc$\,, $\Hoger$\,,
and $\Ho(\Lie^+_{\de})$}
In this section we describe the homotopy versions
for the algebras over
the operads $\calc$\,, $\Ger$\,, and $\Lie^+_{\de}$\,.

\subsection{Description of the operads $\Hocalc$ and $\Hoger$}
To describe the homotopy version of $\calc$-algebras
we use the canonical cofibrant resolution
$Cobar(Bar(\calc))$\,. In other words, we set
\begin{equation}
\label{Calc}
\Hocalc = Cobar(Bar(\calc))\,.
\end{equation}

The cooperad $Bar(\calc)$ will be used throughout
the paper. For this reason we reserve a
short-hand notation
\begin{equation}
\label{bB}
\bB = Bar (\calc)
\end{equation}
for this cooperad.

Recall that, as a cooperad
of graded vector spaces, $\bB = Bar(\calc)$ is freely
generated by the polynomial functor $\bs^{-1}\overline{\calc}$\,,
where $\overline{\calc}$ is the kernel of the augmentation.

\begin{figure}
\begin{center}
\includegraphics[width=0.35\textwidth]{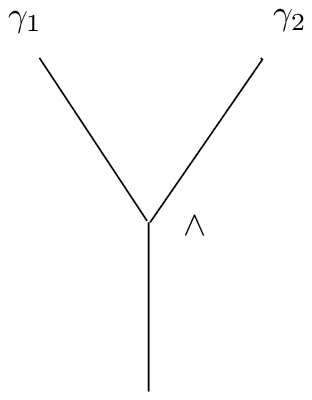}
\hfill
\includegraphics[width=0.35\textwidth]{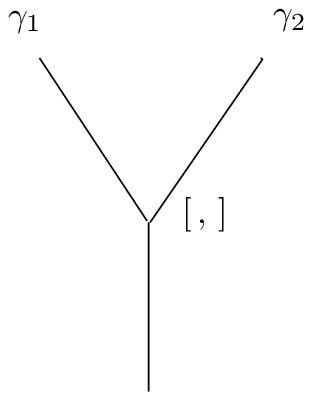} \hfill\\
\parbox[t]{0.45\textwidth}{\caption{The product
$\wedge\in \calc^{\mc}(2,0)$ is used}\label{prod}}
\hfill
\parbox[t]{0.45\textwidth}{\caption{The bracket
$[\,,\,]\in \calc^{\mc}(2,0)$ is used}\label{brack}}
\hfill
\end{center}
\end{figure}

\begin{figure}
\begin{center}
\includegraphics[width=0.45\textwidth]{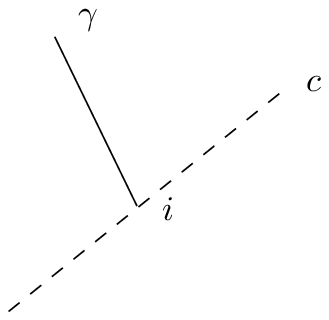}
\hfill
\includegraphics[width=0.45\textwidth]{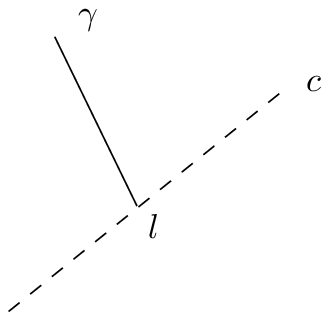} \\
\parbox[t]{0.45\textwidth}{\caption{The commutative module
structure
$i\in \calc^{\ma}(1,1)$ is used}\label{contr}}
\hfill
\parbox[t]{0.45\textwidth}{\caption{The Lie algebra module structure
$l\in \calc^{\ma}(1,1)$ is used}\label{Lieder}}
\end{center}
\end{figure}

\begin{figure}
\psfrag{d}[]{$\delta$}
\begin{center}
\hfill
\includegraphics[width=0.45\textwidth]{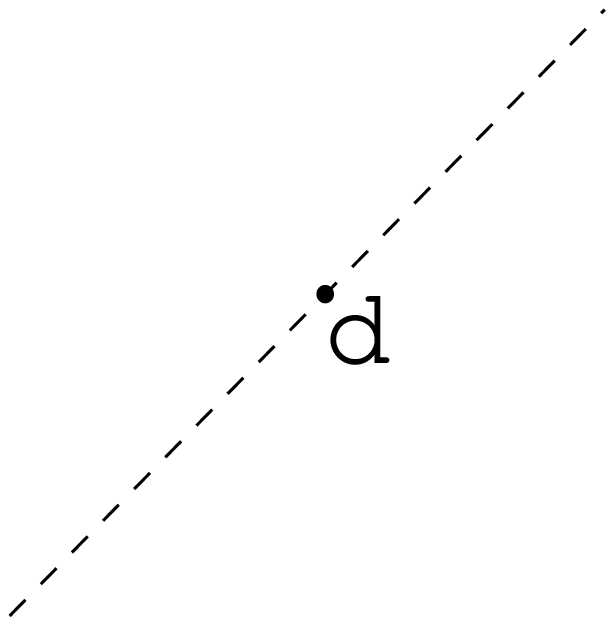}
\hfill
\includegraphics[width=0.45\textwidth]{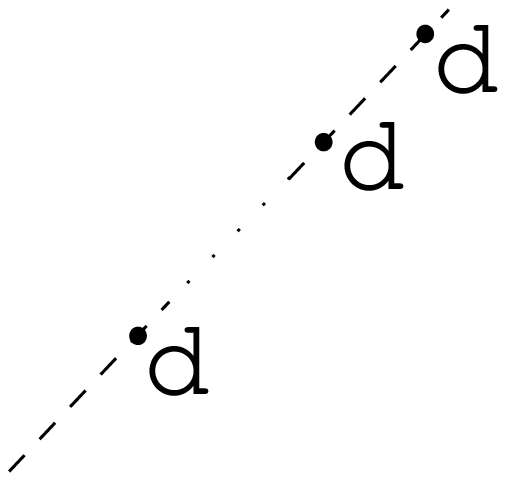} \\
\hfill
\parbox[t]{0.45\textwidth}{\caption{The unary operation
$\de\in \calc^{\ma}(0,1)$ is used}\label{d}}
\hfill
\parbox[t]{0.45\textwidth}{\caption{
The number of $\de$'s on the figure is $m$}
\label{u-m}}
\hfill
\end{center}
\end{figure}

We represent elements of the free coalgebra
$\bbF_{\bB}(\cV, \cW)$ and elements of the cooperad
$\bB$ graphically.
Thus Figures \ref{prod}, \ref{brack} represent the
simplest elements of $\bbF_{\bB}(\cV, \cW)_{\mc}$
with $\ga_1$ and $\ga_2$ being vectors in $\cV$\,.
Figures \ref{contr}, \ref{Lieder} show the simplest
elements of $\bbF_{\bB}(\cV, \cW)_{\ma}$ with $\ga\in \cV$
and $c\in \cW$\,.
Figures \ref{d} and \ref{u-m} represent simple
elements of $\bB^{\ma}(0,1)$\,.
The dashed line in figures \ref{contr},
\ref{Lieder}, \ref{d}, and \ref{u-m} is used to label
the arguments of the color $\ma$ and
the solid line is used to label
the arguments of the color $\mc$\,.

\begin{figure}
\psfrag{d}[]{$\delta$}
\psfrag{di}[]{$\delta\,i$}
\psfrag{id}[]{$i\,\delta$}
\psfrag{D}[]{$\partial^{Bar}$}
\psfrag{p}[]{$+$}
\psfrag{eq}[]{$=$}
\psfrag{g}[]{$\gamma$}
\psfrag{c}[]{$c$}
\psfrag{z}[]{$-(-1)^{|\ga|}$}
\psfrag{l}[]{$l$}
\psfrag{i}[]{$i$}
\begin{center}
\includegraphics[width=0.9\textwidth]{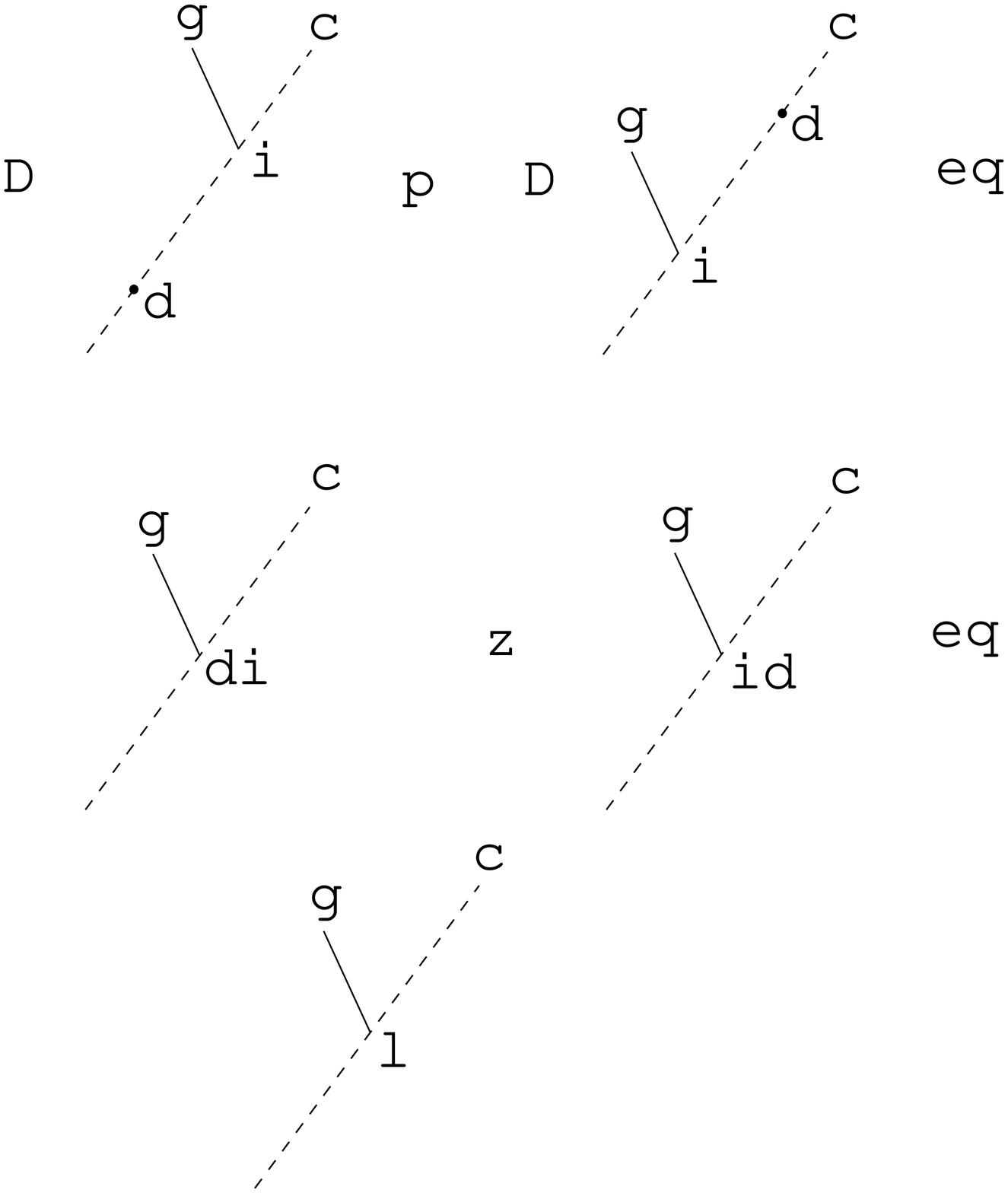}\\[0.5cm]
\parbox[t]{0.6\textwidth}{\caption{A simple computation with $\pa^{Bar}$}
\label{formula}}
\end{center}
\end{figure}

Using this graphical notation we may perform
simple computations in the coalgebra $\bbF_{\bB}(\cV, \cW)$.
For example, using equation (\ref{l-i-delta}), we present
on Figure \ref{formula} a simple computation
with the bar differential $\pa^{Bar}$\,.
Here $\ga\in \cV$ and $c\in \cW$\,.

For the operad $\Ger$ we use a resolution which is
simpler than the canonical one
$Cobar(Bar(\Ger))$\,. More precisely,
as in \cite{BLT}, we set
\begin{equation}
\label{HoGer}
\Hoger = Cobar(\Ger^{\vee})\,.
\end{equation}
Due to koszulity of the operad $\Ger$
the inclusions
\begin{equation}
\label{Ger-Koszul}
\io_{\Ger} : \Ger^{\vee} \hookrightarrow Bar(\Ger)
\end{equation}
and
\begin{equation}
\label{Ger-Koszul1}
Cobar(\io_{\Ger}) :
Cobar(\Ger^{\vee}) \hookrightarrow
Cobar(Bar(\Ger))
\end{equation}
are quasi-isomorphisms of cooperads
and operads, respectively. It is the
second quasi-isomorphism (\ref{Ger-Koszul1})
which allows us to replace the canonical
resolution $Cobar(Bar(\Ger))$ by (\ref{HoGer}).

To get a more tractable description of
algebras over the operads $\Hocalc$ and
$\Hoger$ we introduce the following DGLAs
\begin{equation}
\label{coder-Ger}
\Coder'(\bbF_{\Ger^{\vee}}(\cV)) =
\{Q \in \Coder(\bbF_{\Ger^{\vee}}(\cV)) ~|~
Q\Big|_{\cV}  = 0 \}\,,
\end{equation}
\begin{equation}
\label{coder-calc}
\Coder'(\bbF_{\bB}(\cV,\cW)) =
\{Q \in \Coder(\bbF_{\bB}(\cV, \cW)) ~|~
Q\Big|_{\cV \oplus \cW}  = 0 \}\,,
\end{equation}
where $\Coder(\bbF_{\Ger^{\vee}}(\cV))$
(resp. $\Coder(\bbF_{\bB}(\cV, \cW))$ ) is the
DGLA of coderivations of the free coalgebra
$\bbF_{\Ger^{\vee}}(\cV)$ (resp. the
free coalgebra $\bbF_{\bB}(\cV, \cW)$).
Furthermore, $\cV$ (resp. $\cV \oplus \cW$) is considered
as a subspace of $\bbF_{\Ger^{\vee}}(\cV)$
(resp.  $\bbF_{\bB}(\cV, \cW)$) via the
corresponding coaugmentation.

According to Proposition \ref{GJ} the $\Hoger$-algebra
structures on $\cV$ are in bijection with
the Maurer-Cartan elements of the DGLA
$\Coder'(\bbF_{\Ger^{\vee}}(\cV))$\,.
Similarly, the $\Hocalc$-algebra
structures on the pair $(\cV, \cW)$ are
in bijection with the Maurer-Cartan elements of the DGLA
$\Coder'(\bbF_{\bB}(\cV,\cW))$\,.
Moreover, due to Proposition \ref{coder-cofree}
the Maurer-Cartan element $Q$ of the DGLA (\ref{coder-Ger})
(resp. the DGLA (\ref{coder-calc})) is uniquely determined
by its composition $\rho_{\cV}\circ Q$
(resp. $\rho_{\cV, \cW} \circ Q$)
with the corestriction
$\rho_{\cV}: \bbF_{\Ger^{\vee}}(\cV) \to \cV$
(resp. the corestriction
$\rho_{\cV, \cW}: \bbF_{\bB}(\cV, \cW)
\to \cV\oplus \cW$)\,.

The vector space of the
free coalgebra $\bbF_{\bB}(\cV, \cW)$
splits according to the two colors $(\mc, \ma)$ as
\begin{equation}
\label{bB-mc-ma}
\bbF_{\bB}(\cV, \cW) = \bbF_{\bB}(\cV, \cW)_{\mc}
\oplus \bbF_{\bB}(\cV, \cW)_{\ma}\,,
\end{equation}
where
\begin{equation}
\label{bB-Ger}
\bbF_{\bB}(\cV, \cW)_{\mc} = \bbF_{Bar(\Ger)}(\cV)\,.
\end{equation}

Thus for every $\Hocalc$-algebra $(\cV, \cW)$
the graded vector space $\cV$ is an algebra over
the operad $Cobar(Bar(\Ger))$\,.
Using this algebra structure over $Cobar(Bar(\Ger))$
and the embedding (\ref{Ger-Koszul1}) we get a
$\Hoger$-algebra structure on $\cV$\,.

To describe the relationship between
these algebras we denote by
$Q_{\cV, \cW}$ the Maurer-Cartan element of
the DGLA $\Coder'(\bbF_{\bB}(\cV,\cW))$ corresponding
to the $\Hocalc$-algebra structure
on $(\cV, \cW)$\,.
Next, we denote by $Q_{\cV}$ the Maurer-Cartan element of
the DGLA $\Coder'(\bbF_{\Ger^{\vee}}(\cV))$ corresponding
to the $\Hoger$-algebra structure
on $\cV$\,.

Proposition \ref{how-to-pull} implies that
\begin{equation}
\label{calc-Ger}
\rho_{\cV} \circ Q_{\cV} =
\rho_{\cV \oplus \cW} \circ Q^{\mc}_{\cV, \cW}
\circ \bbF(\io_{\Ger})\,,
\end{equation}
where $\io_{\Ger}$ is the embedding (\ref{Ger-Koszul})
and
$$
Q^{\mc}_{\cV, \cW} =  Q_{\cV, \cW}
\Big|_{\bbF_{\bB}(\cV, \cW)_{\mc}} \,.
$$

Due to Proposition \ref{coder-cofree}
the coderivation $Q_{\cV}$ (resp. the coderivation
$Q_{\cV, \cW}$) is uniquely determined by
the composition $\rho_{\cV} \circ Q_{\cV}$
(resp. $\rho_{\cV, \cW} \circ Q_{\cV, \cW}$)\,.
Thus equation (\ref{calc-Ger}) indeed describes
the relationship between the $\Hocalc$-algebra structure
on $(\cV, \cW)$ and the $\Hoger$-algebra structure
on $\cV$\,.

~\\
{\bf Remark.} The vector space of operations of
the cooperad $\bB$ with no arguments
having color $\mc$ is
\begin{equation}
\label{bB-01}
\bB^{\ma}(0,1) = \bbK[u]\,,
\end{equation}
where $u$ is an auxiliary variable of degree $-2$\,.
The monomial $u^m$ corresponds to the element
of $\bB^{\ma}(0,1)$ which is drawn on Figure
\ref{u-m} (See page \pageref{u-m}).

\subsection{Description of the operad $\Ho(\Lie^+_{\de})$}
The canonical cofibrant resolution
$Cobar(Bar(\Lie^+_{\de}))$ can be simplified.
In this subsection we construct a sub-cooperad
$(\Lie^+_{\de})^{\vee}$ of $Bar(\Lie^+_{\de})$
such that the embedding of operads
$$
Cobar((\Lie^+_{\de})^{\vee}) \hookrightarrow
Cobar(Bar(\Lie^+_{\de}))
$$
is a quasi-isomorphism.
This construction goes along the lines of
\cite{Fresse}, \cite{GK}. (See also
Definition 3.2.1 in \cite{Hinich}.)
It allows us to set
$$
\Ho(\Lie^+_{\de}) =
Cobar((\Lie^+_{\de})^{\vee})\,.
$$

Let us first recall that algebras
over the operad $\Lie^+_{\de}$ are pairs $(\cV, \cW)$
where $\cV$ is a Lie algebra $\cW$ is
a Lie algebra module over $\cV$ and
$\cW$ is equipped with degree $-1$ unary
operation $\de$ which satisfies the identities
\begin{equation}
\label{delta-2}
\de^2 = 0\,.
\end{equation}
and
\begin{equation}
\label{l-delta1}
\delta\, l_a - (-1)^{|a|} l_a \, \de = 0\,,
\end{equation}
where $l: \cV \otimes \cW \to \cW$ is
the action of $\cV$ on $\cW$\,.

Thus the operad $\Lie^+_{\de}$ is generated
by the elementary operations $[\,,\,]$\,,
$l$ and $\de$\,, where $[\,,\,]$ denotes
the Lie bracket. These operations are subject
to the homogeneous quadratic relations:
the Jacobi identity for
the Lie bracket $[\,,\,]$\,, and the compatibility
equation between $l$ and $[\,,\,]$
\begin{equation}
\label{l-brack}
l_a l_b - (-1)^{|a|\, |b|} l_a l_b = l_{[a,b]}
\end{equation}
and, finally, equations (\ref{delta-2}) and
(\ref{l-delta1})\,.

To construct the cooperad $(\Lie^+_{\de})^{\vee}$
we introduce the polynomial functor $S$ spanned linearly
by the elementary operations $[\,,\,]$\,, $l$\,,
and $\de$ of the operad $\Lie^+_{\de}$.

We also introduce the linear span $R$ of
the homogeneous quadratic relations of $\Lie^+_{\de}$
between the elementary operations.

Next, we consider
the free cooperad $\bbT^*(\bs^{-1} S)$ generated
by the polynomial functor $\bs^{-1} S$\,.
The cooperad $\bbT^*(\bs^{-1} S)$ may be viewed
as a sub-cooperad of $Bar(\Lie^+_{\de})$
if we forget about the differential $\pa^{Bar}$\,.

Let us remark that, the cooperad $\bbT^*(\bs^{-1} S)$
is equipped with the natural grading
\begin{equation}
\label{grading}
\bbT^*(\bs^{-1} S) =
\bigoplus_{m=0}^{\infty} \bbT^*_m(\bs^{-1} S)\,,
\qquad
\bbT^*_0(\bs^{-1} S) = \ast\,,
\end{equation}
where $\ast$ is the terminal object (\ref{ast-color})
in the category of $2$-colored cooperads
and $\bbT^*_m(\bs^{-1} S)$ consists of the elements
of degree $m$ in the elementary operations.
Thus, since the relations between the elementary
operations are quadratic, $\bs^{-2} R$ is a subspace
of $\bbT^*_2(\bs^{-1} S)$\,.

First, we construct the cooperad
$(\Lie^+_{\de})^{\vee}$ as
a sub-cooperad of $\bbT^*(\bs^{-1} S)$ and
then we will show
that $(\Lie^+_{\de})^{\vee}$ belongs to the
kernel of the bar
differential $\pa^{Bar}$\,.

We construct $(\Lie^+_{\de})^{\vee}$
by induction on the degree $m$ in (\ref{grading})\,.
The base of the induction is given by
the equations
\begin{equation}
\label{base1}
(\Lie^+_{\de})^{\vee} \cap \bbT^*_0(\bs^{-1} S)
\oplus \bbT^*_1(\bs^{-1} S)=
 \bbT^*_0(\bs^{-1} S) \oplus \bbT^*_1(\bs^{-1} S)\,,
\end{equation}
\begin{equation}
\label{base11}
(\Lie^+_{\de})^{\vee} \cap \bbT^*_2(\bs^{-1} S) =
\bs^{-2} R\,,
\end{equation}
and the step is given by the condition: {\it
a vector $v\in \bbT^*_m(\bs^{-1} S)$ belongs
to $(\Lie^+_{\de})^{\vee}$  provided}
$$
\tDe (v) \in  (\Lie^+_{\de})^{\vee}
\bul  (\Lie^+_{\de})^{\vee} \,.
$$
Here $\D$ is the coproduct:
$$
\D : \bbT^*(\bs^{-1} S) \to  \bbT^*(\bs^{-1} S) \bullet
\bbT^*(\bs^{-1} S)\,,
$$
and
$$
\tDe(v) = \D(v) -
v \otimes (1 \otimes \dots \otimes 1) -
1 \otimes (v \otimes 1 \otimes  \dots \otimes 1)-
1 \otimes (1 \otimes v \otimes 1 \dots \otimes 1) -
\dots
$$
$$
- 1 \otimes (1 \otimes \dots \otimes 1 \otimes v)\,.
$$

By construction $(\Lie^+_{\de})^{\vee} $ is a sub-cooperad
of $\bbT^*(\bs^{-1} S)$\,.

Equation (\ref{base11}) imply immediately that
$$
\pa^{Bar}\, v = 0\,, \qquad
\forall ~~
v \in
(\Lie^+_{\de})^{\vee} \cap
\bbT^*_2(\bs^{-1} S) \,.
$$

Then the compatibility of $\pa^{Bar}$ with
the coproduct $\D$:
$$
\D\, \pa^{Bar} = (
\pa^{Bar} \otimes (1 \otimes \dots \otimes 1) +
1 \otimes (\pa^{Bar} \otimes 1 \otimes  \dots \otimes 1)+
\dots )\, \D
$$
and the inductive definition of $(\Lie^+_{\de})^{\vee} $
imply that
\begin{equation}
\label{pa-Bar-Lie-vee}
\pa^{Bar}\, v = 0\,, \qquad
\forall ~~
v \in
(\Lie^+_{\de})^{\vee}\,.
\end{equation}

Thus $(\Lie^+_{\de})^{\vee}$ belongs
to the kernel of the bar differential
$\pa^{Bar}$ in $Bar(\Lie^+_{\de})$\,.

The following proposition gives us a
description of the coalgebras
over the cooperad $(\Lie^+_{\de})^{\vee}$
\begin{pred}
\label{lie-de-vee}
A pair $(\cV, \cW)$ of graded vector spaces forms
a coalgebra over the cooperad $(\Lie^+_{\de})^{\vee}$
if $(\cV, \cW)$ is a coalgebra over the cooperad
$\La\cocomm^+$ and $\cW$ is equipped with a
degree $2$ endomorphism
$$
\de^{\vee} : \cW \to \cW
$$
satisfying the equation
$$
l^{\vee} \circ  \de^{\vee} =
(1 \otimes \de^{\vee}) l^{\vee}\,,
$$
where $l^{\vee}$ is the coaction of $\cV$ on $\cW$
$$
l^{\vee} : \cW \to \bs^{-1}( \cV \otimes \cW )\,.
$$
\end{pred}
{\bf Proof.}
Let us consider the restricted dual vector space
\begin{equation}
\label{dual-T-star}
[\bbT^*(\bs^{-1} S)]^{*} =
\Hom_{\restr}(\bbT^*(\bs^{-1} S), \bbK)
\end{equation}
of the free cooperad $\bbT^*(\bs^{-1} S)$
with respect to the grading (\ref{grading})\,.
It is not hard to see that
$$
[\bbT^*(\bs^{-1} S)]^{*} = \bbT(\bs S^*)
$$
is the free operad $\bbT(\bs S^*)$ generated by the
suspension $\bs S^*$ of the linear dual $S^*$ of the
polynomial functor $S$\,.

From the construction of $(\Lie^+_{\de})^{\vee}$ it
follows that the restricted dual
$[(\Lie^+_{\de})^{\vee}]^*$ of the cooperad $(\Lie^+_{\de})^{\vee}$
is the quotient of the free operad $\bbT(\bs S^*)$ with respect to the
ideal generated by the polynomial functor of dual relations
\begin{equation}
\label{dual-R}
R^* = \{r \in
\Hom(\bbT^*_2(\bs^{-1} S), \bbK)\,, ~~|~~
r\Big|_{R} =0
\}\,.
\end{equation}

Let $\{[\,,\,]^*, l^*, \de^{*}\}$ be the basis of
$S^*$ which is dual to
the basis $\{[\,,\,], l, \de\}$ of $S$\,.

Dualizing the Jacobi relation for $[\,,\,]$ and the compatibility
(\ref{l-brack}) of $l$ with $[\,,\,]$ we see that the operation
$\bs[\,,\,]^*$ satisfies the axioms of an associative commutative
product and the operation $\bs l^*$ satisfies the axiom of a module over
an associative and commutative algebra. Dualizing the relation
(\ref{l-delta1}) we see that $\bs l^*$ and $\bs \de^*$ are compatible
in the sense of the following relation
\begin{equation}
\label{l-delta-star}
\bs\de^* \,  \bs l^*   = \bs l^* \, (1 \otimes \bs\de^*) \,.
\end{equation}
Finally the presence of the relation (\ref{delta-2}) implies
that we should not impose any additional condition on $\bs \de^*$
besides (\ref{l-delta-star})\,.

Thus a pair $(\tcV, \tcW)$ is an algebra over the operad
$[(\Lie^+_{\de})^{\vee}]^*$ if $(\tcV, \tcW)$ is a
$\La^{-1}\comm^+$-algebra and $\tcW$ is equipped with a degree $2$ endomorphism
$\bs \de^*$ which is compatible with the action of $\tcV$
on $\tcW$ in the sense of (\ref{l-delta-star}).

Taking the dual partner of an algebra over the operad
$[(\Lie^+_{\de})^{\vee}]^*$ we get the statement of
the proposition.  $\Box$

Proposition \ref{lie-de-vee} implies that a free coalgebra
over the cooperad $(\Lie^+_{\de})^{\vee}$
generated by a pair $(\cV, \cW)$ is
\begin{equation}
\label{liede-vee}
\bbF_{(\Lie^+_{\de})^{\vee}}(\cV,\cW)=
\bbF_{\La\cocomm^+}(\cV,\cW[[u]])\,,
\end{equation}
where $u$ is an auxiliary variable of degree $-2$\,.

We claim that
\begin{pred}
\label{Lie-de-vee-Kos}
The operad $\Lie^+_{\de}$ is Koszul.
In other words the embedding
$$
Cobar((\Lie^+_{\de})^{\vee}) \to
Cobar (Bar(\Lie^+_{\de}))
$$
is a quasi-isomorphism of operads.
\end{pred}
{\bf Proof.}
The criterion of Ginzburg and
Kapranov \cite{GK} (theorem $4.2.5$)
reduces this question to computation
of the homology of a free  $\Lie^+_{\de}$-algebra.
More precisely, we need to show that for
every pair $(\cV, \cW)$
of vector spaces
the complex
\begin{equation}
\label{complex-of-free}
\bbF_{(\Lie^+_{\de})^{\vee}} \circ \bbF_{\Lie^+_{\de}} (\cV,\cW)
\end{equation}
has nontrivial cohomology only in degree $0$\,.

Here the differential on the complex (\ref{complex-of-free})
is defined along the lines of \cite{GJ} using the twisting
cochain of the pair $(\Lie^+_{\de}, (\Lie^+_{\de})^{\vee})$\,.
(See Section 2.4 in \cite{GJ} for more details.)

If we split the complex (\ref{complex-of-free})
according to the colors $\mc$ and $\ma$ and
use equation (\ref{liede-vee}) then we get two
complexes:
\begin{equation}
\label{complex-mc}
\bbF_{(\Lie^+_{\de})^{\vee}} \circ \bbF_{\Lie^+_{\de}} (\cV,\cW)_{\mc}
= \bbF_{\La\cocomm} \circ \bbF_{\Lie}(\cV)\,,
\end{equation}
and
\begin{equation}
\label{complex-ma}
\bbF_{(\Lie^+_{\de})^{\vee}} \circ \bbF_{\Lie^+_{\de}} (\cV,\cW)_{\ma}
= \bbF_{\La\cocomm^+} (\bbF_{\Lie}(\cV), T(\cV) \otimes (\cW \oplus \de \cW)
[[u]])_{\ma}\,,
\end{equation}
where $T(\cV)$ denotes the tensor algebra of $\cV$\,,
$\de$ is the unary operation of $\Lie^+_{\de}$ and
$u$ is an auxiliary variable of degree $-2$\,.

The first complex is exactly the Harrison complex of
the free Lie algebra generated by $\cV$ and it is known
that this complex has nontrivial cohomology only
in degree $0$\,.

The second complex is the tensor product of the Harrison
complex of the free module generated by $\cW$ over the free
Lie algebra $\bbF_{\Lie}(\cV)$ and the De Rham complex
$$
(\bbK[[u]] \oplus \de\, \bbK[[u]], \de \frac{\pa}{\pa u} )
$$
of the algebra $\bbK[[u]]$\,. Thus the second complex also
has nontrivial cohomology only in degree $0$\,. $\Box$

This Proposition implies immediately that the
embedding
$$
Cobar(\La(\Lie^+_{\de})^{\vee}) \hookrightarrow
Cobar (Bar(\La\Lie^+_{\de}))
$$
is a quasi-isomorphism of operads.
Thus we may set
\begin{equation}
\label{HoLaLie}
\Ho(\La\Lie^+_{\de}) = Cobar(\La(\Lie^+_{\de})^{\vee})\,.
\end{equation}

We would also like to remark that
equation (\ref{liede-vee}) implies that
\begin{equation}
\label{liede-vee1}
\bbF_{\La(\Lie^+_{\de})^{\vee}}(\cV,\cW)=
\bbF_{\La^2\cocomm^+}(\cV,\cW[[u]])\,,
\end{equation}
where $u$ is an auxiliary variable of degree $-2$\,.

\section{The Kontsevich-Soibelman operad and the
operad of little discs on a cylinder}

\subsection{The Kontsevich-Soibelman operad $\KS$}
Let us describe the auxiliary operad $\cH$ (of sets) of
``natural''\footnote{We are not sure if these
operations are natural
in the sense of category theory.} operations on the
pair
$$
(\Cbu(A), \Cbd(A))\,.
$$
This operad is
going to have a countable set of colors
\begin{equation}
\label{colors}
\Xi = \bbZ^+\, \sqcup\, \bbZ^- \,,
\end{equation}
where $\bbZ^+$ (resp. $\bbZ^-$)
denotes the set of nonnegative (resp. nonpositive)
integers.

The numbers from the set $\bbZ^+$
label the degrees of
the Hochschild cochains and the numbers
from the set $\bbZ^-$ label the
degrees of Hochschild chains.

Using $\cH$ we construct the DG
operad $\KS$ of Kontsevich
and Soibelman. The latter
operad\footnote{In \cite{K-Soi1} this operad
is denoted by $P$.} is
described in sections 11.1, 11.2 and 11.3 of \cite{K-Soi1}.

For the $\Xi$-colored operad $\cH$
we only allow the operations in which
a chain may enter as an argument at most once. If a
chain enters then the result of the
operation is also a chain. Otherwise
the result is a cochain.
We denote the set of operations
producing a cochain from $n$ cochains by
$\cH(n,0)$\,. The set of operations
producing a chain from $n$ cochains and
$1$ chain is denoted by $\cH(n,1)$\,.

$\cH(n,0)$ is the set of equivalence classes
of rooted\footnote{Recall that a tree called {\it rooted} if
if its root vertex has valency $1$\,.} planar trees $T$ with marked vertices.
The equivalence relation is
the finest one in which two such trees are
equivalent if one of them can be obtained
from the other by either:

\begin{itemize}
\item the contraction of
an edge  with unmarked ends or

\item removing an unmarked vertex with only
one edge originating from it and joining the
two edges adjacent to this vertex into one
edge.
\end{itemize}

If a marked vertex is internal then it is reserved
for a cochain which enters as an argument of the operation.
The number of the incoming edges of such vertex
is the degree of the corresponding cochain.
If a marked vertex is terminal then it is reserved
either for a cochain of degree $0$ or for an argument
of the cochain produced by the operation.

The unmarked vertices (both internal and
terminal) are reserved for the
operations of the non-symmetric operad
$\ass$ which controls unital monoids.
For example, an unmarked terminal vertex is reserved for
unit of $A$, an unmarked vertex of valency $2$
is reserved for the identity
transformation on $A$\,, and an unmarked vertex
of valency $3$ is reserved for the associative
product on $A$\,.

The root vertex is special. Since our trees are rooted
this vertex has always valency $1$. It is always marked and reserved
for the outcome of the cochain
produced by the operation corresponding to the tree.

\begin{figure}
\psfrag{Q}[]{$Q$}
\psfrag{P}[]{$P$}
\psfrag{a1}[]{$a_1$}
\psfrag{a2}[]{$a_2$}
\begin{center}
\includegraphics[width=0.3\textwidth]{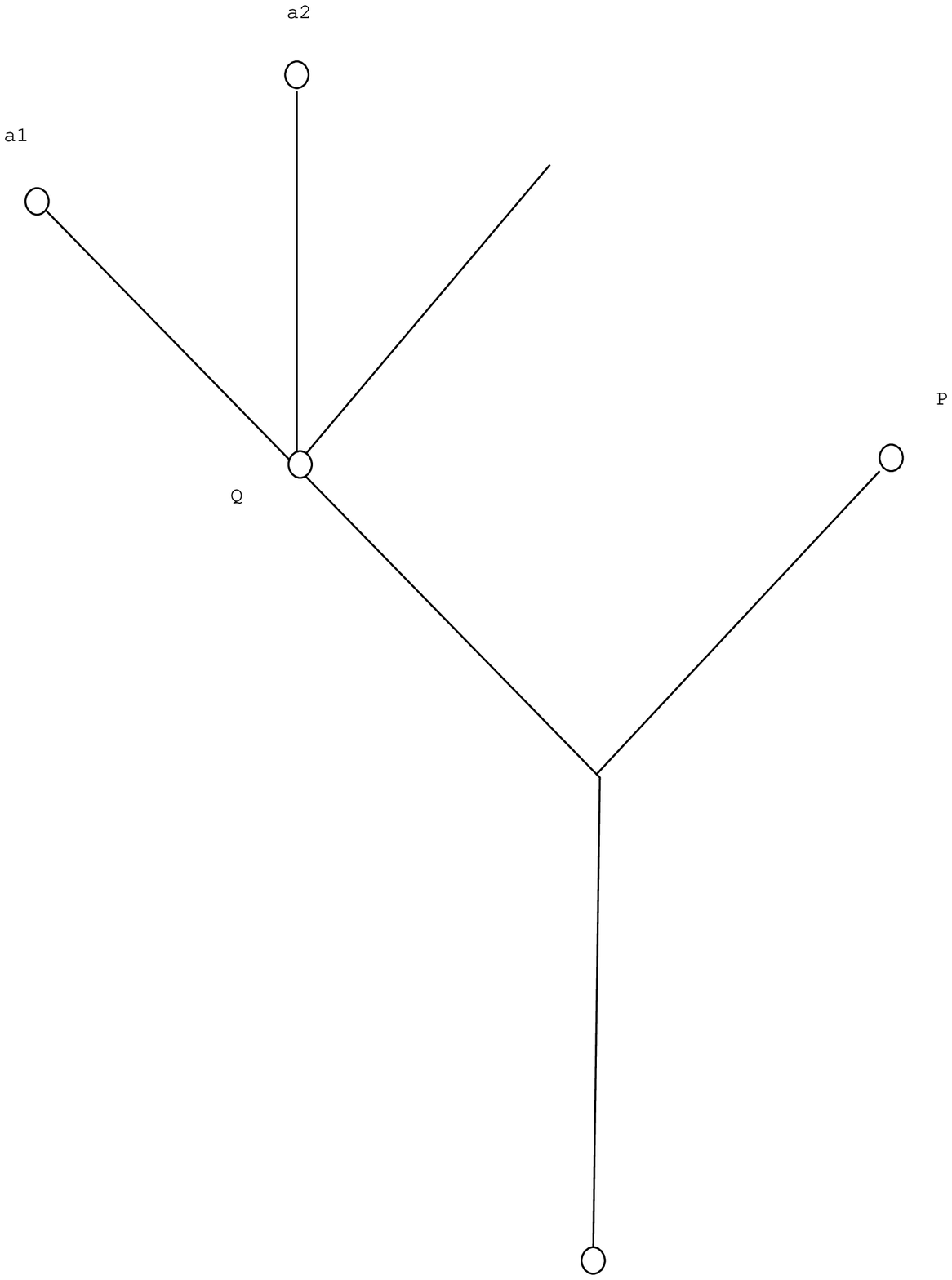}\\[1cm]
\parbox[c]{0.5\textwidth}{\caption{}\label{primer11}}
\end{center}
\end{figure}

The tree on figure \ref{primer11} represents
an operation which produces the $2$-cochain:
$$
a_1\otimes a_2 \to
Q(a_1, a_2, 1) P
$$
from a degree $0$ cochain $P$ and a degree $3$
cochain $Q$.
Marked vertices in this figure are
labeled by small circles. The unmarked
terminal vertex corresponds to the insertion
of the unit into $Q(a_1, a_2, 1)$\,. The unmarked
3-valent vertex gives the product of $P$ and
$Q(a_1, a_2, 1)$\,.

Let us denote by $\cH^{m_a}_{m_r}(n,1)$ the set of operations
producing a chain in $C_{-m_r}(A)$ from $n$ cochains and
a chain in $C_{-m_a}(A)$\,.

$\cH^{m_a}_{m_r}(n,1)$ is described using
forests of rooted trees drawn on the standard cylinder
\begin{equation}
\label{Si}
\Si= S^1 \times [0,1]
\end{equation}
and subject to the following conditions:

\begin{enumerate}

\item every tree of the forest has its root
vertex on the boundary $S^1\times \{0\}$\,,

\item all vertices of the forest lying on the
boundary of the cylinder are marked:

\begin{itemize}

\item the vertices lying on the boundary $S^1\times \{1\}$
are marked by integers $0, 1, \dots, m_a$ in the
counterclockwise order; these vertices are reserved
for the components of the chain which enters as
an argument,

\item the roots are marked by integers
$0, 1, \dots, m_r$  in the same
counterclockwise order; they are reserved for
components of the resulting chain,

\end{itemize}

\item all other marked vertices of the forest
lie on the lateral surface $S^1\times (0,1)$ of the
cylinder and there are exactly $n$ such
marked vertices.

\end{enumerate}

On the set of these forests we introduce
the finest equivalence relation
in which two such forests are
equivalent if one of them can be obtained
from the other by either:

\begin{itemize}

\item isotopy, or

\item the contraction of
an edge with unmarked ends, or

\item removing an unmarked vertex with only
one edge originating from it and joining the
two edges adjacent to this vertex into one
edge.
\end{itemize}
$\cH^{m_a}_{m_r}(n,1)$ is the set of the
corresponding equivalence classes.

As we see from the conditions, all
unmarked vertices lie on the lateral surface
$S^1\times (0,1)$ of the cylinder.
As above, these vertices are
reserved for operations of $\ass$\,.
The marked vertices lying on the lateral surface
$S^1\times (0,1)$ are reserved for cochains.

We allow forests with no marked vertices
lying on the lateral surface $S^1\times (0,1)$\,. Such
forests represent operations which produce a chain
from a chain.

\begin{figure}
\psfrag{Q}[]{$Q$}
\psfrag{P}[]{$P$}
\psfrag{a0}[]{$a_0$}
\psfrag{a1}[]{$a_1$}
\psfrag{a2}[]{$a_2$}
\psfrag{a3}[]{$a_3$}
\psfrag{b0}[]{$b_0$}
\psfrag{b1}[]{$b_1$}
\psfrag{b2}[]{$b_2$}
\psfrag{b3}[]{$b_3$}
\begin{center}
\includegraphics[width=0.5\textwidth]{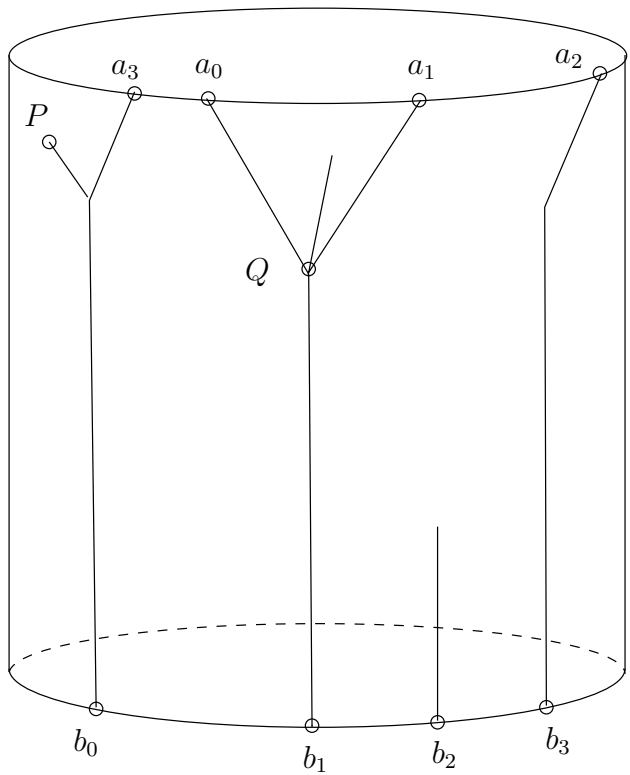}\\[1cm]
\parbox[c]{0.6\textwidth}{\caption{An operation produces which
the chain (\ref{chain})}\label{primer1}}
\end{center}
\end{figure}

Figure \ref{primer1} gives an example of
an operation of $\cH(2,1)$ which produces the
chain
\begin{equation}
\label{chain}
(b_0, b_1, b_2, b_3) =
(P a_3,  Q(a_0,1, a_1), 1, a_2)
\end{equation}
from a degree $0$ cochain $P$, a
degree $3$ cochain $Q$ and a degree $-3$ chain
$(a_0, a_1, a_2, a_3)$\,.
Marked vertices in this figure are
labeled by small circles. The unmarked
3-valent vertex gives the product of $P$
and $a_3$\,, the two unmarked terminal vertices
give units of $A$ and the unmarked 2-valent vertex
gives the identity operation on $A$\,.
The vertices lying on the boundary $S^1\times \{1\}$
are marked by the components of the chain
$(a_0, a_1, a_2, a_3)$ and the roots are marked
by the components of the chain (\ref{chain}).

It is clear how the operad $\cH$ acts on the
pair $(\Cbu(A), \Cbd(A))$\,. From this action it is also
clear how to compose the operations.
For example, the composition of operations
from $\cH(n_1,1)$ and $\cH(n_2,1)$ corresponds
to putting one cylinder on the top of
the other matching the roots of the first cylinder
with the vertices lying on the upper circle
of the second cylinder, and then shrinking
the resulting cylinder to the required height.

Recall that the operad $\cH$ is colored by degrees
of the cochains and degrees of the chains.
It is not hard to see that $\cH(n,0)$ is
a cosimplicial set with respect to the degree of
the resulting cochain and a polysimplicial set with
respect to the degrees of the cochains entering
as arguments.

Similarly, $\cH(n,1)$ is a cosimplicial set with respect
to the degree of the chain entering as an argument
and a polysimplicial set with
respect to the degrees of the cochains entering
as arguments and the degree of the resulting chain.

These poly-simplicial/cosimplicial structure
is compatible with the compositions and we get
\begin{defi}
The DG operad $\KS$ is the realization of
the operad $\cH$ in the category of chain
complexes.
\end{defi}
It follows from the construction that
$\KS$ is a $2$-colored operad which acts on
the pair $(\nCbu(A), \nCbd(A))$\,.

It is not hard to see that the operations $\cup$ (\ref{cup}),
$[\,,\,]_G$ (\ref{Gerst}), $I$ (\ref{I-P}), $L$ (\ref{L-Q})
and $B$ (\ref{B}) come from the action of the operad $\KS$\,.

~\\
{\bf Remark 1.} The operad $\KS$ with its action on
$$
(\nCbu(A), \nCbd(A))
$$
was introduced by Kontsevich and Soibelman
in\footnote{See sections 11.1, 11.2, and
11.3 in \cite{K-Soi1}.}
\cite{K-Soi1} in the case when $A$ is an
$A_{\infty}$-algebra.
Here we recall the construction
of $\KS$ in the case when $A$ is simply an associative
algebra. It is this assumption on $A$ which allows us to utilize
the natural cosimplicial/simplicial
structure on $(\nCbu(A), \nCbd(A))$\,.

~\\
{\bf Remark 2.}
If we restrict ourselves to the subspace
of operations of $\KS$ which do not involve
chains then we get the minimal operad of
Kontsevich and Soibelman described in
\cite{K-Soi}.

\subsection{The operad of little discs on a cylinder}
A ``topological partner'' of $\KS$ is
the operad $\Cyl$ of discs on a cylinder
\cite{K-Soi1}, \cite{TT1}. As well as the operad
of Kontsevich and Soibelman $\Cyl$ is a $2$-colored
operad satisfying the property (\ref{P-top}).

The spaces $\Cyl^{\mc}(n,0)$, $n\ge 1$ are the
spaces of the little disc operad.

To introduce the space $\Cyl^{\ma}(n,1)$ for $n\ge 1$ we consider
cylinders $S^1\times [a,c]$ for $a,c \in \bbR$, $a < c$
with the natural flat metric and
define the topological space $\tCyl_n$\,.

A point of the space $\tCyl_n$
is a cylinder $S^1\times [a,c]$ together with a
configuration of $n\ge 1$ discs on
the lateral surface $S^1\times (a,c)$
and a position of two points $b$ and $t$
lying on the boundaries $S^1\times {a}$ and
$S^1\times {c}$\,, respectively.
The topology on the space $\tCyl_n$ is
defined in the obvious way using the flat metric
on the cylinder.

The space $\tCyl_n$ is equipped with a free
action of the group $S^1 \times \bbR$. The subgroup
$S^1\subset S^1 \times \bbR$ simultaneously rotates
all the cylinders and the subgroup
$\bbR \subset S^1 \times \bbR$ acts by parallel shifts
$$
S^1\times [a,c] \to  S^1\times [a+l,c+l]\,, \qquad
 l\in \bbR\,.
$$

The space $\Cyl^{\ma}(n,1)$ for $n\ge 1$ of the operad
$\Cyl$ is the quotient
\begin{equation}
\label{C-col-2}
\Cyl^{\ma}(n,1) = \tCyl_n / S^1 \times \bbR\,.
\end{equation}

The space $\Cyl^{\ma}(0,1)$ is the space of configurations
of two (possibly coinciding) points $b$ and $t$ on
the circle $S^1$ considered modulo rotations.
Although it is obvious that $\Cyl^{\ma}(0,1)$ is
homeomorphic to the circle $S^1$ we still
define $\Cyl^{\ma}(0,1)$ using
the configuration space in order to better visualize
the operations of the operad.

The insertions of the type
$$
\Cyl^{\mc}(n,0)\times \Cyl^{\mc}(m,0) \to \Cyl^{\mc}(n+m-1,0)
$$
are defined in the same as for the operad of
little squares. The operations of the type
$$
\Cyl^{\mc}(n,0)\times \Cyl^{\ma}(m,1) \to \Cyl^{\ma}(n+m-1,1)
$$
are insertions of the configuration of little discs
of $\Cyl^{\mc}(n,0)$ in a little disc on the lateral surface
of a cylinder.
Finally the operations of the type
$$
\Cyl^{\ma}(n,1)\times \Cyl^{\ma}(m,1) \to \Cyl^{\ma}(n+m,1)
$$
correspond to putting the first cylinder under the second one
while the second cylinder is rotated in such a way that the point $b$
of the second cylinder coincides with the point $t$ of the
first cylinder. The composition involving degenerate configurations
of $\Cyl^{\ma}(0,1)$ are defined in the obvious way.

To describe the operad of homology groups of $\Cyl$
we will need some results about the configuration
spaces of distinct points on the punctured
plane $\bbR^2\setminus \{{\bf 0}\}$\,.

Let us denote by $\Conf_n(\bbR^2 \setminus \{ {\bf 0} \})$
the configuration space of $n$ distinct points on the
punctured plane $\bbR^2 \setminus \{{\bf 0}\}$ and
consider the following projections
\begin{equation}
\label{p-k}
p_k :  \Conf_n(\bbR^2 \setminus \{ {\bf 0} \})
\to  \bbR^2 \setminus \{ {\bf 0} \}\,,
\end{equation}
$$
p_k (\bx_1, \dots, \bx_n) = \bx_k\,.
$$

Due to  E. Fadell and L. Neuwirth \cite{FN} we have
\begin{teo}[Theorem 1, \cite{FN}]
\label{FN-teo}
For every $k=1, 2, \dots, n$ the map
$p_k$ is a locally trivial fibration.
\end{teo}

Using the ideas of E. Fadell and
L. Neuwirth \cite{FN} we show that
\begin{pred}
The map
$$
p : \Conf_n(\bbR^2 \setminus \{{\bf 0}\} ) \to
\Conf_{n-1}(\bbR^2 \setminus \{{\bf 0}\} )\,,
$$
\begin{equation}
\label{p}
p (\bx_1, \bx_2, \dots, \bx_n) =
(\bx_2, \dots, \bx_n)\,.
\end{equation}
is a locally trivial fibration.
Furthermore, the fiber $F_n$ of $p$ is
\begin{equation}
\label{F-n}
F_n = \bbR^2 \setminus
\{{\bf 0}, \bq_2, \dots, \bq_{n}\}  \,,
\end{equation}
where $\bq_2, \dots, \bq_{n}$ are $n-1$ distinct points
of the punctured plane $\bbR^2\setminus  \{{\bf 0}\}$\,.
\end{pred}
{\bf Proof.}
For the open unit disc $D_1$ on $\bbR^2$
centered at the origin
there exists a continuous map
\begin{equation}
\label{theta}
\te : D_1 \times \bar{D}_1 \to \bar{D}_1
\end{equation}
satisfying the following properties:

 -- for all $\bx\in D_1$ the map
 $\te(\bx,\,\, ): \bar{D}_1 \to \bar{D}_1 $ is a
 a homeomorphism having $\pa \bar{D}_1$ fixed.

 -- for all $\bx\in D_1$ we have
 $\te(\bx,\bx)={\bf 0}$\,.

For distinct points
$\bq_2, \dots, \bq_{n}$
on the punctured plane $\bbR^2\setminus  \{{\bf 0}\}$
we choose
open discs
\begin{equation}
\label{discs1}
D_{\bq_2}, \, D_{\bq_3}, \,
\dots, \, D_{\bq_{n}}\,,
\end{equation}
which are centered at
$\bq_2, \bq_3, \dots, \bq_{n}$,
respectively.
Each disc $\bar{D}_{\bq_j}$ avoids the origin
${\bf 0}$ and
for each $i \neq j$
$$
\bar{D}_{\bq_i} \cap
\bar{D}_{\bq_j} = \emptyset\,.
$$

Let us denote respectively by $r_2, \dots, r_{n}$
the radii of the discs (\ref{discs1}) and let
$U_b$ be the following neighborhood of
$\Conf_{n-1}(\bbR^2 \setminus \{\bf 0\} )$:
\begin{equation}
\label{U-b}
U_b = \{(\bx_2, \dots, \bx_n) \,\,|\,\,
\bx_2  \in   D_{\bq_2}\,, \quad
\bx_3  \in   D_{\bq_3}\,, \, \dots \,,
\bx_n  \in   D_{\bq_n} \} \,.
\end{equation}

The desired homeomorphism $h$ from
$p^{-1}(U_b) \to F_n\times U_b$ is given by
the formula:
\begin{equation}
\label{h}
h(\bx_1, \bx_2, \dots, \bx_n) =
\begin{cases}
\displaystyle
(\bx_1, \bx_2, \dots, \bx_n)
\,, \qquad {\rm if} ~~ \bx_1 \notin
\bigcup_{j=2}^n D_{\bq_j}\,, \\[0.3cm]
\displaystyle
(\bq_j + r_j\,
\te\Big(\frac{\bx_j-\bq_j}{r_j},
\frac{\bx_1- \bq_j}{r_j}\Big), \bx_2, \dots, \bx_n)\,,
\quad {\rm if} ~~ \bx_1
 \in \bar{D}_{\bq_j}\,,
\end{cases}
\end{equation}
where $r_j$ is the radius of the $j$-th disc
$D_{\bq_j}$\,.

Since the open subsets of the form $U_b$ (\ref{U-b})
cover $\Conf_{n-1}(\bbR^2 \setminus \{\bf 0\} )$
the map $p$ (\ref{p}) is indeed a locally trivial
fibration. $\Box$

Fiber $F_n$ of $p$ (\ref{p}) is
homotopy equivalent to the wedge
sum $\vee^n S^1$ of $n$ circles.
Hence, the homology groups of $F_n$ (\ref{F-n}) are
\begin{equation}
\label{H-F-n}
H_{\bul} (F_n, \bbK) =
\begin{cases}
\bbK \,, \qquad {\rm if} ~~ \bul = 0\,, \\
\bbK^n\,, \qquad {\rm if} ~~ \bul = 1\,, \\
0 \,, \qquad {\rm otherwise}\,.
\end{cases}
\end{equation}

Let us show that
\begin{pred}
\label{acts-triv}
The fundamental group
$\pi_1(\Conf_{n-1}(\bbR^2 \setminus \{\bf 0\} ))$
acts trivially on the homology groups
$H_{\bul} (F_n, \bbK)$ of $F_n$\,.
\end{pred}
{\bf Proof.} It is obvious that we only need
to consider the action of
$\pi_1(\Conf_{n-1}(\bbR^2 \setminus \{\bf 0\} ))$
on $H_1(F_n, \bbK)$\,.

To get the cycles representing elements
of a basis for $H_1(F_n, \bbK)$ we choose
closed discs
\begin{equation}
\label{discs}
D_{\bf 0}, \, D_{\bq_2}, \, D_{\bq_3}, \,
\dots, \, D_{\bq_{n}}\,,
\end{equation}
which are centered at
${\bf 0}, \bq_2, \bq_3, \dots, \bq_{n}$,
respectively. The discs (\ref{discs})
are chosen in such a way that their
closures
$$
\bar{D}_{\bf 0}, \, \bar{D}_{\bq_2}, \, \bar{D}_{\bq_3}, \,
\dots, \, \bar{D}_{\bq_{n}}
$$
are pairwise disjoint.

The boundaries of the discs (\ref{discs})
are cycles
representing the elements of a basis for
$H_1(F_n, \bbK)$\,.

Let us identify $\bbR^2$ with the complex
plane $\bbC$ and consider the following loop
in $\Conf_{n-1}(\bbR^2 \setminus \{\bf 0\} )$
\begin{equation}
\label{loop}
f(t) =
(e^{2 \pi i t}\bq_2, e^{2 \pi i t} \bq_3, \dots,
e^{2 \pi i t} \bq_{n} )
\,:\, [0,1] \to \Conf_{n-1}(\bbR^2 \setminus \{\bf 0\} ) \,,
\end{equation}
where $(\bq_2, \bq_3, \dots, \bq_{n})$ is
a fixed collection of the distinct points
of $\bbR^2 \setminus \{\bf 0\}$\,.

The loop (\ref{loop}) lifts to the following loop
in $\Conf_{n}(\bbR^2 \setminus \{\bf 0\} )$
\begin{equation}
\label{loop1}
\tf(t) =
(e^{2 \pi i t} \bx_1,
e^{2 \pi i t}\bq_2, e^{2 \pi i t} \bq_3, \dots,
e^{2 \pi i t} \bq_{n} )
\,:\, [0,1] \to
\Conf_{n}(\bbR^2 \setminus \{\bf 0\} ) \,.
\end{equation}
As we go around the loop (\ref{loop1}) the
point $\bx_1$ of the fiber $F_n$ returns to
its original position. Thus the element
$[f]\in \pi_1(\Conf_{n-1}(\bbR^2 \setminus \{\bf 0\} ))$
represented by the loop $f$ (\ref{loop}) acts
trivially on $H_{\bul}(F_n,\bbK )$\,.

Let now $g$ be an arbitrary loop in
$\Conf_{n-1}(\bbR^2 \setminus \{\bf 0\} )$\,.

To find the action of the homotopy class $[g]$ on
$H_{\bul}(F_n,\bbK )$ we need to lift the
map
\begin{equation}
\label{gamma}
\ga (y, t) = g(t): F_n \times [0,1] \to
\Conf_{n-1}(\bbR^2 \setminus \{\bf 0\} )
\end{equation}
to a map
\begin{equation}
\label{gamma-t}
\tga (y, t): F_n \times [0,1] \to
\Conf_{n}(\bbR^2 \setminus \{\bf 0\} )
\end{equation}
which makes the diagram
\begin{equation}
\label{gamma-diag}
\begin{array}{ccc}
F_n \times \{0\} & \hookrightarrow & \Conf_{n}(\bbR^2 \setminus \{ {\bf 0} \} ) \\[0.3cm]
\downarrow & ~^{\tga}\nearrow & \downarrow^{p}  \\[0.3cm]
F_n \times [0,1] & \stackrel{\ga}{\rightarrow} &
\Conf_{n-1}(\bbR^2 \setminus \{{\bf 0}\} )
\end{array}
\end{equation}
commutative.

To construct the lift $\tga$ we divide the segment
$[0,1]$ into small enough subsegments $[t_i, t_{i+1}]$ satisfying the property
\begin{equation}
\label{subseg-prop}
g([t_i, t_{i+1}]) \subset V\,,
\end{equation}
where $V$ is an open subset of $\Conf_{n-1}(\bbR^2 \setminus \{\bf 0\} ) $
of the form (\ref{U-b})\,.
Then each individual lift
\begin{equation}
\label{gamma-t-i}
\tga \Big|_{F_n \times [t_i, t_{i+1}]} :
F_n \times [t_i, t_{i+1}] \to  \Conf_{n}(\bbR^2 \setminus \{ {\bf 0} \} )
\end{equation}
can be constructed using the trivialization (\ref{h}).

With this construction in mind
we consider the compositions $p_k \circ g$ of $g$ with
the projections $p_k$ (\ref{p-k}) for $k \in \{2,\, \dots,\, n\}$\,.
Since the image
$p_k \circ g([0,1])$ of the segment $[0,1]$
 is a compact subset in $\bbR^2 \setminus \{ {\bf 0} \}$ we may choose the
disc $D_{\bf 0}$ in such a way that the closure
 $\bar{D}_{\bf 0}$ avoids
the points of the images $p_k \circ g([0,1])$
for all $k \in \{2,\, \dots,\, n\}$\,.
Therefore the lift $\tga$ (\ref{gamma-t}) can be
chosen in such a way that
$$
\tga(y,t) = y\,, \qquad \forall ~~ y\in \bar{D}_{\bf 0}\subset F_n\,.
$$
Hence the action of the homotopy class $[g]$
on the homology class represented by the
boundary of the disc $D_{\bf 0}$ is trivial.

Let us examine the loop
$p_k\circ g$ in $\bbR^2\setminus \{\bf 0\}$ more
closely.

Since the fundamental group of $\bbR^2\setminus \{\bf 0\}$
is generated by the homotopy class of the loop
$$
l(t) = e^{2 \pi i t} \bq_k : [0,1] \to
\bbR^2\setminus \{\bf 0\}
$$
around the origin there exists an integer
$N \in \bbZ$ such that the loop
$p_k \circ (g * f^N)$ is null-homotopic.
Here $f$ is the loop (\ref{loop}) in
$\Conf_{n-1}(\bbR^2 \setminus \{\bf 0\} )$\,.

But the class of $f$ acts trivially on the
homology of the fiber $F_n$\,. Therefore,
without loss of generality,
we may assume that the composition
$p_k \circ g$ is null-homotopic.

Thus, due to Theorem \ref{FN-teo} of Fadell and
Neuwirth we may further assume that
$p_k \circ g$ is a constant map
\begin{equation}
\label{p-k-g}
p_k \circ g (t) \equiv \bq_k\,.
\end{equation}
In other words, the $k$-th point $\bx_k$ does
not move as we go along the loop $g$\,.

Since for each
$j\neq k$ the image
$p_j \circ g([0,1])$ of the segment $[0,1]$
is a compact subset in $\bbR^2 \setminus \{ {\bf 0} \}$
we may choose the
disc $D_{\bq_k}$ in such a way that the closure
 $\bar{D}_{\bq_k}$ avoids
the points of the images $p_j \circ g([0,1])$
for all $j\neq k$\,.

Thus, using the partition of the segment $[0,1]$
satisfying the property (\ref{subseg-prop}) and
constructing the lift (\ref{gamma-t}) using the
trivializations of the form (\ref{h}) we see that
the lift $\tga$ can be chosen in such a way that
$$
\tga(y,t) = y\,, \qquad \forall~~ y\in \bar{D}_{\bq_k}\subset F_n\,.
$$
Therefore the action of the homotopy class $[g]$
on the homology class represented by the
boundary of the disc $D_{\bq_k}$ is trivial.

The proposition is proved. $\Box$

Generalizing the result of F. R. Cohen \cite{Cohen} we
get the following
\begin{teo}
\label{H-Cyl}
The homology operad  $H_{-\bul}(\Cyl, \bbK)$ of $\Cyl$
with the reversed grading is the operad $\calc$ of calculi.
\end{teo}
{\bf Proof.} Since the operad $\Cyl$ has two colors
algebras over  $H_{-\bul}(\Cyl, \bbK)$ are pairs of
graded vector spaces $(\cV, \cW)$\,.

The components $\Cyl^{\mc}(n,0)$ form the operad of little
discs. Thus, due to Theorem 1.2 in\footnote{We only need the
operations which survive in characteristic zero.} \cite{Cohen}
$\cV$ is a Gerstenhaber algebra.

The space $\Cyl^{\mc}(0,0)$ is empty,
the space $\Cyl^{\mc}(1,0)$ is a point and the space
$\Cyl^{\mc}(n,0)$ for $n > 1$ is homotopy equivalent to
the space $\Conf_n(\bbR^2)$ of configurations of
$n$ distinct points on $\bbR^2$\,.
The map
\begin{equation}
\label{h-equiv}
 \Cyl^{\mc}(n,0) \stackrel{\sim}{\to} \Conf_n(\bbR^2)
\end{equation}
which establishes the equivalence
associates with a configuration of
disjoint discs the configuration of
their centers.

The space $\Conf_2(\bbR^2)$ is, in turn, homotopy
equivalent to $S^1$\,. Thus the generator of
$H_0(S^1)$ represents the commutative product on $\cV$
and the generator of $H_1(S^1)$ represents the
bracket on $\cV$\,. This bracket has degree $-1$ because
we use the reversed grading on the homology groups.

The operations without inputs of color $\mc$
correspond to homology classes of
the space $\Cyl^{\ma}(0,1)$\,. Since this
space is homeomorphic to the circle $S^1$ we have
$$
H_{\bul}(\Cyl^{\ma}(0,1), \bbK) =
\begin{cases}
\bbK\,, \qquad {\rm if} ~~ \bul = 0,1 \,, \\
0 \,, \qquad {\rm otherwise}\,.
\end{cases}
$$
The generator of $H_{0}(\Cyl^{\ma}(0,1))$ corresponds
to the identity transformation of $\cW$ and
the generator $\de$ of $H_{1}(\Cyl^{\ma}(0,1))$ corresponds to
a unary operation on $\cW$\,. The operation
$\de$ has degree $-1$ because we use the reversed
grading on the homology groups.

The identity
$$
\de^2 = 0
$$
follows immediately from the fact that
$$
H_{2}(\Cyl^{\ma}(0,1), \bbK) = 0\,.
$$

Let us construct a homotopy equivalence between the
space $\Cyl^{\ma}(n,1)$ for $n\ge 1$ and the space
\begin{equation}
\label{Conf-col-n}
\Conf_n(\bbR^2 \setminus {\bf 0}) \times S^1\,,
\end{equation}
where $\Conf_n(\bbR^2 \setminus {\bf 0})$ is the
configuration space of $n$ distinct points on the
punctured plane $\bbR^2 \setminus {\bf 0}$\,.

First, we kill the rotation symmetry by fixing the position
of the point $t$ on the upper boundary $S^1\times \{c\}$ of the cylinder
$S^1 \times [a,c]$\,.
Second, we kill translation symmetry by setting $c=0$\,.

Next we assign to each configuration of discs on the
lateral surface $S^1\times (a,0)$ the configuration of centers
of the discs. In this way we get a homotopy equivalence
between the space $\Cyl^{\ma}(n,1)$ and the space
\begin{equation}
\label{Conf-Cyl}
\Conf_n(S^1 \times (a,0)) \times S^1\,,
\end{equation}
where the points on the factor $S^1$ correspond
to positions of the point $b$ on $S^1\times \{a\}$\,.

Finally, using the map
$$
\chi: S^1 \times (a,0) \to \bbR^2 \setminus {\bf 0}\,,
$$
$$
\chi(\vf, y) = \Big( \frac{y}{a-y}\cos(\vf)\,,\, \frac{y}{a-y}\sin(\vf)
\Big)
$$
we get the desired homotopy equivalence
\begin{equation}
\label{homot-equiv}
\Cyl^{\ma}(n,1) \simeq
\Conf_n(\bbR^2 \setminus {\bf 0}) \times S^1\,.
\end{equation}

Let us consider the homology groups of $\Cyl^{\ma}(1,1)$
in more details.

The space $\Cyl^{\ma}(1,1)$ is homotopy equivalent
to
$$
\bbR^2\setminus {\bf 0}\, \times \, S^1
$$
and the latter space is, in turn, homotopy equivalent
to $S^1 \times S^1$\,. Thus
\begin{equation}
\label{H-C-col-1}
\Cyl^{\ma}(1,1) \simeq S^1 \times S^1\,.
\end{equation}
Therefore
$$
H_{-\bul}(\Cyl^{\ma}(1,1), \bbK) =
\bbK\oplus \bs^{-1}\bbK^2 \oplus
 \bs^{-2}\bbK\,.
$$
To identify the cycles representing the homology classes
we parameterize the circle $S^1$ by the
angle variable $\vf\in [0, 2\pi]$ and the torus
$S^1\times S^1$ by the pair of angle variables
$\vf_1, \vf_2 \in [0, 2\pi]$\,.

The zeroth homology space $H_{0}(\Cyl^{\ma}(1,1), \bbK)$
is one-dimensional and its generator corresponds
to the contraction $i$ of elements
of $\cV$ with elements of $\cW$\,.
The second homology space $H_{2}(\Cyl^{\ma}(1,1), \bbK)$
is also one-dimensional. Its generator corresponds to
the composition $\de\, i\, \de$\,.

\begin{figure}
\psfrag{b}[]{$b$}
\psfrag{t}[]{$t$}
\begin{center}
\includegraphics[width=0.4\textwidth]{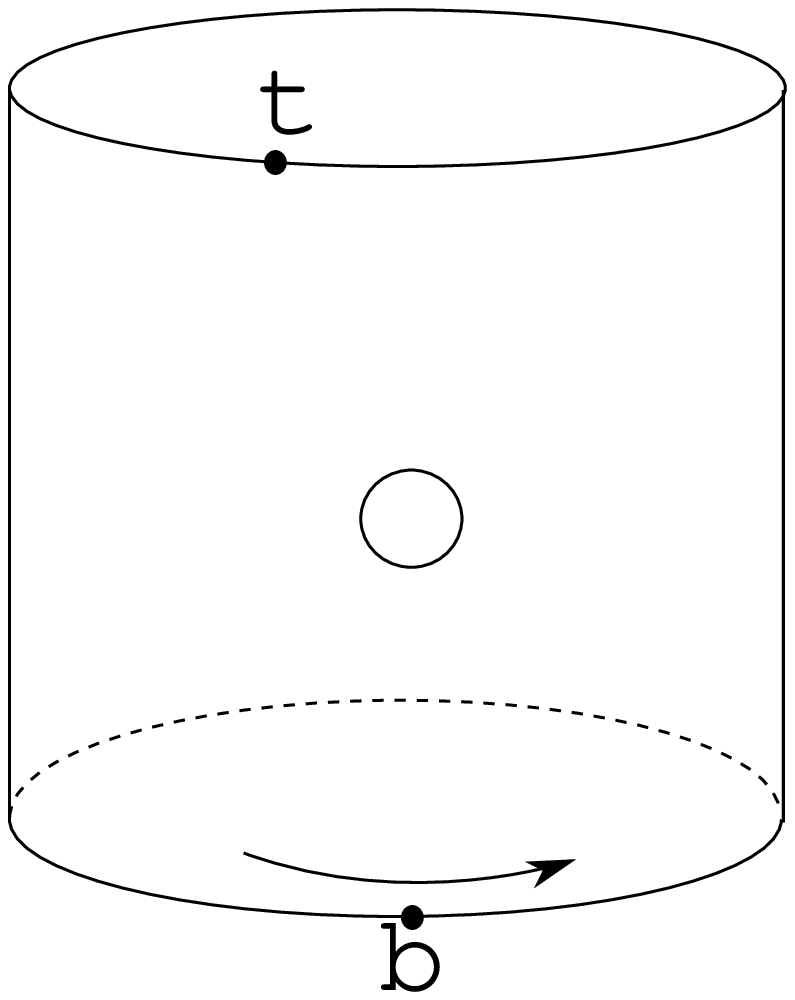}\\[0.5cm]
\parbox[t]{0.9\textwidth}{\caption{How to get the cycle in $\Cyl^{\ma}(1,1)$
representing the operation $i\,\de$}
\label{id}}
\end{center}
\end{figure}

The cycle
\begin{equation}
\label{cyc-id}
\vf \to (0, \vf) : S^1 \hookrightarrow S^1\times S^1
\end{equation}
represents the homology class corresponding to the
composition $i \, \de$\,. In order to get this cycle
in $\Cyl^{\ma}(1,1)$
we need to revolve the point $b$ on the lower boundary
about the vertical axis as it is shown on Figure \ref{id}\,.

\begin{figure}
\psfrag{b}[]{$b$}
\psfrag{t}[]{$t$}
\begin{center}
\includegraphics[width=0.4\textwidth]{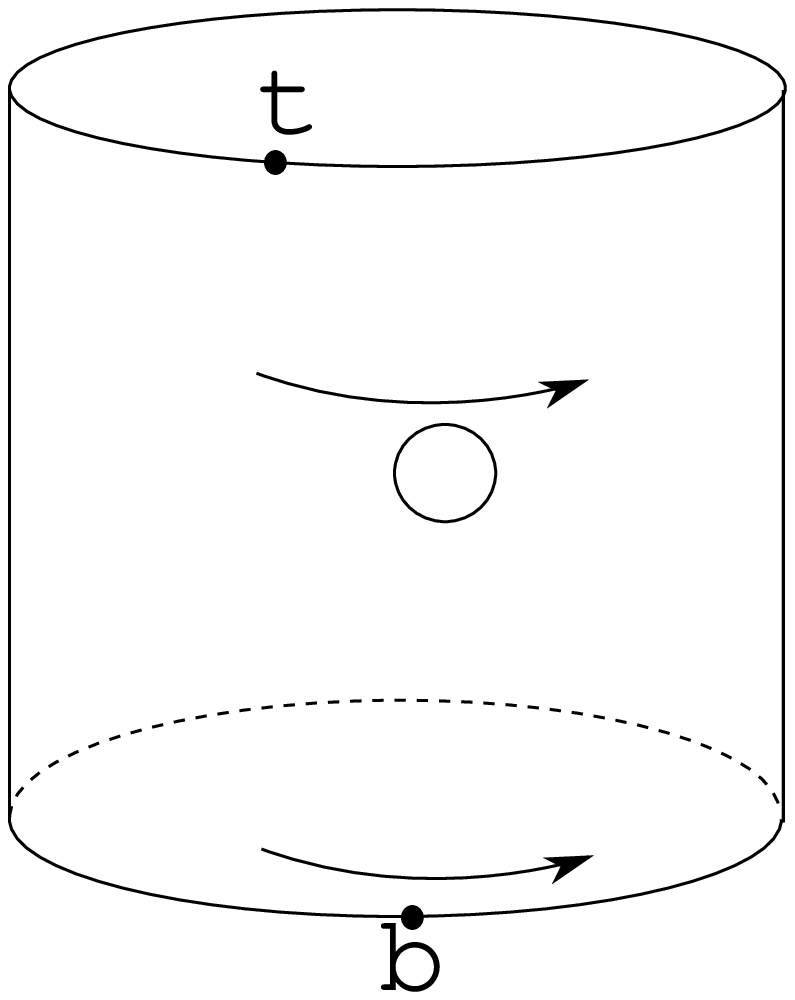}\\[0.5cm]
\parbox[t]{0.9\textwidth}{\caption{How to get the cycle in $\Cyl^{\ma}(1,1)$
representing the operation $\de\, i$}
\label{di}}
\end{center}
\end{figure}

The composition $\de i$ is, in turn, represented by
the diagonal
\begin{equation}
\label{cyc-di}
\vf \to (\vf, \vf) : S^1 \to S^1 \times S^1
\end{equation}
of the torus. To get this cycle we need to revolve simultaneously
the disc and the point $b$ about the vertical axis as it is shown
on Figure \ref{di}\,.

The homology classes $\de i$ and $i \de$
form a basis of $H_{1}(\Cyl^{\ma}(1,1), \bbK)$\,.

\begin{figure}
\psfrag{b}[]{$b$}
\psfrag{t}[]{$t$}
\begin{center}
\includegraphics[width=0.4\textwidth]{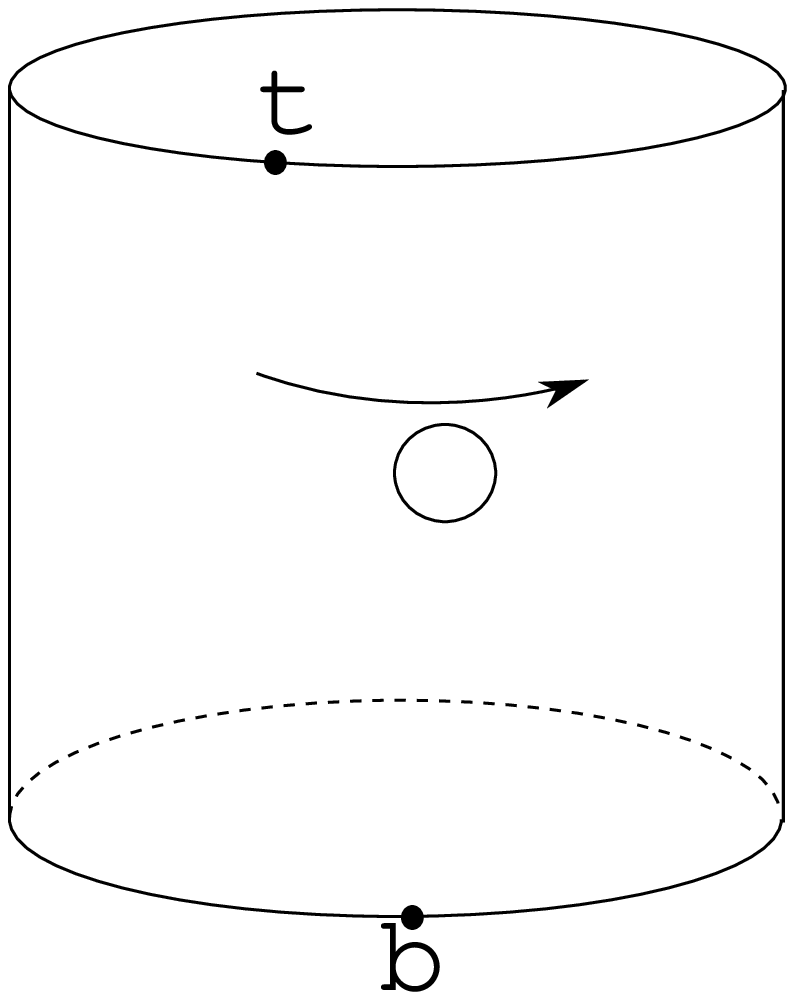}\\[0.5cm]
\parbox[t]{0.9\textwidth}{\caption{How to get the cycle in $\Cyl^{\ma}(1,1)$
representing the operation $l$}
\label{l}}
\end{center}
\end{figure}

We would like to remark that the homology class
represented by the cycle
\begin{equation}
\label{cyc-l}
\vf \to (\vf, 0) : S^1 \hookrightarrow S^1\times S^1
\end{equation}
equals to the combination
$$
\de \, i  - i \, \de\,.
$$
Indeed it is easy to see that the cycles
(\ref{cyc-id}), (\ref{cyc-di}), and (\ref{cyc-l})
form the boundary of the following $2$-simplex
in $S^1\times S^1$
$$
\{(\vf_1, \vf_2)\,\,|\,\, \vf_1 \le \vf_2 \} \subset S^1\times S^1\,.
$$
Thus the cycle (\ref{cyc-l}) represents the homology class
corresponding to the Lie derivative $l$\,.
To get this cycle in $\Cyl^{\ma}(1,1)$ we need to revolve
the disc about the vertical axis as it is shown on Figure
\ref{l}\,.

In general, for $n\ge 1$ the homology
of the space $\Cyl^{\ma}(n,1)$ can be computed
with the help of the homological version
of Lemma 6.2 from \cite{Cohen}\,.
Due to this lemma we have
\begin{equation}
\label{H-Conf}
H_{\bul}(\Conf_n(\bbR^2 \setminus {\bf 0}), \bbK) =
\bigotimes_{j=1}^{n}
H_{\bul} (\vee^j \, S^1, \bbK)\,.
\end{equation}
Using the homotopy equivalence (\ref{homot-equiv}) and
the K\"unneth formula, we deduce that
\begin{equation}
\label{H-Cyl-ma}
H_{\bul}(\Cyl^{\ma}(n,1), \bbK) =
\bigotimes_{j=1}^{n}
H_{\bul} (\vee^j \, S^1, \bbK) \,\, \otimes \,\,
H_{\bul}(S^1, \bbK)\,.
\end{equation}

Let us show that the operad
$H_{-\bul}(\Cyl, \bbK)$ is generated
by operations of $H_{-\bul}(\Cyl^{\mc}(2,0), \bbK)$,
$H_{-\bul}(\Cyl^{\ma}(0,1)$\,, and  $H_{-\bul}(\Cyl^{\ma}(1,1), \bbK)$\,.

Since the case of the operad of little discs was
already considered by F. Cohen \cite{Cohen}
we should only consider the operations of
$H_{-\bul}(\Cyl^{\ma}(n,1), \bbK)$\,.

Due to homotopy equivalence (\ref{homot-equiv})
the homology classes of $\Cyl^{\ma}(n,1)$
are of the forms
\begin{equation}
\label{form0}
\al \otimes 1\,,
\end{equation}
and
\begin{equation}
\label{form1}
\al \otimes \phi\,,
\end{equation}
where
$\al \in H_{-\bul}( \Conf_n(\bbR^2 \setminus {\bf 0}), \bbK ) $\,,
$1$ is the generator of $H_0(S^1, \bbK)$ and
$\phi$ is the generator of $H_0(S^1, \bbK)$\,.

It is obvious that homology classes of the form
(\ref{form1}) are obtained by composing the homology
classes of the form (\ref{form0}) with the generator
$\de $ of $ H_{1}(\Cyl^{\ma}(0,1), \bbK)$\,.

To analyze the homology classes
(\ref{form0}) we consider the Serre spectral
sequence corresponding to the fibration
(\ref{p}). Due to Proposition \ref{acts-triv}
the $E^2$ term of the sequence is
\begin{equation}
\label{E-2}
E^2_{p,q} = H_p(\Conf_{n-1}(\bbR^2 \setminus {\bf 0}), \bbK)
\otimes  H_q(F_n, \bbK)\,,
\end{equation}
where $F_n$ is the fiber (\ref{F-n}) of $p$\,.

Since $F_n$ is homotopy equivalent to the wedge
$\vee^n S^1$ of $n$ circles equation (\ref{H-Conf}) implies
that the vector spaces
$$
\bigoplus_{p}
E^2_{p, \bul-p}
$$
and $H_{\bul}(\Conf_{n}(\bbR^2 \setminus {\bf 0}), \bbK)$ have the
same dimension. Thus, using the fact  that spectral
sequence corresponding to the fibration
(\ref{p}) converges to $H_{\bul}(\Conf_{n}(\bbR^2 \setminus {\bf 0}), \bbK)$\,,
we deduce that this spectral sequence degenerates at $E_2$ and
\begin{equation}
\label{H-Conf1}
H_{\bul}(\Conf_{n}(\bbR^2 \setminus {\bf 0}), \bbK)
=\bigoplus_{p+q= \bul}
 H_p(\Conf_{n-1}(\bbR^2 \setminus {\bf 0}), \bbK) \otimes  H_q(F_n, \bbK)
\end{equation}
Using equation (\ref{H-F-n}) we reduce this expression further to
$$
H_{\bul}(\Conf_{n}(\bbR^2 \setminus {\bf 0}), \bbK) =
$$
\begin{equation}
\label{H-Conf11}
H_{\bul}(\Conf_{n-1}(\bbR^2 \setminus {\bf 0}), \bbK)\otimes  H_0(F_n, \bbK) \,\,
\oplus \,\,
 H_{\bul-1}(\Conf_{n-1}(\bbR^2 \setminus {\bf 0}), \bbK)
\otimes  H_1(F_n, \bbK)\,.
\end{equation}
Let $v_1, v_2, \dots, v_n\in \cV$ and $w\in \cW$ be the
arguments of an operation corresponding to a homology class
in (\ref{H-Conf11}).
It is not hard to see that
the generator of  $H_0(F_n, \bbK)= \bbK$ corresponds
to the contraction $i$ with $v_1$ and the generators of  $H_1(F_n, \bbK) = \bbK^n$
correspond\footnote{This picture is very reminiscent of the consideration
of Hochschild-Serre spectral sequence in the proof of Proposition 4.1 in
\cite{Dima-Disc}.}
to the brackets $[v_1, v_j]$ for $j \in \{2, 3, \dots, n\}$
and the Lie derivative $l_{v_1}$ along $v_1$\,.

Thus the homology classes of $H_{-\bul}(\Cyl^{\ma}(n,1), \bbK)$
are all produced by the operadic compositions of the classes
in $H_{-\bul}(\Cyl^{\ma}(n-1,1), \bbK)$\,,
$H_{-\bul}(\Cyl^{\ma}(1,1), \bbK)$\,,
and $H_{-\bul}(\Cyl^{\mc}(2,0), \bbK)$\,.

This inductive argument allows us to conclude that
the operad $H_{-\bul}(\Cyl, \bbK)$ is generated by
the classes
\begin{equation}
\label{gener-H-Cyl}
\begin{array}{c}
\begin{array}{ccc}
i\in H_{0}(\Cyl^{\ma}(1,1), \bbK)\,, &
l \in H_{1}(\Cyl^{\ma}(1,1), \bbK)\,, &
\de \in H_{1}(\Cyl^{\ma}(0,1), \bbK)\,,
\end{array}\\[0.5cm]
\begin{array}{cc}
\wedge \in H_0(\Cyl^{\mc}(2,0), \bbK)\,, &
[\,,\,] \in H_1(\Cyl^{\mc}(2,0), \bbK)\,.
\end{array}
\end{array}
\end{equation}

Using equation (\ref{free-calc-ma}) and the fact \cite{Klyach} that
$\dim \Lie(n) = (n-1)!$ it is not hard to show that
$$
\calc^{\ma}(n,1) \cong
\calc^{\ma}(n-1,1) \otimes (\bbK \oplus \bs^{-1} \bbK^n)
$$
as graded vector spaces.

On the other hand, equation (\ref{H-Cyl-ma}) gives
us the same isomorphism
$$
H_{-\bul}( \Cyl^{\ma}(n,1), \bbK) \cong
H_{-\bul}( \Cyl^{\ma}(n-1,1), \bbK)
 \otimes (\bbK \oplus \bs^{-1} \bbK^n)
$$
of graded vector spaces for $H_{-\bul}(\Cyl, \bbK)$\,.

Therefore, the dimensions of graded components
of $H_{-\bul}(\Cyl^{\ma}(n,1), \bbK)$
and $\calc^{\ma}(n,1)$ are equal.

Thus, in order the complete the proof of the theorem,
we need to show the operations (\ref{gener-H-Cyl})
satisfy the identities of calculus algebra
(see Definition \ref{calc}).

\begin{figure}
\psfrag{b}[]{$b$}
\psfrag{t}[]{$t$}
\psfrag{p}[]{$+$}
\psfrag{s}[]{$\sim$}
\begin{center}
\includegraphics[width=0.8\textwidth]{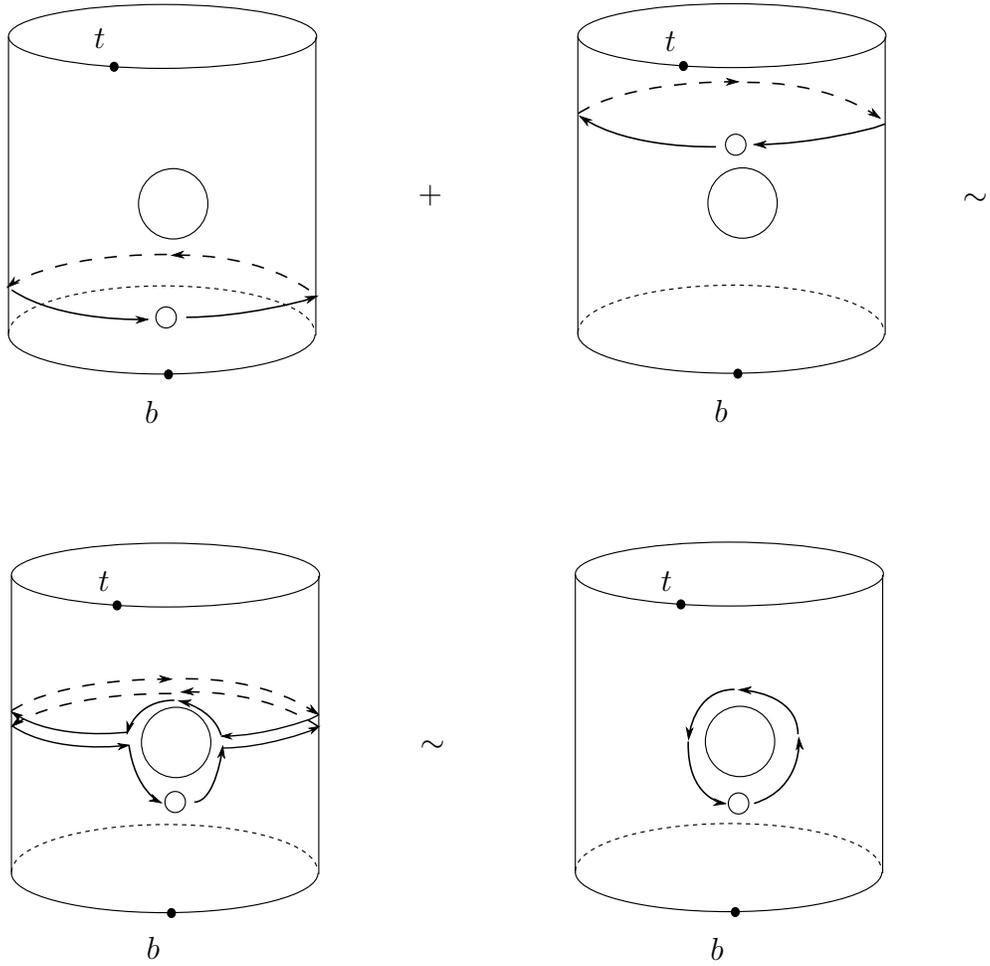}\\[0.5cm]
\parbox[t]{0.7\textwidth}{\caption{How to check identity (\ref{i-l-brack})}
\label{iden}}
\end{center}
\end{figure}

The identities of the Gerstenhaber algebra were already
checked in \cite{Cohen}. The identity (\ref{l-i-delta})
was checked above. On Figure \ref{iden}
we show how to check the identity
\begin{equation}
\label{i-l-brack}
[i, l] = i_{[\,,\,]}\,.
\end{equation}
The remaining identities can be checked in
the similar way.

Theorem \ref{H-Cyl} is proved. $\Box$

\subsection{Required results from \cite{K-Soi1}}
We will need
\begin{teo}[M. Kontsevich, Y. Soibelman, \cite{K-Soi1}]
\label{KS-teo}
The operad $\KS$ is quasi-isomorphic to the
operad of singular chains of the topological
operad $\Cyl$.
The homology operad $H_{-\bul}(\KS, \bbK)$ is generated
by the classes of the operations  $\cup$ (\ref{cup}),
$[\,,\,]_G$ (\ref{Gerst}), $I$ (\ref{I-P}), $L$ (\ref{L-Q})
and $B$ (\ref{B})\,.
\end{teo}
Furthermore, (See Proposition 11.3.3 on page 50
in \cite{K-Soi1})
\begin{pred}
\label{discs-formal}
The operad of singular chains of the topological
operad $\Cyl$ is formal.
\end{pred}
Combining these two statements with Theorem \ref{H-Cyl}
we get the following corollary.
\begin{cor}
\label{vot-ono}
The pair $(\nCbu(A), \nCbd(A))$ is a
homotopy calculus algebra. The operations
of this algebra structure are expressed
in terms of operations of $\cH$\,.
The induced
calculus structure on the pair
$(HH^{\bul}(A), HH_{\bul}(A))$ coincides
with the one in \cite{DGT}\,.
\end{cor}

\subsection{A useful property of the operad  $\KS$}
For the $2$-colored operad $\KS$ of chain complexes we
have the following
\begin{pred}
\label{low-degrees}
The elements of $\KS^{\mc}(k,0)$ have degrees
\begin{equation}
\label{deg1}
\deg  \ge 1 - k
\end{equation}
and the elements of $\KS^{\ma}(k,1)$ have
degrees
\begin{equation}
\label{deg11}
\deg \ge -1-k\,.
\end{equation}
\end{pred}
{\bf Proof.} Let us start with $\KS^{\mc}(k,0)$ for $k=1$\,.
Operations of $\KS^{\mc}(1,0)$ produce a
Hochschild cochain from a Hochschild cochain.
In order to prove that these operations have
nonnegative degrees we need to show that the
operations from $\cH(1,0)$ which lower the number the
number of arguments of the cochain do not contribute
to the realization of $\cH$\,.

Let $T$ be a tree representing an operation from
$\cH(1,0)$ and let $v_a$ be the vertex adjacent to the root.
If $v_a$ is marked then it is marked by the
only cochain $P$ which enters as an argument.
All other marked vertices are necessarily terminal and
they are reserved for the arguments of the cochain
produced by the operation.

In order to lower the number of arguments,
we need to insert the unit into the
cochain $P$.
The insertion of the unit is a degeneracy of the
simplicial structure on $\cH(1,0)$\,.
Hence all operations from $\cH(1,0)$ which
lower the degree of the cochain
do not contribute to the realization.

If $v_a$ is unmarked then, starting with the marked
vertex $v_P$ reserved for the cochain $P$,
we can form the proper maximal subtree with
$v_P$ being the vertex adjacent to the root.
In order to contribute to the realization
 the operation corresponding to this
subtree has to have a nonnegative degree.
Hence, so does the operation corresponding to
the whole tree.

We proved (\ref{deg1}) for $k=1$\,.

Let us take it as a base of the induction
and assume that (\ref{deg1}) is proved
for all $m < k $\,.

We consider a tree $T$ which
represents an operation from $\cH(k,0)$
and denote the vertex adjacent to the root of
$T$ by $v_a$\,.
Let us consider the case when the vertex $v_a$ is
marked. Say, $v_a$ is reserved for a cochain $P_1$ of
degree $q_1$\,.

Then the tree $T$ has exactly $q_1$ maximal proper subtrees
whose root vertex is $v_a$\,. We denote these subtrees
by $T_1, T_2, \dots, T_{q_1}$\,. The number $q_1$ splits
into the sum
\begin{equation}
\label{q1}
q_1 = p_{n} + p_{y}\,,
\end{equation}
where $p_{n}$ is the number of the subtrees
with no vertices reserved for cochains and $p_{y}$ is the number of
the subtrees in which at least one vertex is reserved
for a cochain.

Let $P$ denote the cochain produced by the operation
in question and let $r$ by the number of arguments
of this cochain. We will find an estimate for $r$
using the inductive assumption.

Every subtree with no vertices reserved for cochains
has to give at least one argument for the cochain $P$.
Otherwise, we have to insert
the unit as an argument of the cochain $P_1$\,.
In this case
the operation in question is a obtained from
another operation by degeneracy.
Therefore this operation would not contribute to
the realization of $\cH$\,.

If $T_j$ is a subtree with exactly $k_j$ vertices
reserved for cochains then, obviously, $k_j < k$\,. Hence,
applying the assumption of the induction, we get
that the number of arguments of $P$ coming from $T_j$ is
greater or equal
$$
(1- k_j) + q_j\,,
$$
where $q_j$ is the total degree of all cochains
entering as arguments of the operation corresponding
to the subtree $T_j$\,.

Thus the number $r$ of the arguments of the cochain $P$
produced by the operation in question can be estimated by
\begin{equation}
\label{r-ineq}
r \ge p_{n} + \sum_{j=1}^{p_y}(1- k_j + q_j)\,.
\end{equation}
This inequality can be rewritten as
$$
r - (\, p_n + p_y + \sum_{j=1}^{p_y} q_j \,) \ge
- \sum_{j=1}^{p_y} k_j\,.
$$
Due to equation (\ref{q1}) the sum
$$
p_n + p_y + \sum_{j=1}^{p_y} q_j
$$
is the total degree of all cochains entering
as arguments of the operation.
Furthermore, since the vertex $v_a$ is reserved
for one of the cochains
$$
\sum_{j=1}^{p_y} k_j = k-1
$$
and the desired inequality (\ref{deg1}) is
proved in this case.

Let us now consider the case when
the vertex $v_a$ is unmarked.

The valency of the vertex $v_a$ has to be at least $2$\,.
Otherwise the operation will have the empty set
of arguments. If the valency of this vertex is $2$ then we
remove it using the equivalence transformation.

Thus, without loss of generality, we may assume
that the vertex $v_a$ has at least
$2$ incoming edges. Let us denote by $s$
the number of these incoming edges and
let $T_1, \dots, T_s$ be the maximal proper subtrees of $T$
whose root vertex is $v_a$\,.

If the vertices reserved for cochains belong to
only one of these subtrees $T_i$ then excising
the subtrees $T_j$ for $j\neq i$ we get another operation
whose degree is less or equal the degree of the original
operation. Since in the new tree the unmarked vertex $v_a$ has
the valency $2$ we may remove this vertex by the corresponding
equivalence transformation.

If in this modified tree the vertex adjacent to the root is
marked then we deduce the desired inequality to the case
considered above.

Otherwise, we should only consider
the case when the vertex
$v_a$ is unmarked, its valency is at least $3$ and
each maximal proper subtree $T_j$ of $T$
with root vertex $v_a$ has at least one vertex reserved
for a cochain.

Let $s$ be the number of the maximal proper subtrees
of $T$ whose root vertex is $v_a$\,. Since the number of these
subtrees is greater or equal than $2$ therefore
 every subtree $T_j$ represents an operation
with the number of arguments $k_j < k$\,.

Let $\deg(T_j)$ be the degree of the operation
corresponding to the $j$-th subtree $T_j$\,.
Applying the assumption of the induction we
get the inequality
\begin{equation}
\label{m-j}
\deg(T_j) \ge 1 - k_j\,.
\end{equation}
It is clear that the degree $\deg(T)$ of the operation
corresponding to the tree $T$ is
the sum of degrees of the operations corresponding to
the subtrees $T_1, \dots, T_s$\,. Therefore
$$
\deg(T) \ge \sum_{i=1}^s (1- k_i)\,.
$$
On the other hand $\displaystyle \sum_{i=1}^s k_i = k$\,.
Therefore,
$$
\deg(T) \ge s - k
$$
and inequality (\ref{deg1}) follows from the
fact that $s \ge 2$\,.

To prove the second inequality (\ref{deg11}) we
denote by
 $\cH^{m_a}_{m_r}(k,1)$ the set of operations
producing a chain
\begin{equation}
\label{chain-r}
(b_0, b_1, \dots, b_{m_r})
\end{equation}
in $C_{-m_r}(A)$ from $k$ cochains and
a chain
\begin{equation}
\label{chain-a}
(c_0, c_1, \dots, c_{m_a})
\end{equation}
in $C_{-m_a}(A)$\,.

Let $F$ be a forest on the cylinder (\ref{Si})
representing an operation from $\cH^{m_a}_{m_r}(k,1)$
which contribute to the realization of $\cH(k,1)$\,.
Our purpose is to prove the inequality
\begin{equation}
\label{ineq-ch}
-m_r \ge -m_a + q -1 - k \,,
\end{equation}
where $q$ is the total degree of
all $k$ cochains of the operation.

This inequality is equivalent to
\begin{equation}
\label{ineq-ch1}
m_a \ge m_r + q -1 - k\,.
\end{equation}

By construction the forest $F$ has exactly $m_r$
trees. Let us denote these trees by $S_1, \dots, S_{m^n_r}$\,,
$T_1, \dots, T_{m^y_r}$ where the trees  $S_1, \dots, S_{m^n_r}$
have no vertices reserved for cochains and each tree
$T_i$ has at least one vertex reserved for a cochain.
Obviously,
\begin{equation}
\label{m-r}
m_r = m^n_r + m^y_r\,.
\end{equation}

The roots of the trees $S_1, \dots, S_{m^n_r}$\,,
$T_1, \dots, T_{m^y_r}$ are marked by components of
the chain (\ref{chain-r}).
If the root of the tree $S_i$ is marked by the
component $b_j$ of the for $j \neq 0$ then $S_i$ has to have
at least one terminal vertex marked by a component
of the chain (\ref{chain-a}). Otherwise we have to
insert the unit as the $j$-th component of (\ref{chain-r})
for $j\neq 0$\,. In this case the operation in question
is a composition of the another operation and a degeneracy.
Hence, this operation would not contribute to the
realization of $\cH(k,1)$\,.

Let us denote by $m_i$ the number of the
terminal vertices of the tree $T_i$ marked by
components of the chain (\ref{chain-a}).
If the tree $T_i$ has exactly $k_i$ vertices reserved
for the cochains then $T_i$ represents an operation from
$\cH(k_i,0)$\,. Furthermore,
if the operation corresponding to the forest $F$
contributes to the realization then so does
the operation corresponding to the tree $T_i$\,.
Hence, the number $m_i$ can be estimated using
the inequality (\ref{deg1})
$$
m_i \ge q_i +1 - k_i\,,
$$
where $q_i$ is the total degree of all cochains
of the operation corresponding to the tree
$T_i$\,.

Thus we get the following inequality for $m_a$
\begin{equation}
\label{ineq-m-a}
m_a \ge (m^n_r -1) + \sum_{i=1}^{m^y_r} (q_i +1 - k_i)\,,
\end{equation}
where the first term $(m^n_r -1)$ in the right hand side
comes from estimate of the number of
the marked terminal vertices of the
trees $S_1, \dots, S_{m^n_r}$\,.

Inequality (\ref{ineq-m-a}) can be rewritten as
\begin{equation}
\label{ineq-m-a1}
m_a \ge  m^n_r + m^y_r  + q - 1 - k\,,
\end{equation}
where $q$ is the total degree of
all $k$ cochains of the operation in question.

Due to equation (\ref{m-r}) inequality (\ref{ineq-m-a1})
coincides exactly with the desired inequality
(\ref{ineq-ch1}).

The proposition is proved. $\Box$

\section{The homotopy calculus on the pair
$(\nCbu(A), \nCbd(A))$\,.}
Since the operad $\La\Lie^+_{\de}$ is a suboperad of
$\calc$\,, Corollary \ref{vot-ono} implies that
the pair
$$
(\nCbu(A), \nCbd(A))
$$
carries a natural $Ho(\La\Lie^+_{\de})$-algebra
structure. In this section we show that
the homotopy calculus on $(\nCbu(A), \nCbd(A))$
can be modified in such a way that
its $Ho(\La\Lie^+_{\de})$-algebra part becomes
the $\La\Lie^+_{\de}$-algebra structure given by
the operations $[\,,\,]_G$ (\ref{Gerst}),
$L$ (\ref{L-Q}), and $B$ (\ref{B}).

Due to Proposition \ref{GJ}
a homotopy calculus structure on the pair
$(\nCbu(A), \nCbd(A))$ is a Maurer-Cartan element
\begin{equation}
\label{Q}
Q \in \Coder'( \bbF_{\bB}(\nCbu(A), \nCbd(A)) )\,.
\end{equation}

In other words, $Q$ is a degree $1$ coderivation of
the coalgebra $\bbF_{\bB}(\nCbu(A), \nCbd(A))$
satisfying the condition
$$
Q\Big|_{\nCbu(A) \oplus \nCbd(A)} = 0\,,
$$
and the equation
\begin{equation}
\label{MC-Q}
[\pa^{Bar} + \pa^{Hoch},Q] + Q \circ Q =0\,,
\end{equation}
where $\nCbu(A) \oplus \nCbd(A)$ is considered as subspace
of $\bbF_{\bB}(\nCbu(A), \nCbd(A))$ via coaugmentation
and $\pa^{Hoch}$
is the differential coming from the Hochschild (co)boundary
operator on $(\nCbu(A), \nCbd(A))$\,.

It is convenient to reserve a notation for
the Lie algebra
$\Coder'(\bbF_{\bB}(\nCbu(A), \nCbd(A)))$
\begin{equation}
\label{bbL}
\bbL = \Coder'(\bbF_{\bB}(\nCbu(A), \nCbd(A)))\,.
\end{equation}

Proposition \ref{coder-cofree} implies that
the coderivation $Q$ is uniquely determined by its
composition with the corestriction
$\rho : \bbF_{\bB}(\nCbu(A), \nCbd(A)) \to \nCbu(A)\oplus \nCbd(A) $
\begin{equation}
\label{q}
q = \rho \circ Q : \bbF_{\bB}(\nCbu(A), \nCbd(A)) \to
\nCbu(A)\oplus \nCbd(A) \,,
\end{equation}
while equation (\ref{MC-Q}) is equivalent to
\begin{equation}
\label{MC-Qq}
[\pa^{Hoch}, q] + q \circ \pa^{Bar}
+ q \circ Q  = 0\,.
\end{equation}

The coalgebra
$\bbF_{\bB}(\nCbu(A), \nCbd(A))$
is equipped with a natural increasing
filtration\footnote{See Equation (\ref{bB-01}).}
$$
\nCbu(A) \oplus \nCbd(A)[[u]] =
\cF^1\, \bbF_{\bB}(\nCbu(A), \nCbd(A))
\subset
\cF^2\, \bbF_{\bB}(\nCbu(A), \nCbd(A))
\subset \dots
$$
\begin{equation}
\label{filtr-calc}
\begin{array}{c}
\displaystyle
\cF^m\, \bbF_{\bB}(\nCbu(A), \nCbd(A))_{\mc}
=
\bigoplus_{k\le m}
\bB^{\mc}(k,0)\,\otimes_{S_{k}}\,(\nCbu(A))^{\otimes\, k}\,, \\[0.3cm]
\displaystyle
\cF^m\, \bbF_{\bB}(\nCbu(A), \nCbd(A))_{\ma}
=
\bigoplus_{k+1\le m}
\bB^{\ma}(k,1)\,\otimes_{S_{k}}\,(\nCbu(A))^{\otimes\, k} \otimes
\Cbd(A)\,,
\end{array}
\end{equation}
where $u$ is an auxiliary variable of degree $-2$\,.

Using (\ref{filtr-calc}) we endow the
Lie algebra (\ref{bbL})
with the following decreasing filtration
$$
\bbL
=\cF^0\, \bbL
\supset \cF^1\, \bbL
\supset
\cF^2\, \bbL
\supset \dots
$$
\begin{equation}
\label{filtr-coder}
\cF^m\, \bbL
= \{\, Y\in \bbL
~~|~~
Y\Big|_{\cF^m\, \bbF_{\bB}(\nCbu(A), \nCbd(A))} = 0
\, \}\,.
\end{equation}

Since
$$
\bbL =
\lim_{m} \bbL /
\cF^m\, \bbL
$$
the Lie subalgebra $\cF^1 \bbL^0$ is pronilpotent.

Therefore $\cF^1 \bbL^0$ integrates to a prounipotent
group
\begin{equation}
\label{bbG}
\bbG  = \exp(\cF^1 \bbL^0)
\end{equation}
which acts on the Maurer-Cartan elements of $\bbL$\,.
This action is defined by the formula:
\begin{equation}
\label{action}
\exp(Y)\, Q =
\exp([\,\,,Y])\, Q + f([\,\,,Y])\, [\pa^{Bar}
+ \pa^{Hoch}, Y]\,,
\end{equation}
where $f$ is the power series of the function
$$
f(x) = \frac{e^x - 1}{x}
$$
at the point $x=0$\,.

Let $Q$ be a  Maurer-Cartan element of the DGLA $\bbL$
which corresponds to the homotopy calculus
structure on the pair $(\nCbu(A), \nCbd(A))$
which comes from the action of the operad $\KS$\,.

For every $Y \in \cF^1 \bbL$ the homotopy
calculus on the pair $(\nCbu(A), \nCbd(A))$
corresponding to the Maurer-Cartan element $\exp(Y)\, Q$ is
quasi-isomorphic to the  homotopy calculus
corresponding the original Maurer-Cartan element $Q$\,.
Indeed, the desired $\Hocalc$-quasi-isomorphism
is expressed in terms of $Y$ as
\begin{equation}
\label{q-iso-Y}
\exp([\,\,,Y]) : (\bbF_{\bB}(\nCbu(A), \nCbd(A)), Q)
\to (\bbF_{\bB}(\nCbu(A), \nCbd(A)), \exp(Y)\, Q)\,.
\end{equation}

Thus we get a family of mutually quasi-isomorphic
homotopy calculus structures on $(\nCbu(A), \nCbd(A))$\,.
Let us denote this family by $\bS_{\calc}$\,.

We claim that
\begin{teo}
\label{valid-beer}
The family $\bS_{\calc}$ contains a
homotopy calculus structure whose
$Ho(\La\Lie^+_{\de})$-algebra part is
the $\La\Lie^+_{\de}$-algebra structure given by
the operations $[\,,\,]_G$ (\ref{Gerst}),
$L$ (\ref{L-Q}), and $B$ (\ref{B}).
\end{teo}
{\bf Proof.}
According to Proposition \ref{GJ} and equation (\ref{HoLaLie})
the $\Ho(\La\Lie^+_{\de})$ is given by a Maurer-Cartan element $M$ of the DGLA
\begin{equation}
\label{esche-odna}
\Coder'(\,\bbF_{\La(\Lie^+_{\de})^{\vee}}(\nCbu(A), \nCbd(A))\, )\,.
\end{equation}
Due to Proposition \ref{coder-cofree} this Maurer-Cartan element is, in turn,
uniquely determined by the composition with the corestriction
$\rho$ onto $\nCbu(A) \oplus \nCbd(A)$
\begin{equation}
\label{m}
m =
\rho \circ M : \bbF_{\La(\Lie^+_{\de})^{\vee}}(\nCbu(A), \nCbd(A))
\to \nCbu(A) \oplus \nCbd(A)\,.
\end{equation}

Finally, Proposition \ref{how-to-pull} shows that
the map $m$ is related to the map $q$ (\ref{q}) by the
equation
\begin{equation}
\label{m-q}
m = q \Big|_{\bbF_{\La(\Lie^+_{\de})^{\vee}}(\nCbu(A), \nCbd(A))}\,,
\end{equation}
where $\La(\Lie^+_{\de})^{\vee}$ is considered as a sub-cooperad
of $\bB = Bar(\calc)$ via the chain of embeddings:
$$
\La(\Lie^+_{\de})^{\vee} \hookrightarrow Bar(\La\Lie^+_{\de})
\hookrightarrow Bar(\calc)\,.
$$

The $\Ho(\La\Lie)$-algebra structure on $\nCbu(A)$
is induced by the $\Hoger$-algebra structure which,
in turn, comes from the action of the minimal operad
of Kontsevich and Soibelman \cite{K-Soi} on $\nCbu(A)$\,.
It was proved in \cite{BLT} (see Theorem 2) that this
$\Ho(\Lie)$-algebra is in fact a genuine Lie algebra
given by $[\,,\,]_G$\,.

Thus it remains to take care about the
operations involving a chain.

Due to (\ref{liede-vee1}) all the operations of
the $\Ho(\La\Lie^+_{\de})$-algebra involving a chain
are combined into a single degree $1$ map
$$
m^{\ma} =
$$
\begin{equation}
\label{m-ma}
m \Big|_{\bbF_{\La^2\cocomm^+}(\nCbu(A),\nCbd(A)[[u]])_{\ma}} :
\bbF_{\La^2\cocomm^+}(\nCbu(A),\nCbd(A)[[u]])
\end{equation}
$$
\to \nCbd(A)\,.
$$
such that
\begin{equation}
\label{m-00}
m^{\ma}\Big|_{\nCbd(A)} = 0\,.
\end{equation}
The latter equation follows from the fact that
the coderivation $M$ belongs to the DGLA (\ref{esche-odna}).

In other words, we have the infinite collection of operations:
\begin{equation}
\label{m-kn}
m^{\ma}_{k,n} : S^k(\nCbu(A))\otimes \nCbd(A) \to \nCbd(A)
\end{equation}
of degrees $1- 2 k- 2n$\,, where $S^k(\nCbu(A))$ is the
$k$-th component of the symmetric algebra $S(\nCbu(A))$ of $\nCbu(A)$\,.

Due to Proposition \ref{low-degrees} the operation
$m^{\ma}_{k,n}$ vanish if $k+2n > 2$\,. Furthermore, due to
equation (\ref{m-00}) we have $m^{\ma}_{0,0}=0$\,.
Thus we need to analyze only there operations:
one unary operation
\begin{equation}
\label{m-01}
m^{\ma}_{0,1} : \nCbd(A) \to \nCbd(A)
\end{equation}
of degree $-1$\,, one binary operation
\begin{equation}
\label{m-10}
m^{\ma}_{1,0} : \nCbu(A) \otimes \nCbd(A) \to \nCbd(A)
\end{equation}
of degree $-1$ and one ternary operation
\begin{equation}
\label{m-20}
m^{\ma}_{2,0} :  S^2(\nCbu(A))\otimes \nCbd(A) \to \nCbd(A)
\end{equation}
of degree $-3$\,.

Theorems \ref{H-Cyl} and \ref{KS-teo} imply that the operation (\ref{m-01})
differs from Connes' operator $B$ by an exact operation.
Namely,
$$
m^{\ma}_{0,1}(c) = B(c) + \pa^{Hoch}\, \beta(c) - \beta(\pa^{Hoch}\, c) \,,
$$
where $\beta$ is an operation in $\KS$
$$
\beta : \nCbd(A) \to \nCbd(A)
$$
of degree $-2$\,.

Using Proposition \ref{low-degrees} we deduce that $\beta$
is zero. Hence $m^{\ma}_{0,1}= B$\,.

\begin{figure}
\psfrag{b1}[]{$b_1$}
\psfrag{b2}[]{$b_2$}
\psfrag{P}[]{$P$}
\psfrag{c}[]{$c$}
\psfrag{eq}[]{$=$}
\psfrag{l}[]{$l$}
\psfrag{i}[]{$i$}
\begin{center}
\includegraphics[width=0.7\textwidth]{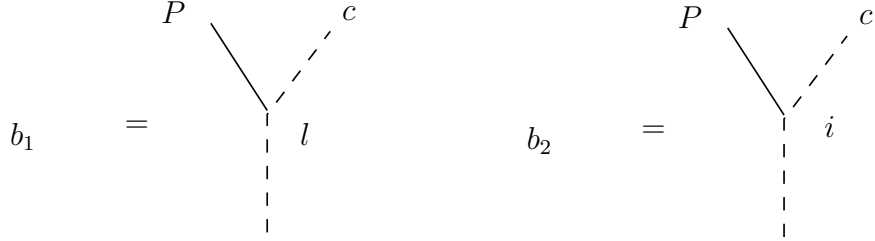}\\[1cm]
\parbox[c]{0.5\textwidth}{\caption{Elements $b_1$ and $b_2$}\label{b1b2}}
\end{center}
\end{figure}

Due to Proposition \ref{how-to-pull} the operation
(\ref{m-10}) is expressed in terms of $q$ (\ref{q}) as
\begin{equation}
\label{L-1}
m^{\ma}_{1,0}(P, c) = q(b_1)\,,
\end{equation}
where $P\in \nCbu(A)$\,, $c\in \nCbd(A)$\,, and
the element $b_1 \in \bbF_{\bB}(\nCbu(A), \nCbd(A))$
is depicted on Figure \ref{b1b2}.

Theorems \ref{H-Cyl} and \ref{KS-teo} imply that $m^{\ma}_{1,0}$
differs from the action $L$ of cochains
on chains by an exact operation. In other words,
\begin{equation}
\label{L-11}
m^{\ma}_{1,0} (P, c) = -(-1)^{|P|}L_P \,c + \pa^{Hoch} \psi(P, c)
- \psi(\pa^{Hoch}P, c) - (-1)^{|P|} \psi(P, \pa^{Hoch} c)\,,
\end{equation}
where $|P|$ is the degree of $P$ and
$$
\psi: \nCbu(A) \otimes \nCbd(A) \to \nCbd(A)
$$
is an operation in $\KS^{\ma}(1,1)$ of degree $-2$\,.

We remark that $\psi$ may be considered as a map
\begin{equation}
\label{psi}
\psi: \La(\Lie^+)^{\vee}(1,1)\otimes \nCbu(A) \otimes \nCbd(A)
\to \nCbd(A)
\end{equation}
of degree $0$\,. Our purpose is to extend $\psi$ ``by zero'' to
the whole vector space of the coalgebra
$\bbF_{\bB}(\nCbu(A) , \nCbd(A))$\,. This extension
depends on the choice of basis in $\calc^{\ma}(1,1)$\,.

\begin{figure}
\psfrag{b3}[]{$b_3$}
\psfrag{b4}[]{$b_4$}
\psfrag{P}[]{$P$}
\psfrag{c}[]{$c$}
\psfrag{eq}[]{$=$}
\psfrag{ld}[]{$l\,\delta$}
\psfrag{id}[]{$i\,\delta$}
\begin{center}
\includegraphics[width=0.7\textwidth]{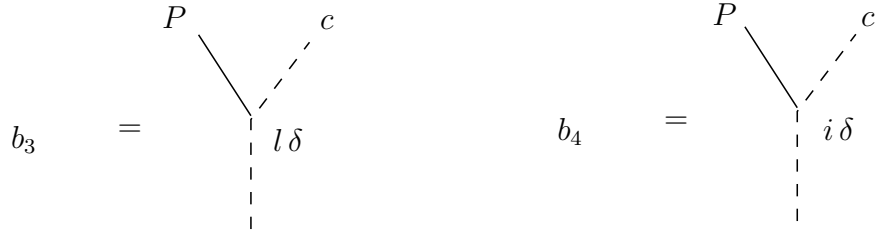}\\[1cm]
\parbox[c]{0.5\textwidth}{\caption{Elements $b_3$ and $b_4$}\label{b3b4}}
\end{center}
\end{figure}

\begin{figure}
\psfrag{bl}[]{$b_{\lambda}$}
\psfrag{d}[]{$\delta$}
\psfrag{P}[]{$P$}
\psfrag{c}[]{$c$}
\psfrag{eq}[]{$=$}
\psfrag{la}[]{$\lambda$}
\begin{center}
\includegraphics[height=0.7\textwidth]{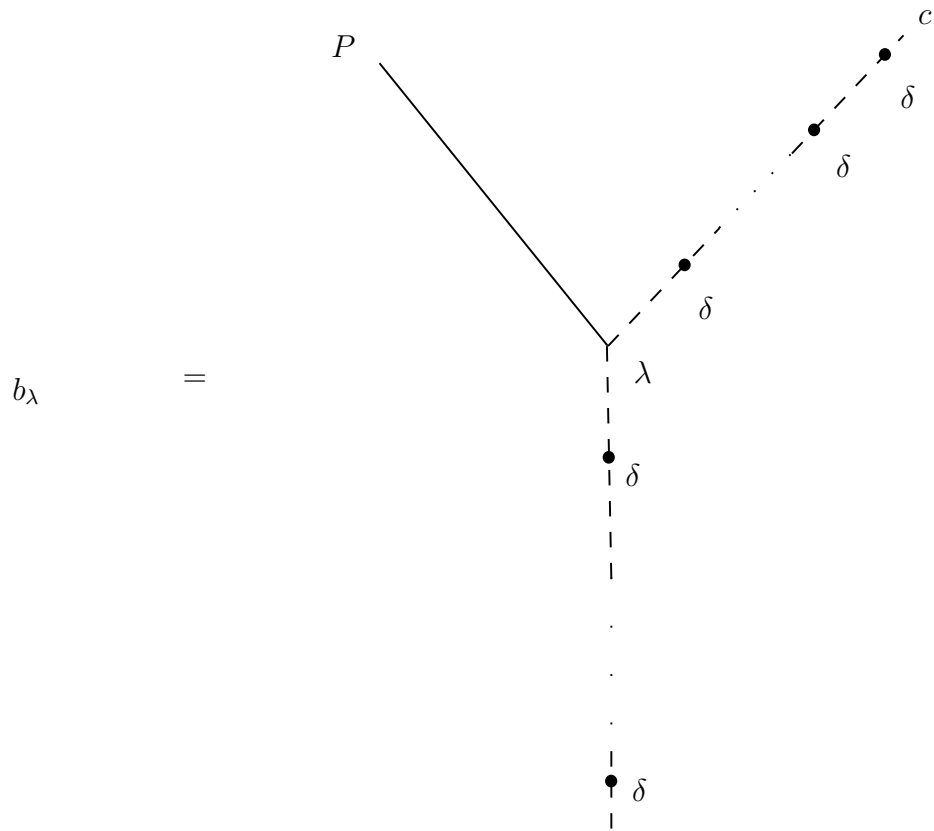}\\[1cm]
\parbox[c]{0.5\textwidth}{\caption{The number of $\de$'s is $\ge 1$}\label{bl}}
\end{center}
\end{figure}

We choose the basis
\begin{equation}
\label{basis}
\{ l, ~ i, ~ i\,\de, ~ l \, \de\}
\end{equation}
and extend $\psi$ to $ \bB^{\ma}(1,1)\otimes \nCbu(A) \otimes \nCbd(A)$
\begin{equation}
\label{ono}
\psi: \bB^{\ma}(1,1)\otimes \nCbu(A) \otimes \nCbd(A)
\to \nCbd(A)\,,
\end{equation}
as
$$
\psi(b_1) = \psi(P, c)\,, \qquad \psi(b_2) = \psi(b_3)= \psi(b_4) = 0\,,
$$
$$
\psi(b_{\la})=0\,,
$$
where the elements
$$
b_1, b_2, b_3, b_4, b_{\la}\in \bB^{\ma}(1,1)\otimes \nCbu(A) \otimes \nCbd(A)
$$
are depicted
on Figures \ref{b1b2}, \ref{b3b4}, \ref{bl}\,,
$\la$ is an arbitrary element of the basis (\ref{basis}),
and $P \in \nCbu(A)$\,, $c\in \nCbd(A)$\,.

Next we extend $\psi$ by zero to the whole vector space
of the coalgebra
\begin{equation}
\label{ona}
\begin{array}{c}
\bbF_{\bB}(\nCbu(A), \nCbd(A)) = \\[0.4cm]
\displaystyle
\bigoplus_n \bB^{\mc}(n,0)\otimes_{S_n} (\nCbu(A))^{\otimes\, n}
\,\, \oplus \,\,
\bigoplus_n \bB^{\ma}(n,1)\otimes_{S_n} (\nCbu(A))^{\otimes\, n}
\otimes \nCbd(A)\,.
\end{array}
\end{equation}

Then according to Proposition \ref{coder-cofree}
the equation
$$
\rho \circ \Psi =  \psi
$$
defines a derivation $\Psi$ of the
coalgebra (\ref{ona}). The derivation $\Psi$ has
degree $0$ since so does the map $\psi$\,.
Furthermore, it is obvious that
$$
\Psi \in \cF^1\, \bbL\,.
$$

Applying the element $\exp(-\Psi)$ of the group
$\bbG$ (\ref{bbG}) to the Maurer-Cartan element $Q$ (\ref{Q})
we adjust the component $m^{\ma}_{1,0}$ by killing
this additional exact term $\pa^{Hoch} \psi(P, c)
- \psi(\pa^{Hoch}P, c) -
(-1)^{|P|} \psi(P, \pa^{Hoch} c) $ in (\ref{L-11})\,.
In doing this we do not change the unary operations
because $\Psi\in \cF^1\, \bbL$\,.

Thus we are left with only one non-vanishing
operation (\ref{m-20})\,.

The Maurer-Cartan equation (\ref{MC-Qq}) implies that
$m^{\ma}_{2,0}$ should be closed with respect to the
differential $\pa^{Hoch}$\,.

Since the degree of $m^{\ma}_{2,0}$ is $-3$\,,
using Theorems \ref{H-Cyl} and \ref{KS-teo}\,,
we deduce that up to $\pa^{Hoch}$-exact
terms the operation $m^{\ma}_{2,0}$ is made of the
following ``building blocks'':
$$
L_{[P_1, P_2]_G}\, B\, c\,, \qquad L_{P_1} L_{P_2}\, B\, c\,,
\qquad  L_{P_2} L_{P_1}\, B\, c\,,
$$
where $P_1, P_2 \in \nCbu(A)$ and $c\in \nCbd(A)$\,.

Using the symmetry in the arguments $P_1$, $P_2$
and the compatibility with $\pa^{Hoch}$ it is not hard
to show that (up to $\pa^{Hoch}$-exact
terms) the most general expression for $m^{\ma}_{2,0}$
is
\begin{equation}
\label{L-2}
m^{\ma}_{2,0} (P_1, P_2, c) = (-1)^{|P_1|}\mu L_{[P_1, P_2]_G}\, B\, c\,,
\end{equation}
where $\mu\in \bbK$\,, $P_1, P_2\in \nCbu(A)$ and
$c\in \nCbd(A)$\,.

If necessary, we apply the above trick with the
action (\ref{action}) of the group (\ref{bbG})
to modify $Q$ (\ref{Q}) so that equation (\ref{L-2})
indeed holds.

To kill $m^{\ma}_{2,0}$ we will introduce the map
\begin{equation}
\label{y}
y: \bB^{\ma}(1,1)\otimes \nCbu(A) \otimes \nCbd(A) \to
\nCbd(A)\,.
\end{equation}
This map is defined by the equations
$$
y(b_1) = -\mu L_P\, B\, c\, \,, \qquad y(b_2) = y(b_3)= y(b_4) =0\,,
$$
$$
y(b_{\la})=0\,,
$$
where the elements
$$
b_1, b_2, b_3, b_4, b_{\la}\in \bB^{\ma}(1,1)\otimes \nCbu(A) \otimes \nCbd(A)
$$
are depicted on Figures \ref{b1b2}, \ref{b3b4}, \ref{bl},
$\la$ is an arbitrary element of the basis (\ref{basis}) in
$\calc^{\ma}(1,1)$ and
$$
P \in \nCbu(A)\,, \qquad  c\in \nCbd(A)\,.
$$
Then we extend $y$ by zero to the whole vector space
of the coalgebra (\ref{ona}). It is not hard to see that
$y$ is of degree $0$\,.

According to Proposition \ref{coder-cofree} the equation
$$
\rho \circ Y = y
$$
defines a degree $0$ coderivation $Y$ of the coalgebra
(\ref{ona})\,. Furthermore,
\begin{equation}
\label{Hoch-Y}
[\pa^{Hoch}, Y] = 0\,,
\end{equation}
and
\begin{equation}
\label{YF1}
Y \in \cF^1\, \bbL\,.
\end{equation}

Applying the element $\exp(Y)$ of the group $\bbG$ (\ref{bbG})
to the Maurer-Cartan element $Q$ (\ref{Q}) we get another homotopy
calculus structure
on $(\nCbu(A), \nCbd(A))$\,. This homotopy calculus structure is
determined by the Maurer-Cartan element $\exp(Y)\, Q$\,.

Let us denote by $\tm^{\ma}_{0,1}$\,, $\tm^{\ma}_{1,0}$\,,
and $\tm^{\ma}_{2,0}$ the operations (\ref{m-01}), (\ref{m-10}),
(\ref{m-20})
of the $Ho(\La\Lie^+_{\de})$-algebra corresponding to the
new homotopy calculus $\exp(Y)\, Q$\,.
Since $Y\in \cF^1 \bbL$ the unary operation cannot change.
Thus
$$
\tm^{\ma}_{0,1} = m^{\ma}_{0,1} = B\,.
$$
For the binary operation we have
\begin{equation}
\label{tm-10}
\tm^{\ma}_{1,0} (P,c) = \rho \circ \exp(Y)\, Q (b_1)\,,
\end{equation}
where $\rho$ is the projection from $\bbF_{\bB}(\nCbu(A), \nCbd(A))$
onto $ \nCbu(A) \oplus \nCbd(A) $ and
the element $b_1 \in \bbF_{\bB}(\nCbu(A), \nCbd(A))_{\ma}$
is depicted on Figure \ref{b1b2}.

Using (\ref{action}), (\ref{Hoch-Y}), and (\ref{YF1}) we
simplify equation (\ref{tm-10}) as follows
$$
\tm^{\ma}_{1,0} (P,c)=
m^{\ma}_{1,0} (P,c) + [(q\circ  Y - y \circ Q) + y \circ \pa^{Bar}]\, b_1\,.
$$
It is obvious that $\pa^{Bar} b_1 = 0$\,. Furthermore,
it is not hard to see that $q \circ Y(b_1)=0$ and
$y \circ Q (b_1)=0$\,. Thus the binary operation (\ref{m-10})
is also unchanged.

\begin{figure}
\psfrag{b}[]{$[\,,\,]$}
\psfrag{m}[]{$-$}
\psfrag{P1}[]{$P_1$}
\psfrag{P2}[]{$P_2$}
\psfrag{c}[]{$c$}
\psfrag{l}[]{$l$}
\psfrag{z1}[]{$(-1)^{|P_1|}$}
\psfrag{z2}[]{$-(-1)^{|P_1||P_2|}$}
\begin{center}
\includegraphics[height=0.7\textwidth]{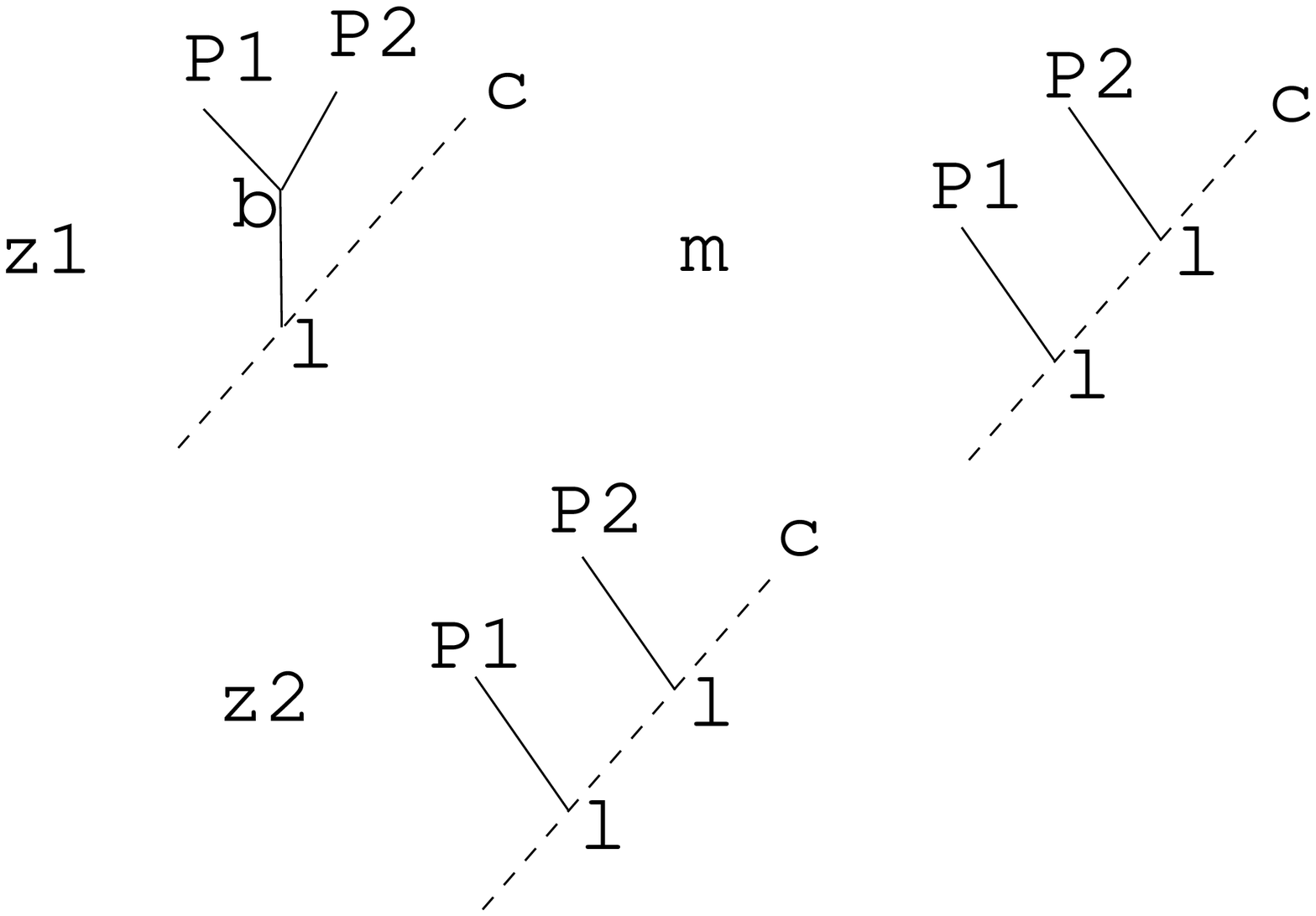}\\[1cm]
\parbox[c]{0.7\textwidth}{\caption{Here $P_1, P_2 \in \nCbu(A)$ and
$c\in \nCbd(A)$}\label{b}}
\end{center}
\end{figure}

For the ternary operation $\tm^{\ma}_{2,0}$ we have
$$
\tm^{\ma}_{2,0} (P_1,P_2,c) = \rho \circ \exp(Y)\, Q (b)=
$$
$$
= m^{\ma}_{2,0} (P_1,P_2,c) + (q\circ  Y - y \circ Q)(b)
$$
$$
+ y \circ \pa^{Bar} (b) - y \circ \pa^{Bar} \circ  Y(b) +
\frac{1}{2} y \circ Y \circ \pa^{Bar}(b)\,,
$$
where the element $b \in \bbF_{\bB}(\nCbu(A), \nCbd(A))_{\ma}$
is depicted on Figure \ref{b}
and, as above, we used (\ref{Hoch-Y}) and (\ref{YF1}).

It is obvious that $\pa^{Bar}(b)=0$ and it is not hard
to check that $\pa^{Bar}\circ Y(b)=0$\,. Furthermore a direct
computation shows that
$$
(q\circ  Y - y \circ Q)(b) = - (-1)^{|P_1|}\mu L_{[P_1, P_2]_G}\, B\, c\,.
$$
Thus
\begin{equation}
\label{tm-20}
\tm^{\ma}_{2,0} = 0
\end{equation}
and Theorem \ref{valid-beer} is proved. $\Box$

\section{Formality theorem}

\subsection{Enveloping algebra of a Gerstenhaber algebra}
Let $(\cV, \cW)$ be a pair of graded vector spaces.
For our purposes we will need to consider $\calc$-algebras
on $(\cV, \cW)$ with a fixed Gerstenhaber algebra structure
$(\cV, \wedge, [\,,\,])$ on $\cV$\,. For such $\calc$-algebras
we call $\cW$ {\it a $\calc$-module over $(\cV, \wedge, [\,,\,])$}.

The category of $\calc$-modules over $\cV$ is equivalent
to a category of ordinary modules over the enveloping
algebra \cite{TT} of the Gerstenhaber algebra $\cV$\,. In
this section we
recall the construction of this enveloping algebra
and describe its properties.

Let us start with the following definition:
\begin{defi}
For a Gerstenhaber algebra $(\cV, \wedge, [\,,\,])$
we define an associative algebra $\cY_0(\cV)$
which is generated by two sets of symbols
\begin{equation}
\label{generators}
l_{v},\,\, i_{v} \qquad  v \in \cV
\end{equation}
of degrees
\begin{equation}
\label{degrees}
|i_{v}| = |v|\,, \qquad |l_{v}|= |v|-1\,.
\end{equation}
These symbols are $\bbK$-linear in $v$ and
they are subject to the following relations
\begin{equation}
\label{Y-0relations}
\begin{array}{c}
\begin{array}{cc}
i_{v_1 \cdot v_2} = i_{v_1} i_{v_2}\,,  &
[i_{v_1}, l_{v_2}] = i_{[v_1, v_2]}\,,
\end{array} \\[0.7cm]
l_{v_1 \cdot v_2} = l_{v_1} i_{v_2} +
(-1)^{|v_1|} i_{v_1} l_{v_2}\,,
\qquad
[l_{v_1}, l_{v_2}] = l_{[v_1, v_2]}\,.
\end{array}
\end{equation}
\end{defi}

Furthermore,
\begin{defi}[\cite{TT}]
If  $(\cV, \cdot, [\,,\,])$ is a Gerstenhaber algebra
then the associative algebra $\cY(\cV)$
is generated by
 symbols (\ref{generators}) and
an element $\de$ of degree $-1$\,.
The symbols (\ref{generators}) are $\bbK$-linear
in $v$ and they are
subject to the following relations
\begin{equation}
\label{Y-relations}
\begin{array}{cc}
\de^2 = 0\,,    &  [\de, i_v] = l_v\,,  \\[0.3cm]
i_{v_1 \cdot v_2} = i_{v_1} i_{v_2}\,, &
[i_{v_1}, l_{v_2}] = i_{[v_1, v_2]}\,.
\end{array}
\end{equation}
\end{defi}
It is obvious that the category of $\calc$-modules
over $\cV$ is equivalent to the category of
ordinary modules over the associative algebra
$\cY(\cV)$\,.

Let us give the following natural definition
\begin{defi}
\label{regular}
A DG commutative algebra $\cV$ is called regular
if the module $\Om^1(\cV)$ of its K\"ahler differentials
is cofibrant.
\end{defi}
{\bf Remark.} Notice that if $\cV$ is a commutative
algebra concentrated in degree $0$
then the above condition of
regularity means that the $\cV$-module
$\Om^1(\cV)$ is projective.

We claim that
\begin{pred}
\label{Y-alm-exact}
Let $\cV$ be a DG Gerstenhaber algebra and
$\cR \to \cV$ be its cofibrant resolution.
If the corresponding DG commutative algebra
$\cV$ is regular then the induced map
\begin{equation}
\label{Y-map}
\cY(\cR) \to \cY(\cV)
\end{equation}
is a quasi-isomorphism of DG associative algebras.
\end{pred}
{\bf Proof.} Due to the obvious equality
of chain complexes
\begin{equation}
\label{cY-cY-0}
\cY(\cV) = \cY_0(\cV) \oplus
 \cY_0(\cV) \, \de
\end{equation}
it suffices to show that the map
\begin{equation}
\label{Y-0-map}
\cY_0(\cR) \to \cY_0(\cV)
\end{equation}
is a quasi-isomorphism.

For this purpose we introduce the following Lie-Rinehart
algebra structure \cite{Rinehart} on the pair $(\cV,\Om^1(\cV))$\,, where
$\Om^1(\cV)$ is the module of K\"ahler differentials of
the DG commutative algebra $\cV$\,. The Lie bracket $\{\,,\,\}$ on
$\Om^1(\cV)$ and the action $\ml$ of $\Om^1(\cV)$ on
$\cV$ are defined in terms of the Lie bracket on $\cV$
as follows
\begin{equation}
\label{brack-OmV}
\{ a_1\, d a_2, b_1 \, d b_2 \} =
(-1)^{|a_2|+1} a_1 [a_2, b_1]\, d b_2
+(-1)^{(|a_2|+1) |b_1|} a_1 b_1\, d([a_2, b_2]) +
\end{equation}
$$
(-1)^{(|a_1|+|a_2|+1)(|b_1|+|b_2|+1) + |b_2|}
\, b_1 [b_2, a_1] \, d a_2 \,,
$$
\begin{equation}
\label{action-OmV}
\ml_{a_1\, d a_2}(b) = (-1)^{|a_2|+1}a_1 [a_2, b]\,.
\end{equation}
The identities of a Gerstenhaber algebra imply that
equations (\ref{brack-OmV}) and (\ref{action-OmV})
indeed define a Lie-Rinehart algebra on the pair
$(\cV, \Om^1(\cV))$\,.

Next we remark that for every (DG) Gerstenhaber algebra $\cV$
the associative algebra $\cY_0(\cV)$ is nothing but
the enveloping algebra of the Lie-Rinehart algebra
$(\cV, \Om^1(\cV))$\,. Indeed, the required isomorphism
is defined on generators as
$$
a \mapsto i_a\,,
\qquad
d b \mapsto l_b\,, \qquad
a,b \in \cV\,.
$$

Then the PBW-filtration on $\cY_0(\cV)$ is
\begin{equation}
\label{Y-0-filtr}
\cV \cong \cF^0 \cY_0(\cV)
\subset  \cF^1 \cY_0(\cV) \subset
 \cF^2 \cY_0(\cV) \subset \dots
\end{equation}
where
$$
\cF^{k}\cY_0(\cV)
$$
is spanned by monomials in which the number
of symbols $l_{v}\,,$  $v \in \cV$ is less or
equal to $k$\,.

Since the DG commutative algebra $\cV$ is regular
we can apply the PBW-theorem \cite{Rinehart} to the Lie-Rinehart
algebra $(\cV, \Om^1(\cV))$\,. Using this theorem (see Theorem 3.1 in
\cite{Rinehart}) we conclude that
the associated graded algebra is isomorphic
to the symmetric algebra $S_{\cV}(\Om^1(\cV))$
\begin{equation}
\label{Assoc-Gr}
\bigoplus_k \cF^{k} \cY_0(\cV)/
 \cF^{k-1} \cY_0(\cV)
 \cong S_{\cV}(\Om^1(\cV))\,.
\end{equation}

Since $\cR$ is a cofibrant resolution of $\cV$ the
same argument with PBW theorem from \cite{Rinehart}
works for $\cY_0(\cR)$\,.
The map (\ref{Y-0-map}) is obviously compatible
with the filtrations (\ref{Y-0-filtr}) on $\cY_0(\cR)$
and $\cY_0(\cV)$\,. Furthermore, these filtrations
are cocomplete
$$
\cY_0(\cV) = \colim_{k} \, \cF^k \cY_0(\cV)\,,
\qquad
\cY_0(\cR) = \colim_{k} \, \cF^k \cY_0(\cR)\,.
$$
Hence, in order to prove that
the map (\ref{Y-0-map}) is a quasi-isomorphism,
we need to show that so is the map
\begin{equation}
\label{Om-map}
S_{\cR}(\Om^1(\cR))
\to S_{\cV}(\Om^1(\cV))\,.
\end{equation}
This statement follows from the regularity of
$\cV$\,. Thus the proposition is proved. $\Box$

\subsection{Sheaves of Hochschild (co)chains on
an algebraic variety}

Let $X$ be a smooth algebraic variety
over $\bbK$ with the structure sheaf $\cO_X$\,.
We denote by $\Vb_X$ the sheaf of
polyvector fields and by $\OmbX$ be the sheaf
of exterior forms with reversed grading.
$\cD_X$ denotes the sheaf of differential operators
on $X$ and $\cD_{\OmbX}$ denotes the sheaf of
differential operators on the sheaf of (graded) commutative
algebras $\OmbX$.

In the affine case $X = \Spec(A)$ we will use the short-hand
notation for the corresponding modules
of global sections $\Vb(A) = \G(X,\Vb_{X})$\,, $\Omb(A)=\G(X, \OmbX)$\,,
$\cD(A) = \G(X, \cD_{X})$\,, and finally
$\cD(\Omb(A))= \G(X, \cD_{\OmbX})$\,.

The pair $(\Vb_X, \OmbX)$ is a calculus algebra
with respect to the operations: the exterior product $\wedge$
on $\Vb_X$\,, the Schouten-Nijenhuis bracket $[\,,\,]_{SN}$
on $\Vb_X$\,, the contraction $\cI$ of a form with
a polyvector, the Lie derivative $\cL$ of a form along
a polyvector, and finally the de Rham differential
$d$ on the exterior forms.

Using the contraction $\cI$ we define the
natural $\cO_X$-linear pairing $\lan\,,\,\ran$ between
the sheaves $\Vb_X$ and $\OmbX$
$$
\lan \,,\, \ran: \Vb_X \otimes_{\cO_X} \OmbX \to \cO_X
$$
\begin{equation}
\label{pairing}
\lan \ga , \eta \ran =
\begin{cases}
\cI_{\ga} \eta\,, \qquad {\rm if} ~~ |\eta|= - |\ga| \,, \\
0 \,, \qquad {\rm otherwise}\,.
\end{cases}
\end{equation}
Here $\ga$ and $\eta$ are local sections of $\Vb_X$ and
$\OmbX$, respectively.

An appropriate version of Hochschild cochain complex for
the structure sheaf $\cO_X$ is the sheaf of polydifferential
operators \cite{Swan}, \cite{Y}. We will denote this sheaf
by $\Cbu(\cO_X)$\,. For example $C^0(\cO_X)$ is the structure
sheaf $\cO_X$ and $C^1(\cO_X)$ is the sheaf $\cD_X$ of differential
operators on $X$\,.

Let us also denote by $\nCbu(\cO_X)$
the sheaf of normalized polydifferential operators.
These are the polydifferential operators satisfying the
property
$$
P(\dots, 1, \dots) = 0\,.
$$

Similarly an
appropriate version of Hochschild chain complex for
the structure sheaf $\cO_X$ is the sheaf of polyjets
\cite{Tsygan}:
\begin{equation}
\label{polyjets}
\Cbd(\cO_X) = \cHom_{\cO_X}(C^{-\bul}(\cO_X), \cO_X)\,,
\end{equation}
where $\cHom$ denotes the sheaf-Hom and $\Cbu(\cO_X)$
is considered with its natural left $\cO_X$-module
structure. For example $C_{0}(\cO_X)$ is the structure sheaf
$\cO_X$ and $C_{-1}(\cO_X)$ is the sheaf of $\infty$-jets.

There are natural analogs of the degenerate Hochschild chains
$$
c = (c_0, \dots, 1, \dots)
$$
and these degenerate chains
form a subsheaf $C^{{\rm degen}}_{\bullet}(\cO_X)$
of $\Cbd(\cO_X)$\,. Furthermore the subsheaf $C^{{\rm degen}}_{\bullet}(\cO_X)$
is closed with respect to the Hochschild boundary operator $\pa^{Hoch}$\,.

We define the sheaf $\nCbd(\cO_X)$ of normalized Hochschild
chains as the quotient sheaf $\Cbd(\cO_X)/C^{{\rm degen}}_{\bullet}(\cO_X) $\,.
It is not hard to show that
\begin{equation}
\label{n-polyjets}
\nCbd(\cO_X) = \cHom_{\cO_X}(C^{-\bullet}_{norm}(\cO_X), \cO_X)\,.
\end{equation}

As well as for Hochschild complexes of an associative
algebra the inclusion
$$
\nCbu(\cO_X) \hookrightarrow \Cbu(\cO_X)
$$
and the projection
$$
\Cbd(\cO_X) \to \nCbd(\cO_X)
$$
are quasi-isomorphisms of complexes of sheaves.
Furthermore, the action of the Kontsevich-Soibelman
operad $\KS$ on the pair of sheaves $(\nCbu(\cO_X), \nCbd(\cO_X))$
is well-defined\footnote{Obvious extensions of the operations
on Hochschild chains to the operations on polyjets is discussed in
details in \cite{CDH}.}.

Thus the pair $(\nCbu(\cO_X), \nCbd(\cO_X))$ is a sheaf of
homotopy calculi.

Let us recall that the embedding
\begin{equation}
\label{HKR1}
\vr_{\mc}: \Vb_X \hookrightarrow \nCbu(\cO_X)
\end{equation}
is called the Hochschild-Kostant-Rosenberg map.
It is known \cite{HKR} that $\vr_{V}$ is
a quasi-isomorphism of complexes of sheaves where
the sheaf $\Vb_X$ is considered with the zero differential.

The corresponding quasi-isomorphism for Hochschild chains
\begin{equation}
\label{HKR11}
\vr_{\ma}: \nCbd(\cO_X) \to \OmbX
\end{equation}
is called the Connes-Hoch\-schild-Kos\-tant-Ro\-sen\-berg
map. This map is defined by the equation
\begin{equation}
\label{HKR111}
\lan \ga,
\vr_{\ma}(c) \ran = c(\vr_{\mc}(\ga))\,,
\end{equation}
where $c$ is a local section of $C_{-m}^{{\rm norm}}(\cO_X)$\,,
$\ga$ is a local section of $V^m_X$ and the pairing
$\lan\, , \, \ran$ is defined in (\ref{pairing}).

It is known \cite{Boris-book} that
the maps $\vr_{\mc}$ and $\vr_{\ma}$
are compatible with the operations of the Cartan calculus on the
pair $(\Vb_X, \OmbX)$ and the operations
$\cup$ (\ref{cup}), $[\,,\,]_G$ (\ref{Gerst}), $I$ (\ref{I-P}),
$L$ (\ref{L-Q}) and $B$ (\ref{B}) on the pair
$$
(\nCbu(\cO_X), \nCbd(\cO_X))
$$
{\it up to homotopy}.
We upgrade this observation to the following theorem.
\begin{teo}
\label{main}
If $X$ is a smooth algebraic variety
over a field $\bbK$ of characteristic zero
then the sheaf $(\nCbu(\cO_X), \nCbd(\cO_X))$ of
homotopy calculi
is quasi-isomorphic to the sheaf
$(\Vb_X, \Om^{-\bul}_X)$ of calculi.
\end{teo}
The proof of this theorem is given in Subsection \ref{proof}.

\subsection{Morita equivalence}
In this subsection we will show that the sheaf of algebras
$\cY(\Vb_X)$ is Morita equivalent to the sheaf $\cD_X[\md]/(\md^2)$
where $\md$ is an auxiliary variable of degree $-1$ commuting
with all the differential operators.

First we remark that $\cY_0(\Vb_X)$-module structure
on the sheaf $\OmbX$ gives us a natural map
\begin{equation}
\label{Y-D-Om-map}
\cY_0(\Vb_X) \to \cD_{\OmbX}
\end{equation}
between the sheaves of associative algebras.
We claim that
\begin{pred}
\label{Y-D-Om}
The map (\ref{Y-D-Om-map})
is an isomorphism of sheaves of associative algebras.
\end{pred}
{\bf Proof.}
We need to show that (\ref{Y-D-Om-map}) gives us an
isomorphism on stalks
$$
\cY_0(\Vb_X)_x \to (\cD_{\OmbX})_x
$$
for every point $x\in X$\,.

Thus, since $X$ is smooth, it suffices to show that
the map
\begin{equation}
\label{Y-D-A}
\cY_0(\Vb(A)) \to \cD_{\Omb(A)}
\end{equation}
is an isomorphism for every local regular (commutative)
algebra $A$ over $\bbK$\,.

It is easy to see that the associative algebra
$\cY_0(\Vb(A))$ is generated by symbols:
\begin{equation}
\label{symbols}
i_a,\,\, i_v, \,\, l_a,\,\, l_v\,,
\end{equation}
where $a\in A$ and $v\in V^1(A)$\,.

Under the map (\ref{Y-D-A}) the symbols go to
$$
i_a \to \cI_a\,, \qquad  i_v \to \cI_v\,,
\qquad l_a \to \cL_a \,, \qquad l_v \to \cL_v\,.
$$
Since the images of the symbols (\ref{symbols}) satisfy the
same relations therefore the map (\ref{Y-D-A}) is
injective.

To show that (\ref{Y-D-A}) is surjective we remark that
the algebra $\cD_{\Omb(A)}$ is generated by differential
operators of the form
$$
a\, \cdot\,, \qquad  d\, b\, \cdot\,, \qquad a,b\in A
$$
and derivations $\Der_{\bbK}(\Omb(A))$ of $\Omb(A)$\,.

Since $a\, \cdot\, = \cI_a$ and $d\,b \, \cdot = \cL_b$
it remains to show that every derivation $W\in \Der_{\bbK}(\Omb(A))$
belongs to the image of the map (\ref{Y-D-A}).

The regularity of $A$ implies that the $A$-modules
$\Om^{1}(A)$ and $\Omb(A)$ are free. More precisely, if
$x_1, \dots, x_n$ is a regular system of parameters in $A$ then
the $A$-module $\Om^1(A)$ is freely generated by the $1$-forms
\begin{equation}
\label{1forms}
d\,x^i\,, \qquad i =1, 2, \dots, n\,,
\end{equation}
and the $A$-module $\Omb(A)$ is freely generated by the forms
\begin{equation}
\label{forms}
d x^{i_1}\, d x^{i_2}\, \dots \, d x^{i_k}\,,
\qquad 1 \le i_1 < i_2 < \dots < i_k \le n\,.
\end{equation}

Dually the $A$-module $V^1(A) = \Der_{\bbK}(A)$ is freely generated by
\begin{equation}
\label{vectors}
e_1, e_2, \dots, e_n\,,
\end{equation}
where $e_i$ is a derivation of $A$ defined by the equation
$$
\cI_{e_i}(d x^j) = \de^j_i\,.
$$

Since the $A$-module $\Omb(A)$ is freely generated by
the forms (\ref{forms}) therefore every derivation
$W\in \Der_{\bbK}(\Omb(A))$ is
uniquely determined by its values on the elements of
$A$ and the $1$-forms (\ref{1forms})\,.
In general we have
$$
W(a) = \sum_{1 \le i_1 < i_2 < \dots < i_k \le n}
W_{i_1 \dots i_k}(a)\, d x^{i_1} d x^{i_2} \dots
d x^{i_k}
$$
and
$$
W(d x^i) = \sum_{1 \le i_1 < i_2 < \dots < i_k \le n} W^i_{i_1 \dots i_k}
d x^{i_1} d x^{i_2} \dots d x^{i_k}\,,
$$
where $ W_{i_1 \dots i_k}$ are derivations of $A$ over $\bbK$
and $W^i_{i_1 \dots i_k}\in A$\,.

Let $W_1$ be the following derivation of $\Omb(A)$:
$$
W_1 =  \sum_{1 \le i_1 < i_2 < \dots < i_k \le n}
d x^{i_1} d x^{i_2} \dots d x^{i_k} \, \cL_{W_{i_1 \dots i_k}}\,.
$$
It is obvious that the difference $W-W_1$ is $A$-linear.
Hence $W_2 = W - W_1$ is uniquely determined by its values on
$1$-forms (\ref{1forms}):
$$
W_2(d x^i) = \sum_{1 \le i_1 < i_2 < \dots < i_k \le n} \tW^i_{i_1 \dots i_k}
d x^{i_1} d x^{i_2} \dots d x^{i_k}\,,
$$
where $\tW^i_{i_1 \dots i_k}\in A$\,.
Thus the derivation $W$ can be rewritten as
$$
W=  \sum_{1 \le i_1 < i_2 < \dots < i_k \le n}
d x^{i_1} d x^{i_2} \dots d x^{i_k} \, \cdot\, \cL_{W_{i_1 \dots i_k}} +
 \sum_{1 \le i_1 < i_2 < \dots < i_k \le n}
d x^{i_1} d x^{i_2} \dots d x^{i_k} \, \cdot\, \cI_{\tW_{i_1 \dots i_k}}\,,
$$
where $\tW_{i_1 \dots i_k}$ is the derivation of $A$ defined by
$$
\tW_{i_1 \dots i_k} = \sum_i \tW^i_{i_1 \dots i_k} e_i\,.
$$
This shows that the map (\ref{Y-D-A}) is surjective.
Hence so is the map (\ref{Y-D-Om-map}). $\Box$

Using the pairing (\ref{pairing}) we introduce the following map
of sheaves of $\cO_X$-bimodules:
\begin{equation}
\label{map-r}
r\,:\, \OmbX \otimesO \cD_X \otimesO  \Vb_X \, \to \, \cD_{\OmbX}\,,
\end{equation}
$$
r(\eta, D, \ga)(\vf) = \eta \, D \lan \ga,\vf \ran\,,
$$
where $\eta$ and $\vf$ are local sections of $\OmbX$\,,
$D$ is a local section of $\cD_X$\,, $\ga$ is a local
section of $\Vb_X$\,, and $\cD_X$ is considered as
sheaf of bimodules over $\cO_X$\,.

We claim that
\begin{pred}
\label{map-r-pred}
The map $r$ in (\ref{map-r}) is an isomorphism of
sheaves of $\cO_X$-bimodules.
\end{pred}
{\bf Proof.}
Again, it suffices to prove that the map $r$ gives
an isomorphism on stalks
$$
(\OmbX \otimesO \cD_X \otimesO  \Vb_X)_x \to (\cD_{\OmbX})_x
$$
for every point $x\in X$\,.

Hence, due to smoothness of $X$\,, we need to show
that the map
\begin{equation}
\label{map-r-A}
r_A: \Omb(A) \otimes_A \cD(A) \otimes_A  \Vb(A) \to \cD_{\Omb(A)}
\end{equation}
is an isomorphism for every local regular (commutative)
algebra over $\bbK$\,.

Since $\Omb(A)$ and $\Vb(A)$ are free modules over $A$
the injectivity of $r_A$ follows easily from the injectivity
of the restriction
$$
r_A \Big|_{\cD(A)} : \cD(A) \to \cD_{\Omb(A)}\,.
$$

Due to Proposition \ref{Y-D-Om} the algebra
$\cD_{\Omb(A)}$ is generated by differential
operators of the form
\begin{equation}
\label{ops}
\cI_a = a\, \cdot\,, \qquad  \cL_a = d a\, \cdot\,, \qquad a\in A\,,
\end{equation}
\begin{equation}
\label{ops1}
\cI_v\,, \qquad \cL_v\,, \qquad v\in \Der_{\bbK}(A)\,.
\end{equation}
Thus, to show that $r_A$ is surjective it suffices to
prove that the operators (\ref{ops}) and (\ref{ops1})
belong to the image of $r_A$\,.

If we choose a regular system of parameters
$x_1, \dots, x_n$ in $A$ then
using the generators
(\ref{forms}) and (\ref{vectors}) we may rewrite
the operators  $\cI_a = a \,\cdot$\,,  $\cL_a = d a\, \cdot$
as follows:
$$
a \cdot \vf = \sum_{k; \, 1 \le i_1 < i_2 < \dots < i_k \le n}
d x^{i_1} d x^{i_2} \dots d x^{i_k}\, a \lan e_{i_k}\wedge \dots
\wedge e_{i_2} \wedge e_{i_1}, \vf \ran \,,
$$
and
$$
d a \cdot \vf = \sum_{k; \, 1 \le i_1 < i_2 < \dots < i_k \le n}
d a \, d x^{i_1} d x^{i_2} \dots d x^{i_k} \, \,
\lan e_{i_k} \wedge \dots \wedge e_{i_2} \wedge e_{i_1}, \vf \ran \,.
$$
Thus the operators (\ref{ops}) belong to the image of $r_A$\,.

Every derivation $v\in \Der_{\bbK}(A)$ can be uniquely
written as
$$
v = \sum_i v^i e_i\,,
$$
where $v^i \in A$\,.

Using this decomposition we rewrite the operators $\cI_v$ and  $\cL_v$
as
$$
\cI_v \vf =  \sum_{k; \, 1 \le i_1 < i_2 < \dots < i_k \le n}
\sum_{s=1}^k (-1)^{s-1}
v^{i_s} \, d x^{i_1} \dots \widehat{d x^{i_s}} \dots d x^{i_k} \, \,
\lan e_{i_k} \wedge \dots \wedge e_{i_2} \wedge e_{i_1}, \vf\ran\,,
$$
\begin{equation}
\label{L-v}
\cL_v \vf =  \sum_{k; \, 1 \le i_1 < i_2 < \dots < i_k \le n}
 d x^{i_1} d x^{i_2}  \dots  d x^{i_k}\, v \,
\lan e_{i_k} \wedge \dots \wedge e_{i_2} \wedge e_{i_1}, \vf \ran  +
\end{equation}
$$
\sum_{k; \, 1 \le i_1 < i_2 < \dots < i_{k+1} \le n}
\sum_{t,s=1}^{k+1} (-1)^{s-t} e_{i_t}(v^{i_s})\,
d x^{i_1} \dots \widehat{d x^{i_s}} \dots
d x^{i_{k+1}} \,\,
\lan e_{i_{k+1}} \wedge \dots \wedge
\widehat{e_{i_t}}\wedge \dots \wedge e_{i_1}, \vf \ran\,,
$$
where the symbol $\widehat{~}$ over $d x^{i}$ (resp. $e_i$)
means that the $1$-form $d x^{i}$ (resp. the vector $e_i$) is
omitted. The vector $v$ in the right hand side
of (\ref{L-v}) is considered as a differential operator on $A$\,.

Thus the operators (\ref{ops1}) also belong to the image of
$r_A$\,. This concludes the proof of the proposition. $\Box$

We remark that the sheaf
\begin{equation}
\label{cP}
\cP = \OmbX \otimesO \cD_X
\end{equation}
has the natural left $\cD_{\OmbX}$-module
structure and the natural right $\cD_X$-module structure.
Similarly, the sheaf
\begin{equation}
\label{cQ}
\cQ = \cD_X \otimesO \Vb_X
\end{equation}
has the natural left $\cD_X$-module structure
and a natural right $\cD_{\OmbX}$-module
structure.

Proposition \ref{map-r-pred} implies that
\begin{equation}
\label{Morita-PQ}
\cP \otimes_{\cD_X} \cQ  \cong \cD_{\OmbX}\,.
\end{equation}
Furthermore, it is obvious that
\begin{equation}
\label{Morita-QP}
\cQ \otimes_{\cD_{\OmbX}} \cP  \cong \cD_{X}\,.
\end{equation}

Thus we arrive at the following statement
\begin{pred}
\label{Morita0}
The sheaves $\cP$ (\ref{cP}) and $\cQ$ (\ref{cQ})
establish a Morita equivalence between the sheaves
of associative algebras $\cD_{\OmbX}$ and $\cD_X$\,. $\Box$
\end{pred}

Combining this statement with Proposition \ref{Y-D-Om}
we conclude that the sheaves $\cY_0(\Vb_X)$ and
$\cD_{X}$ are Morita equivalent.

In order to get the sheaf of algebras $\cY(\Vb_X)$
from $\cY_0(\Vb_X)$ we need to tensor $\cY_0(\Vb_X)$
with the constant sheaf
$$
\bbK[\de]/(\de^2)\,, \qquad \deg(\de) = -1
$$
and impose the equation
\begin{equation}
\label{de-hde}
[\de, P] = \hde(P)\,,
\end{equation}
where $\hde$ is the derivation of $\cY_0(\Vb_X)$ defined by
\begin{equation}
\label{hat-delta}
\hde(i_{\ga}) = l_{\ga}\,, \qquad \hde(l_{\ga}) = 0\,,
\end{equation}
$P$ is a local section of $\cY_0(\Vb_X)$ and $\ga$ is a
local section of $\Vb_X$\,.

On the other hand, for every local section $P$ of $\cY_0(\Vb_X)$
we have
$$
\hde(P) = [d, P]\,,
$$
where $d$ is the de Rham differential.

The de Rham differential $d$ is a global section of the
sheaf $\cD_{\OmbX}$. Hence, due to Proposition \ref{Y-D-Om},
$d$ is a global section of $\cY_0(\Vb_X)$\,.

Thus, switching from $\de$ to
\begin{equation}
\label{md}
\md = \de - d
\end{equation}
we get the following isomorphism of the sheaves
of algebras
\begin{equation}
\label{Y-Y-0}
\cY(\Vb_X) \cong \cY_0(\Vb_X)[\md]/(\md^2)\,,
\end{equation}
where $\md$ has degree $-1$ and
$$
[\md, P] =0
$$
for every local section $P$ of $\cY_0(\Vb_X)$\,.

Combining this observation with Proposition \ref{Morita0}
we arrive at the following statement
\begin{pred}
\label{Morita}
Let $\cP$ and $\cQ$ be the sheaves defined in (\ref{cP})
and (\ref{cQ}), respectively. The sheaves
$$
\cP[\md]/(\md^2)\,, \qquad \cQ[\md]/(\md^2)
$$
establish a Morita equivalence between the
sheaf of associative algebras $\cY(\Vb_X)$
and the sheaf
$$
\cD_X[\md]/(\md^2)\,,
$$
where $\md$ has degree $-1$ and
$$
[\md, D] =0
$$
for every local section $D$ of $\cD_X$\,. $\Box$
\end{pred}

\subsection{Proof of Theorem \ref{main}}
\label{proof}

Let us recall that, due to Proposition \ref{GJ},
a homotopy calculus structure on the pair
$$
(\nCbu(\cO_X), \nCbd(\cO_X))
$$
is a Maurer-Cartan element $Q$ of the DGLA
\begin{equation}
\label{coder-calc1}
\Coder'( \bbF_{\bB}(\nCbu(\cO_X), \nCbd(\cO_X)) )\,,
\end{equation}
where $\bB$ is as above the cooperad $Bar(\calc)$ and
the DGLA (\ref{coder-calc1}) consists of the
coderivations of $\bbF_{\bB}(\nCbu(\cO_X), \nCbd(\cO_X))$
satisfying the condition
$$
Q \Big|_{\nCbu(\cO_X) \oplus \nCbd(\cO_X)} =0\,.
$$

The codifferential on
the sheaf of coalgebras $\bbF_{\bB}(\nCbu(\cO_X), \nCbd(\cO_X))$
is the sum
\begin{equation}
\label{pa}
\pa = \pa^{Bar} + \pa^{Hoch}\,,
\end{equation}
where $\pa^{Bar}$ comes from the bar differential on
$\bB= Bar(\calc)$ and $\pa^{Hoch}$ comes from the Hochschild
(co)boundary operators on $\nCbu(\cO_X)$ and $\nCbd(\cO_X)$\,.

Using the Maurer-Cartan element $Q$ we shift the codifferential
$\pa$ by $[Q, \,]$ and get the new codifferential on
$\bbF_{\bB}(\nCbu(\cO_X), \nCbd(\cO_X))$:
\begin{equation}
\label{pa-Q}
\pa^Q = \pa^{Bar} + \pa^{Hoch} + [Q,\,]\,.
\end{equation}
Let us denote the resulting sheaf of DG $\bB$-coalgebras by
$C_Q$:
\begin{equation}
\label{C-Q}
C_Q = \Big( \bbF_{\bB}(\nCbu(\cO_X), \nCbd(\cO_X)),
\pa^Q \Big)\,.
\end{equation}

We use $C_Q$ to get the canonical free resolution\footnote{
The construction of this free resolution is known in topology
as {\it the rectification} \cite{BVogt}. We describe this construction
in more details in \cite{BLT} (See Proposition 3 therein).} $\cR$
of the sheaf $(\nCbu(\cO_X), \nCbd(\cO_X))$ of homotopy calculi.
As a sheaf of calculi,
\begin{equation}
\label{cR}
\cR = \bbF_{\calc}(C_Q)
\end{equation}
and the differential on $\cR$ is
\begin{equation}
\label{d-cR}
\pa^{\cR} = \pa^{\tw} + \pa^Q\,,
\end{equation}
where $\pa^{Q}$ comes from the differential on
$C_Q$ and $\pa^{\tw}$ is defined using the twisting cochain
between operad $\calc$ and cooperad $\bB=Bar(\calc)$\,.
(See Section 2.3 in \cite{GJ} on twisting cochain and
the construction of the differential $\pa^{\tw}$ for
algebras over an abstract operad.)

The sheaf of DG $\calc$-algebras $\cR$ splits according
to the colors $(\mc, \ma)$ as
$$
\cR= \cR_{\mc} \oplus \cR_{\ma}\,,
$$
where
\begin{equation}
\label{cR-ma}
\cR_{\ma} =
\bbF_{\calc}\Big(\bbF_{Bar(\calc)}(\nCbu(\cO_X), \nCbd(\cO_X)) \Big)_{\ma}\,,
\end{equation}
and
\begin{equation}
\label{cR-mc}
\cR_{\mc} =
\bbF_{\Ger}\circ \bbF_{Bar(\Ger)}(\nCbu(\cO_X))\,.
\end{equation}
Thus the sheaf $\cR_{\mc}$ with the differential
$$
\pa^{\cR} \Big|_{\cR_{\mc}}
$$
is a free resolution of the sheaf $\nCbu(\cO_X)$ of
homotopy Gerstenhaber algebras.
This resolution can be simplified. More precisely,
we may consider the subsheaf
\begin{equation}
\label{Ger-simple0}
\bbF_{\Ger}\circ \bbF_{\Ger^{\vee}}(\nCbu(\cO_X))
\subset \cR_{\mc}
\end{equation}
with the differential obtained by restricting the
one on (\ref{cR-mc}). Then, using
the fact that the inclusion
$\io_{\Ger}$ (\ref{Ger-Koszul}) is a quasi-isomorphism
of cooperads one can show that the inclusion
\begin{equation}
\label{Ger-simple-map}
\bbF_{\Ger}(\io_{\Ger})\, : \,
\bbF_{\Ger}\circ \bbF_{\Ger^{\vee}}(\nCbu(\cO_X))\,
\stackrel{\sim}{\hookrightarrow} \,
\bbF_{\Ger}\circ \bbF_{Bar(\Ger)}(\nCbu(\cO_X))
\end{equation}
is a quasi-isomorphism of sheaves of DG Gerstenhaber algebras.

Thus the sheaf (\ref{Ger-simple0}) is also a free resolution
of the sheaf $\nCbu(\cO_X)$ of homotopy Gerstenhaber algebras.

We denote the differential on the sheaf
(\ref{Ger-simple0}) by $\pa^{\cR}_{\mc}$
and reserve the notation $\cR(\nCbu(\cO_X))$
for this resolution
\begin{equation}
\label{Ger-simple}
\cR(\nCbu(\cO_X)) =
(\, \bbF_{\Ger}\circ \bbF_{\Ger^{\vee}}(\nCbu(\cO_X))\, , \,
\pa^{\cR}_{\mc} \,)\,.
\end{equation}

The quasi-isomorphism (\ref{Ger-simple-map}) provides the sheaf
$\cR_{\ma}$ (\ref{cR-ma}) with a (DG) $\calc$-module structure
over the sheaf $\cR(\nCbu(\cO_X))$ (\ref{Ger-simple}).
Thus, in order to prove Theorem \ref{main}, we need to
show that the sheaf $(\cR(\nCbu(\cO_X)), \cR_{\ma})$ of
calculi is quasi-isomorphic to the sheaf
$(\Vb_X, \OmbX)$\,.

In paper \cite{BLT} we constructed a chain of quasi-isomorphisms
of sheaves of DG Gerstenhaber algebras which connects the sheaf
$\cR(\nCbu(\cO_X))$ to the sheaf $\Vb_X$\,.

In this construction we use the sheaf of $\Ger^{\vee}$-coalgebras:
\begin{equation}
\label{DT}
\Xi_X = \bbF_{\La^{2}\cocomm} \circ
\bbF_{\La\coLie^+} (\cO_X, V^1_X)\,,
\end{equation}
where $\coLie^+$ is the
cooperad which governs pairs ``a Lie coalgebra $+$ its comodule.''

We also use the canonical free resolution
\begin{equation}
\label{R-V-X}
(\bbF_{\Ger}\circ \bbF_{\Ger^{\vee}}(\Vb_X),
 \pa^{\cR}_V)
\end{equation}
of the sheaf of Gerstenhaber algebras $\Vb_X$\,.
Here the differential $\pa^{\cR}_V$ on the sheaf (\ref{R-V-X}) comes
from the twisting cochain \cite{GJ} of the pair $(\Ger, \Ger^{\vee})$\,.

It is obvious that $\Xi_X$ is a subsheaf of
$$
\bbF_{\Ger^{\vee}}(\nCbu(\cO_X))\,.
$$
and a subsheaf of
$$
\bbF_{\Ger^{\vee}}(\Vb_X)\,.
$$
Due to this observation we have two inclusions of
the sheaves of free Gerstenhaber algebras
\begin{equation}
\label{incl1}
\si_1 :\bbF_{\Ger}(\Xi_X) \hookrightarrow
\cR(\nCbu(\cO_X))\,,
\end{equation}
and
\begin{equation}
\label{incl11}
\si_2 :\bbF_{\Ger}(\Xi_X) \hookrightarrow
\bbF_{\Ger}\circ \bbF_{\Ger^{\vee}}(\Vb_X)\,.
\end{equation}

It was shown in \cite{BLT} that the sheaf $\bbF_{\Ger}(\Xi_X)$
is closed both with respect to the differential  $\pa^{\cR}_{\mc}$
on $\cR(\nCbu(\cO_X))$ and the differential $\pa^{\cR}_V$ on the
sheaf (\ref{R-V-X}). Furthermore, the restriction of the differential
$\pa^{\cR}_{\mc}$ to $\bbF_{\Ger}(\Xi_X)$ coincides with the restriction
of the differential $\pa^{\cR}_{V}$\,. In other words, the sheaf
of free Gerstenhaber algebras  $\bbF_{\Ger}(\Xi_X)$
is equipped with a canonical differential.

Composing $\si_2$ (\ref{incl11}) with the
quasi-isomorphism
$$
\bbF_{\Ger}\circ \bbF_{\Ger^{\vee}}(\Vb_X)
\, \stackrel{\sim}{\to} \, \Vb_X
$$
we arrive at the following pair of maps
of sheaves of Gerstenhaber algebras:
\begin{equation}
\label{chain-Ger}
\begin{array}{ccccc}
\Vb_X & \stackrel{\la}{\leftarrow} &
\bbF_{\Ger}(\Xi_X)
& \stackrel{\si_1}{\hookrightarrow} & \cR(\nCbu(\cO_X))\,.
\end{array}
\end{equation}
In \cite{BLT} it was shown that both $\la$ and $\si_1$
are quasi-isomorphisms of complexes of sheaves.

The quasi-isomorphism $\si_1$ in (\ref{chain-Ger}) provides the sheaf
$\cR_{\ma}$ (\ref{cR-ma}) with a (DG) $\calc$-module structure
over the sheaf $\bbF_{\Ger}(\Xi_X)$\,.
Thus, in order to prove Theorem \ref{main} we need
to show that the sheaf of calculi $(\bbF_{\Ger}(\Xi_X), \cR_{\ma})$
is quasi-isomorphic to the sheaf $(\Vb_X, \OmbX)$\,.

For this purpose we introduce the bar resolution of the
the sheaf $\cR_{\ma}$ of $\cY(\bbF_{\Ger}(\Xi_X))$-modules:
\begin{equation}
\label{BR-a}
\cB\cR_{\ma} = \bigoplus_{k \ge 1} \bs^{1-k} \,
\cY(\bbF_{\Ger}(\Xi_X))^{\otimes \, k} \otimes \cR_{\ma}\,.
\end{equation}

The map $\la$ in (\ref{chain-Ger}) induces the
following map of sheaves of associative algebras
\begin{equation}
\label{cY-la}
\cY(\la) : \cY(\bbF_{\Ger}(\Xi_X)) \to
\cY(\Vb_X)\,.
\end{equation}
Considering the map $\cY(\la)$ on the level of
stalks at a point $x\in X$ we get the map
of associative algebras
\begin{equation}
\label{cY-la-stalk}
\cY(\bbF_{\Ger}(\Xi(A))) \to \cY(\Vb(A))\,,
\end{equation}
where $A$ is the local algebra at the point $x$\,,
$\Xi(A) = \Xi_{\Spec(A)}$ and $\Vb(A)= \Vb_{\Spec(A)}$\,.

Since the variety $X$ is smooth the local algebra $A$
and hence the graded commutative algebra $\Vb(A)$ is regular.
Furthermore, the Gerstenhaber algebra $\bbF_{\Ger}(\Xi(A))$
is a free resolution of $\Vb(A)$\,. Thus, due to
Proposition \ref{Y-alm-exact}, the map (\ref{cY-la-stalk})
a quasi-isomorphism. Hence (\ref{cY-la}) is a quasi-isomorphism
of complexes of sheaves.

Recall that $\cB\cR_{\ma}$ is the free resolution (\ref{BR-a}) of
the sheaf $\cR_{\ma}$ of $\cY(\bbF_{\Ger}(\Xi_X))$-modules. Therefore,
applying the functor
$$
\otimes_{\cY(\bbF_{\Ger}(\Xi_X))}\, \cB\cR_{\ma}
$$
to the quasi-isomorphism (\ref{cY-la}) we get the quasi-isomorphism
of sheaves of $\cY(\bbF_{\Ger}(\Xi_X))$-modules:
\begin{equation}
\label{cY-mod-map}
\cB\cR_{\ma} \, \stackrel{\sim}{\rightarrow} \,
\cY(\Vb_X) \, \otimes_{\cY(\bbF_{\Ger}(\Xi_X))} \, \cB\cR_{\ma}\,.
\end{equation}

Thus we need to show that the sheaf of $\cY(\Vb_X)$-modules
$\cY(\Vb_X) \, \otimes_{\cY(\bbF_{\Ger}(\Xi_X))} \, \cB\cR_{\ma}$
is quasi-isomorphic to $\OmbX$\,.

For this purpose we remark that the sheaf of
(graded) commutative algebras $\OmbX$ can be realized as
a subsheaf of $\cY(\Vb_X)$. Indeed,
\begin{equation}
\label{Omb-cY}
\OmbX = \cY_0(\cO_X) \subset \cY(\Vb_X)\,.
\end{equation}

Using the global section $\bfone \in \G(X, C^{{\rm norm}}_0(\cO_X))$ we
introduce the following global section
\begin{equation}
\label{cycle-E}
E = 1_{\cY(\Vb_X)} \, \otimes \, 1_{\cY(\bbF_{\Ger}(\Xi_X))}\,
\otimes \, \bfone \in  \G(X,
\cY(\Vb_X) \, \otimes_{\cY(\bbF_{\Ger}(\Xi_X))} \, \cB\cR_{\ma} )
\end{equation}
of the sheaf
$\cY(\Vb_X) \, \otimes_{\cY(\bbF_{\Ger}(\Xi_X))} \, \cB\cR_{\ma}$\,.
Here $\bfone \in \G(X, C^{\rm norm}_0(\cO_X))$ is also considered as a
global section of the sheaf
$$
\cR_{\ma} =
\bbF_{\calc}\Big(\bbF_{\bB}(\nCbu(\cO_X), \nCbd(\cO_X)) \Big)_{\ma}\,,
$$
via the unit of the operad $\calc$ and the coaugmentation of
the cooperad $\bB$\,.

It is obvious that $E$ is closed with respect to the total differential
on
$$
\G(X, \cY(\Vb_X) \, \otimes_{\cY(\bbF_{\Ger}(\Xi_X))} \, \cB\cR_{\ma})\,.
$$
Furthermore,
for every point $x\in X$ the
cohomology class of the germ $E_x$ corresponds to the cohomology
class of the germ $\bfone_x$\,.

Using the cycle $E$ and equation (\ref{Omb-cY}) we define the following
map of sheaves
$$
\nu :\OmbX \to
\cY(\Vb_X) \, \otimes_{\cY(\bbF_{\Ger}(\Xi_X))} \, \cB\cR_{\ma}\,,
$$
\begin{equation}
\label{nu}
\nu (a\, d b_1 \, d b_2 \, \dots \, d b_m) =
(i_a\, l_{b_1}\, l_{b_2}\, \dots \, l_{b_m} \,,\, E)\,,
\end{equation}
where $a, b_1, b_2, \dots, b_m $ are local sections
of the structure sheaf $\cO_X$\,.

It is obvious that $\nu$ is compatible with
$\OmbX$-module structures.

We claim that
\begin{pred}
\label{nu-q-iso}
The map
$\nu$ (\ref{nu})
is a quasi-isomorphism of complexes of sheaves.
\end{pred}
{\bf Proof.}
Indeed, let us consider the corresponding map of stalks
$$
\nu_x :(\OmbX)_x \to
\Big(\cY(\Vb_X) \, \otimes_{\cY(\bbF_{\Ger}(\Xi_X))} \, \cB\cR_{\ma}
\Big)_x
$$
at a point $x\in X$\,.
Let $a, b_1, b_2, \dots, b_m$ be germs of functions on $X$
at $x$\,.

Since the cohomology class of the germ $E_x$ corresponds
to the cohomology class of the germ $\bfone_x\in C^{\rm norm}_0(\cO_X)_x$
therefore the cohomology class of the element
$$
(i_a\, l_{b_1}\, l_{b_2}\, \dots \, l_{b_m} \,,\, E_x)
$$
corresponds to the cohomology class of the Hochschild cycle
\begin{equation}
\label{vot-chain}
\sum_{\si\in S_m} (-1)^{|\si|}
(a, b_{\si(1)}, b_{\si(2)}, \dots, b_{\si(m)})
= I_a\, L_{b_1}\, L_{b_2}, \dots\, L_{b_m} \bfone_x
\in C^{{\rm norm}}_m(A)\,.
\end{equation}
Furthermore, under Connes-Hochschild-Kostant-Rosenberg
map (\ref{HKR11}) the chain (\ref{vot-chain}) goes to the form
$a\, d b_1 \, d b_2 \, \dots \, d b_m$\,. Thus $\nu_x$
induces isomorphism on the level of cohomology groups.  $\Box$

Recall that according to Proposition \ref{Y-D-Om}
and equation (\ref{Y-Y-0}) the sheaf of associative algebras
$\cY(\Vb_X)$ is isomorphic to the sheaf
$$
\cD_{\OmbX}[\md]/(\md^2)\,,
$$
where $\cD_{\OmbX}$ is the sheaf of differential operators
on exterior forms and $\md$ is an auxiliary variable of
degree $-1$ which commutes with local sections
of $\cD_{\OmbX}$\,. Let us also recall that due
to the isomorphism (\ref{map-r}) the sheaf
$\cD_X \otimesO \Vb_X[\md]/(\md^2)$ has a natural structure
of left $\cD_{\OmbX}[\md]/(\md^2)$-module.

It is not hard to see that for every sheaf
$\cM$ of $\cD_{\OmbX}[\md]/(\md^2)$-modules
we have the natural isomorphism
of sheaves of $\cO_{X}$-modules
\begin{equation}
\label{cQotimes}
\cD_X \otimesO \Vb_X[\md]/(\md^2) \,\,
\otimes_{\cD_{\OmbX}[\md]/(\md^2)} \,\, \cM
\cong \cO_X \otimes_{\OmbX} \cM\,,
\end{equation}
where the $\OmbX$-module structure on $\cO_X$ is given
by the equation
$$
f\, \eta = \lan \eta,  f \ran\,,
\qquad \eta\in \G(U, \OmbX)\,,
\quad f\in \G(U, \cO_X)\,,
$$
and the pairing $\lan\,,\,\ran$ is defined in
(\ref{pairing}).

Having in mind Proposition \ref{Morita},
we apply the
functor
$$
\cD_X \otimesO \Vb_X[\md]/(\md^2) \,\,
\otimes_{\cD_{\OmbX}[\md]/(\md^2)}
$$
to the map $\nu$ (\ref{nu})
and get the following quasi-isomorphism of
sheaves of $\cO_X$-modules
\begin{equation}
\label{q-iso1}
\tnu \,\, : \,\, \cO_X \,
\stackrel{\sim}{\rightarrow} \,
\cD_X \otimesO \Vb_X[\md]/(\md^2) \,
\otimes_{\cY(\bbF_{\Ger}(\Xi_X))} \, \cB\cR_{\ma}\,,
\end{equation}
where the right $\cY(\bbF_{\Ger}(\Xi_X))$-module
structure on $\cD_X \otimesO \Vb_X[\md]/(\md^2)$ is
obtained from the right $\cY(\Vb_X)$-module structure
via the map $\la$ in (\ref{chain-Ger}).

It is not hard to see that the $\cD_X$-module structure
on
\begin{equation}
\label{sheaf}
\cO_X \cong
\cD_X \otimesO \Vb_X[\md]/(\md^2) \,
\otimes_{\cD_{\OmbX}[\md]/(\md^2)} \, \OmbX
\end{equation}
is the standard one. Furthermore, the element $\md$ acts
on the sections of the sheaf $\cO_X$ in (\ref{sheaf}) by zero simply because
$\md$ has degree $-1$ and $\cO_X$ is concentrated in the single
degree $0$\,.

Let us denote the target of the map $\tnu$ (\ref{q-iso1})
by $\cG$:
\begin{equation}
\label{cG}
\cG =
\cD_X \otimesO \Vb_X[\md]/(\md^2) \,
\otimes_{\cY(\bbF_{\Ger}(\Xi_X))} \, \cB\cR_{\ma}\,.
\end{equation}
We also denote by $\pa^{\cG}$ the total differential
on this sheaf.

For the next proposition we will need the Cech resolution
$\ccCb(\cG)$ of the sheaf of $\cD_X[\md]/(\md^2)$-modules
$\cG$ in the category of sheaves.
\begin{pred}
\label{pravilnoe-pivo}
The map (\ref{q-iso1}) extends to an $A_{\infty}$ quasi-isomorphism
$\Ups$ from the sheaf of  $\cD_X[\md]/(\md^2)$-modules $\cO_X$ to the
the Cech resolution $\ccCb(\cG)$ of the sheaf
$\cG$ (\ref{cG}).
\end{pred}
{\bf Proof.} First, let us prove that the map $\tnu$
(\ref{q-iso1}) is compatible with the action of
the sheaf $\cD_X[\md]/(\md^2)$ up to homotopy
on the level of stalks.

Let $A$ be the stalk $(\cO_X)_x$ of the structure sheaf $\cO_X$
at a point $x$\,. Since $X$ is smooth $A$ is a local regular
algebra.

Let, as above, $x_1, \dots, x_n$ be a regular system of
parameters in $A$\,. The module $\Om^1(A)$ of K\"ahler
differentials is freely generated by the
$1$-forms (\ref{1forms}) and the module of derivations
$\Der(A)$ is freely generated by (\ref{vectors}).

Since the algebra $\cD(A)[\md]/(\md^2)$
is generated by $A$, derivations of $A$ and the element $\md$
it suffices to show that the action of the element
$\md$ and a derivation $v$ of $A$ sends the element
\begin{equation}
\label{el}
1_{\cD(A) \otimesO \Vb(A)[\md]/(\md^2)} \,\,
\otimes_{\cY(\Vb(A))}\,\, E_x
\end{equation}
to cohomologically trivial elements of the
chain complex
\begin{equation}
\label{cG-x}
\cG_x =
\cD(A) \otimesO \Vb(A)[\md]/(\md^2) \,
\otimes_{\cY(\Vb(A))}\, \cY(\Vb(A))
\otimes_{\cY(\bbF_{\Ger}(\Xi(A)))} \, (\cB\cR_{\ma})_x\,.
\end{equation}

The isomorphism
$$
\cY(\Vb(A)) \cong  \Omb(A) \otimes_A \cD(A) \otimes_A  \Vb(A)[\md]/(\md^2)
$$
allows us to consider $\md$ and the derivation
$v\in \Der(A)$ as elements of the algebra $\cY(\Vb(A))$\,.
Thus we need to show that the cocycles
$$
v\, E_x
$$
and
$$
\md\, E_x
$$
are cohomologically trivial.

Since the map  $\tnu$ (\ref{q-iso1}) is a quasi-isomorphism
of complexes of sheaves, the cohomology of the complex
(\ref{cG-x}) is concentrated in the degree $0$\,.
Hence, the cocycle $\md\, E_x$ is a coboundary because
it has degree $-1$\,.

Next, using the generators (\ref{1forms}) and (\ref{vectors})
we rewrite the cocycle $v\, E_x$ as
\begin{equation}
\label{v-Y}
v\, E_x = \frac{1}{n!}\, l_v \prod_{k=1}^n \Big( k -
\sum_{j=1}^n l_{x^j} i_{e_j} \Big) \, E_x\,,
\end{equation}
where the element
$$
\frac{1}{n!}\prod_{k=1}^n \Big( k - \sum_{j=1}^n l_{x^j} i_{e_j}
\Big)\in
\cD(\Omb(A))
$$
operates as a projection on the degree $0$ forms.

Since the cohomology class of $E_x$ corresponds
to the cohomology class of $\bfone_x = 1 \in C^{\rm norm}_0(A)$
the cohomology class of the element (\ref{v-Y}) corresponds
to the class of
$$
\frac{1}{n!}
L_v \prod_{k=1}^n \Big( k - \sum_{j=1}^n L_{x^j} I_{e_j} \Big)\,\, \bfone_x
\in C^{{\rm norm}}_0(A)\,.
$$

It is easy to see that
$$
\frac{1}{n!}
L_v \prod_{k=1}^n \Big( k - \sum_{j=1}^n L_{x^j} I_{e_j} \Big)\,\, \bfone_x =0\,.
$$

Thus the map $\tnu$ (\ref{q-iso1}) is indeed compatible
with the action of the sheaf $\cD_X[\md]/(\md^2)$ up to homotopy
on the level of stalks.

An $A_{\infty}$ morphism from the sheaf of
$\cD_X[\md]/(\md^2)$-modules $\cO_X$ to the sheaf
$\ccCb(\cG)$ is the degree $0$ element
\begin{equation}
\label{Ups}
\Ups \in \bigoplus_{k \ge 0}\,
\bs^k \, \Hom \Big( (\cD_X[\md]/(\md^2)) ^{\otimes \, k} \otimes
\cO_X \,\,, \,\, \ccCb(\cG) \Big)
\end{equation}
satisfying the cocycle condition
\begin{equation}
\label{Ups-eq}
(\pa^{\cG} + \cpa + \cD^{Hoch}) \Ups = 0\,,
\end{equation}
where $\pa^{\cG}$ is the differential on the sheaf (\ref{cG}),
$\cpa$ is the Cech differential and
$\cD^{Hoch}$ is the Hochschild coboundary operator
of the sheaf $\cD_X[\md]/(\md^2)$ with values in the sheaf of
bimodules $\cHom(\cO_x, \ccCb(\cG))$\,.

In other words, $\Ups$ can be defined by the
infinite collection of maps
\begin{equation}
\label{Ups-k}
\Ups_k \in
\Hom \Big( (\cD_X[\md]/(\md^2)) ^{\otimes \, k} \otimes
\cO_X \,\,, \,\, \ccCb(\cG) \Big)\,, \qquad
k =0, 1, 2 \dots
\end{equation}
such that $\Ups_k$ has degree $-k$ and
all the maps satisfy the equations
\begin{equation}
\label{Ups-k-eq}
(\pa^{\cG} + \cpa)\Ups_{k+1}  + \cD^{Hoch}\, \Ups_k = 0\,.
\end{equation}
Our purpose is to show that there is exists an $A_{\infty}$-morphism
$\Ups$ with
$$
\Ups_0 = \tnu\,.
$$

Let us show that there exists $\Ups_1$
satisfying the equation
\begin{equation}
\label{Ups-0-eq}
(\pa^{\cG} + \cpa)\Ups_{1}  + \cD^{Hoch}\, \Ups_0 = 0\,.
\end{equation}
We will find $\Ups_1$ by induction in  degrees of the
Cech complex. In general,
\begin{equation}
\label{Ups-1}
\Ups_1 = \sum_{q=0}^{\infty} \Ups^q_1\,,
\end{equation}
where
\begin{equation}
\label{Ups-q-1}
\Ups^q_1 \in
\Hom \Big( \cD_X[\md]/(\md^2) \otimes
\cO_X \,\,, \,\, \ccC^q(\cG^{-1-q}) \Big)
\end{equation}
and the equation (\ref{Ups-0-eq}) is equivalent to
\begin{equation}
\label{Ups-0-0-eq}
\pa^{\cG} \, \Ups^0_1  + \cD^{Hoch}\, \Ups_0 = 0\,.
\end{equation}
\begin{equation}
\label{Ups-0-q-eq}
\pa^{\cG}\, \Ups^{q+1}_1 + \cpa\,\Ups^q_1 = 0\,,
\qquad
q= 0,1,2, \dots\,.
\end{equation}
If we set $\Ups_0= \tnu$ then
there exists a map $\Ups^0_1$ satisfying equation (\ref{Ups-0-0-eq})
because $\tnu$ is compatible with the action of $\cD_X[\md]/(\md^2)$
up to homotopy on the level of stalks.

Due to equation (\ref{Ups-0-0-eq}) the element
$\cpa\,\Ups^0_1$ is closed with respect to $\pa^{\cG}$
$$
\pa^{\cG}\, (\cpa\,\Ups^0_1) = 0\,.
$$
But
$$
\cpa\,\Ups^0_1 \in
\Hom \Big( \cD_X[\md]/(\md^2) \otimes
\cO_X \,\,, \,\, \ccC^1(\cG^{-1}) \Big)\,.
$$
Thus, using the fact that the cohomology of the stalk
$\cG_x$ is concentrated in degree $0$ we conclude that
there exists the next map $\Ups^1_1$ in (\ref{Ups-q-1})
satisfying the equation
$$
\pa^{\cG}\, \Ups^{1}_1 + \cpa\,\Ups^0_1 = 0\,.
$$
This is the base of the induction.

Let us now assume that for $m > 0$ there exists the collection
of maps $\Ups^q_1$ (\ref{Ups-q-1}) for $q < m$ satisfying
equation (\ref{Ups-0-q-eq}) for $q<m-1$\,. Then, due to equation
(\ref{Ups-0-q-eq}) for $q = m-2$ the map $\cpa\, \Ups^{m-1}_1$
is closed with respect to the differential $\pa^{\cG}$:
$$
\pa^{\cG}\,(\cpa\, \Ups^{m-1}_1) = 0\,.
$$
But
$$
\cpa\,\Ups^{m-1}_1 \in
\Hom \Big( \cD_X[\md]/(\md^2) \otimes
\cO_X \,\,, \,\, \ccC^m(\cG^{-m}) \Big)\,.
$$
Thus, using the fact that the cohomology of the stalk
$\cG_x$ is concentrated in degree $0$ we conclude that
there exists the next map $\Ups^m_1$ in (\ref{Ups-q-1})
satisfying equation (\ref{Ups-0-q-eq}) for $q = m-1$\,.

We proved the existence of the map $\Ups_1$
in (\ref{Ups-k}) satisfying equation (\ref{Ups-k-eq})
for $k=0$\,.

Now we proceed by induction on $k$ in (\ref{Ups-k})
and (\ref{Ups-k-eq}).

Let us assume that $\Ups_k$ (\ref{Ups-k}) are constructed
for $k < m$ and equation (\ref{Ups-k-eq}) holds for $k < m-1$\,.
Then equation (\ref{Ups-k-eq}) for $k = m-2$ implies that
the element
\begin{equation}
\label{on-samyj}
\cD^{Hoch} \Ups_{m-1} \in \Hom\Big( (\cD_X[\md]/(\md^2)) ^{\otimes \, (m-1)} \otimes
\cO_X \,\,, \,\, \ccCb(\cG) \Big)
\end{equation}
is closed with respect to the differential $\pa^{\cG} + \cpa$\,.

Since the sheaf $\ccCb(\cG)$ is acyclic with respect to the
functor of global sections the map $\tnu$ (\ref{q-iso1})
induces the quasi-isomorphism between
the chain complex
\begin{equation}
\label{Rhom-cG}
\Hom\Big( (\cD_X[\md]/(\md^2)) ^{\otimes \, (m-1)} \otimes \cO_X
\,\,, \,\, \ccCb(\cG) \Big)
\end{equation}
and the chain complex
\begin{equation}
\label{Rhom-cO}
\Hom\Big( (\cD_X[\md]/(\md^2)) ^{\otimes \, (m-1)} \otimes \cO_X
\,\,, \,\, \ccCb(\cO_X) \Big)\,.
\end{equation}
It is obvious that the
cohomology of the latter complex is concentrated only
in non-negative degrees.

On the other hand the cocycle (\ref{on-samyj}) has
the negative degree $-m+1$\,. Hence there exists the next
map $\Ups_m$ satisfying equation (\ref{Ups-k-eq}) for $k = m-1$\,.

Proposition \ref{pravilnoe-pivo} is proved. $\Box$

Thus the sheaf (\ref{sheaf}) of
$\cD_X[\md]/(\md^2)$-modules is quasi-isomorphic to the sheaf
$$
\cG = \cD_X \otimesO \Vb_X[\md]/(\md^2) \,
\otimes_{\cY(\bbF_{\Ger}(\Xi_X))} \, \cB\cR_{\ma}\,.
$$

Combining this observation with Proposition \ref{Morita}
we see that the sheaves
$$
\OmbX
$$
and
$$
\cY(\Vb_X) \, \otimes_{\cY(\bbF_{\Ger}(\Xi_X))} \, \cB\cR_{\ma}
$$
are quasi-isomorphic as sheaves of $\cY(\Vb_X)$-modules.

Theorem \ref{main} is proved. $\Box$

\section{Applications and generalizations}
Let, as above, $X$ be a smooth algebraic variety over a field
$\bbK$ of characteristic zero. The homotopy calculus algebra
on the pair $(\nCbu(X),\nCbd(X))$ gives us a $\comm^+$-module
structure on the pair $(\Vb_X, \OmbX)$\,.
Theorem \ref{KS-teo} implies that this $\comm^+$-module
structure on $(\Vb_X, \OmbX)$ is given by the $\wedge$-product
of polyvector fields and contraction of polyvectors with forms.

According to \cite{Swan} and \cite{Y} the Hochschild
cohomology $HH^{\bul}(X)$ of the variety $X$
is the hypercohomology of the sheaf $\nCbu(\cO_X)$:
$$
HH^{\bul}(X) = \bbH^{\bul}(\nCbu(\cO_X))\,.
$$
Furthermore, according to \cite{Cald}, the Hochschild homology
$HH_{\bul}(X)$ of the variety $X$ is the hypercohomology of
the sheaf $\nCbd(\cO_X)$
$$
HH_{\bul}(X) = \bbH^{\bul}(\nCbd(\cO_X))\,.
$$

Thus, using Theorem \ref{main} we get the following generalization
of Corollary 2 from \cite{BLT}
\begin{cor}
For every smooth algebraic variety $X$ over a field
$\bbK$ of characteristic zero the $\comm^+$-algebras
$$
( H^{\bul}(X,  \Vb_X),  H^{\bul}(X, \OmbX))
$$
and
$$
(\,HH^{\bul}(X), HH_{\bul}(X) \,)
$$
are isomorphic. $\Box$
\end{cor}
This statement is the existence part of Caldararu's conjecture
\cite{Cald1} on the Hochschild structure of a smooth algebraic
variety. The cohomological part of this conjecture was proved
in \cite{CV}. As far as we know, D. Calaque, C. Rossi, and
M. Van den Bergh are currently writing an article \cite{CRV} with a proof
of homological part of Caldararu's conjecture.

Combining Theorem \ref{valid-beer} with Theorem \ref{main}
we deduce the statement of cyclic formality conjecture
(see Conjecture 3.3.2 in \cite{Tsygan})
from \cite{Tsygan} for an arbitrary smooth
algebraic variety over a field $\bbK$ of characteristic zero:
\begin{cor}[T. Willwacher, \cite{W}]
\label{Will}
If $X$ is a smooth algebraic variety a field $\bbK$ of characteristic
zero then the sheaf of $\La\Lie^+_{\de}$-algebras
$(\Cbu(\cO_X), \Cbd(\cO_X))$ is formal. $\Box$
\end{cor}
{\bf Remark.} Strictly speaking the methods
used by T. Willwacher in \cite{W} require an additional
assumption $\bbR \subset \bbK$\,. Theorems
\ref{valid-beer} and \ref{main} allow us to remove the
assumption $\bbR \subset \bbK$ from the statement of
Corollary \ref{Will}.

The proof of Theorem \ref{main} can be easily modified for
the following two cases:

\begin{itemize}

\item $X$ is complex manifold with $\cO_X$ being the
sheaf of holomorphic functions,

\item $X$ is a real manifold with $\cO_X$ being
the sheaf of $C^{\infty}$ functions.

\end{itemize}

Thus we get the following obvious modification of
Theorem \ref{main}
\begin{teo}
\label{main1}
If $X$ is a complex manifold (resp. real manifold)
with $\cO_X$ being the sheaf of holomorphic functions
(resp. the sheaf of $C^{\infty}$ real functions) then
the sheaf
$$
(\nCbu(\cO_X), \nCbd(\cO_X))
$$
of homotopy calculi
is quasi-isomorphic to the sheaf
$(\Vb_X, \Om^{-\bul}_X)$ of calculi.
\end{teo}

For $C^{\infty}$ real case we also get the following
statement
\begin{cor}
\label{glob-sect}
If $X$ is a real manifold with $\cO_X$ being the sheaf of
$C^{\infty}$ functions then the homotopy calculus algebra
$$
\Big(\G(X, \nCbu(\cO_X)), \G(X, \nCbd(\cO_X)) \Big)
$$
is quasi-isomorphic to the calculus algebra
$$
\Big(\G(X, \Vb_X), \G(X, \OmbX) \Big)\,.
$$
\end{cor}
{\bf Proof.}
In the $C^{\infty}$ real case the chain of
quasi-isomorphisms connecting the sheaves
$$
(\nCbu(\cO_X), \nCbd(\cO_X))
$$
and
$$
(\Vb_X, \Om^{-\bul}_X)
$$
consists of soft sheaves.
Hence, applying the functor $\G(X, \,\,)$ of global
sections we get the desired result. $\Box$

We would like to mention recent papers \cite{CW}
and \cite{CF}. In paper \cite{CF} A. Cattaneo and G. Felder consider the DG
Lie algebra module $CC^-_{-\bul}(\cO_X)$ of negative cyclic chains
over the DGLA $\Cbu(\cO_X)$ of Hochschild cochains on a $C^{\infty}$ real manifold
equipped with a volume form. Using an interesting modification of the Poisson sigma
model A. Cattaneo and G. Felder construct a curious $L_{\infty}$ morphism
(not a quasi-isomorphism!)
from this DG Lie algebra module to a DG Lie algebra module modeled on polyvector
fields using the volume form.  A. Cattaneo and G. Felder also apply
this result to a construction of a specific
trace on the deformation quantization algebra of a unimodular Poisson manifold.
Although this trace can be constructed using the formality quasi-isomorphism for
Hochschild chains \cite{Sh}, \cite{W} the relation of the $L_{\infty}$
morphism of A. Cattaneo and G. Felder to the formality quasi-isomorphism is
a mystery.

Paper \cite{CW} is devoted to the proof of Kontsevich's cyclic formality
conjecture for cochains formulated in paper \cite{Sh1}. We suspect that
the statement of this conjecture may be related to Theorem \ref{main}
and Corollary \ref{Will} via the Van den Bergh duality \cite{VB}
between Hochschild cohomology and Hochschild homology.

~\\

\noindent\textsc{Department of Mathematics,
University of California at Riverside, \\
900 Big Springs Drive,\\
Riverside, CA 92521, USA \\
\emph{E-mail address:} {\bf vald@math.ucr.edu}}

~\\

\noindent\textsc{Mathematics Department,
Northwestern University, \\
2033 Sheridan Rd.,\\
Evanston, IL 60208, USA \\
\emph{E-mail addresses:} {\bf tamarkin@math.northwestern.edu},
{\bf tsygan@math.northwestern.edu}}

\end{document}